\input amstex
\documentstyle{amsppt}

\hsize=4.75in
\vsize=8in
\nologo

\def\J{{\Cal J}}
\def\G{{\Cal G}}
\def\H{{\Cal H}}
\def\I{{\Cal I}}

\def\O{{\Cal O}}

\def\K{{\Cal K}}

\def\T{{\Cal T}}

\def\F{{\Cal F}}
\def\L{{\Cal L}}

\def\A{{\Cal A}}
\def\B{{\Cal B}}
\def\C{{\Cal C}}

\def\bN{{\Bbb N}}

\def\ov{\overline}

\def\comp{{\Bbb {C}}}

\def\cst{{C${}^*$}}

\rightheadtext {Hilbert $C^*$--bimodules}
\leftheadtext {T. Kajiwara, C. Pinzari, Y. Watatani}
\topmatter
\title
Jones index theory for Hilbert $C^*$--bimodules and its equivalence
with conjugation theory
\endtitle
\author
Tsuyoshi Kajiwara\footnote{Partially supported by
Grants-in-Aid
for Scientific Research 12640210 and 14340050 from the Ministry of
Education, Culture, Sports, Science and Technology of Japan.}
\\ Department of Environmental and Mathematical Sciences\\ 
Okayama University, Tsushima, 700-8530, Japan\\ \\ Claudia Pinzari
 \\
Dipartimento di Matematica\\
 Universit\`a di Roma la Sapienza,
P.le A. Moro, 2, 00185 Roma, Italy\\ \\ Yasuo Watatani${}^1$\\ 
Department of Mathematical Sciences\\
 Kyushu University, Hakozaki, 
812-8581 Japan
\endauthor
\affil
\endaffil
\endtopmatter
\document

\heading
{\bf Abstract}
\endheading

{We introduce the notion of finite right (respectively left) 
 numerical  index
on a $C^*$--bimodule ${}_AX_B$  with a bi-Hilbertian structure. 
This notion 
is based 
on a Pimsner--Popa-type inequality.
The right (respectively left)
index element of $X$ can be constructed in the centre of the
enveloping von Neumann algebra
of $A$ (respectively $B$). $X$ is called of finite right 
index if the right index element lies in the multiplier algebra of $A$.
 In this case we can perform the Jones basic construction.
Furthermore the $C^*$--algebra of bimodule mappings with a right adjoint
is a continuous field of  $C^*$--algebras over a 
compact Hausdorff space, whose fiber dimensions are bounded above
by the
index.
We show that if $A$ is unital, the right index element belongs to $A$ if and
only if $X$ is finitely generated as a right module.
A finite index bimodule is a bi-Hilbertian $C^*$--bimodule
which is at the same time of finite right and left index.

We study the relationship between the notion of finite index and the 
notion of conjugation in a 
 tensor 2-$C^*$--category, in the sense of Longo and Roberts. 
We show that bi-Hilbertian, finite index $C^*$--bimodules, when regarded
as objects
of the tensor 2-$C^*$--category of right Hilbertian
  $C^*$--bimodules, are precisely those objects with a conjugate in that
category. }
\bigskip

\heading
Introduction
\endheading

The theory of conjugation 
in abstract  tensor $C^*$--categories
appeared 
in the algebraic    formulation of Quantum Field
Theory
\cite{H}.

In this connection,  Doplicher and Roberts showed in \cite{DR1}, \cite{DR2}
that any    symmetric tensor $C^*$--category with conjugation can be embedded
into a category of finite dimensional Hilbert spaces, and therefore the category
is isomorphic to the representation category of a compact group. 

However, some tensor $C^*$--categories with a unitary braiding, arising
from low dimensional QFT, can not be embedded
into categories of Hilbert spaces \cite{LR}.

R. Longo and J.E. Roberts  recently
studied in \cite{LR}  conjugation  in 
general  
tensor $C^*$--categories, and they showed that this notion is closely
 related to the  Jones index theory
 for subfactors \cite{J}.

A different, but related, approach to 
conjugation (or duality) in a tensor 
$C^*$--category has been recently studied  by   Yamagami  in \cite{Y}.

 An interesting open problem is to
decide    which tensor $C^*$--categories with conjugation can be 
embedded into categories of Hilbert $C^*$--bimodules.
A related, easier, problem is  to ask   which sort of
bimodules
should appear.

In this paper we investigate the latter problem. 
We introduce
 Jones index theory for general
Hilbert bimodules over pairs of  $C^*$--algebras, and we show that
bimodules of
finite index are precisely those endowed with a conjugate object in the
same category.

Therefore if a tensor $C^*$--category with conjugation can be embedded into a 
category of right Hilbert  
$C^*$--bimodules, these should be of finite Jones index.

In \cite{KW1}  the first and third--named 
authors studied Hilbert $C^*$--bimodules with finite Jones index in the
 case where the $C^*$--algebras are unital and the
bimodules are finitely generated as right as well as left
modules. It turns out that a bimodule of finite index over {\it unital}
$C^*$--algebras, in
the 
sense introduced in this paper (Def. 2.23), is automatically finitely
generated as a bimodule, and therefore it is of the kind described in
\cite{KW1} (Cor. 2.25).

If $A$ and $B$ are 
$C^*$--algebras, an object of our category  is a right 
Hilbert $B$--module $X$ with an action 
of $A$ on the left given by a nondegenerate $^*$--homomorphism of $A$ into
the 
$C^*$--algebra of $B$--module maps of $X$ into itself with adjoint.
We will refer to such a bimodule   as a 
  {\it right\/} Hilbert $A$--$B$ 
$C^*$--bimodule.
The space of intertwiners from ${}_AX_B$ to ${}_AY_B$ is the set of
adjointable bimodule maps.

Strictly speaking, the 
notion of tensor $C^*$--category is not correct for the category of right
Hilbert 
$C^*$--bimodules.  In fact, the set of objects is not a unital semigroup. 
We  rather need the framework of     tensor 
2-$C^*$--categories, where a theory of conjugation can be developed without any
substantial difficulty, as sketched in the appendix of
\cite{LR}.

On the other hand, the notion of Jones index 
leads us   to introduce bi-Hilbertian  structures on
$C^*$--bimodules. ${}_AX_B$  is called bi-Hilbertian if it
is at the same
time a right and a left Hilbert $C^*$-bimodule in such a way that the two
Banach space norms arising from the  two inner products are equivalent.
Examples are Rieffel's imprimitivity bimodules \cite{R1}. 

A bi-Hilbertian $C^*$--bimodule will be called of {\it finite right
(or left) numerical
index} if a suitable Pimsner--Popa-type inequality relating the two
Banach space
norms holds (see Definitions 2.8 and 2.9). A bimodule of finite numerical
index is a bimodule which is at the same time of finite right and left
numerical index.

If $X$ is of finite right (or left) numerical index, we  construct the
right
(or left) index element of $X$ as a positive
central element of $A''$ (or $B''$).

 While the tensor product 
of two imprimitivity bimodules is still an imprimitivity bimodule,
 a tensor product of bi-Hilbertian bimodules can not
be made, in general, into a bi-Hilbertian bimodule in the natural way. 
However,
if ${}_AX_B$ has finite right numerical
index and
${}_BY_C$ has
finite left
numerical  index  then the algebraic tensor product bimodule $X\odot_B
Y$ can be
completed
into a bi-Hilbertian bimodule in the natural way (Prop. 2.13). 
One can not
assume less: the left and right seminorms on
$X\odot_B\ell^2(B)_{\K\otimes B}$  
 are equivalent if and only if $X$ is of finite right numerical 
index.

Typical examples of bimodules of finite right numerical index arise, of
course, from conditional expectations between $C^*$--algebras satisfying a 
Pimsner--Popa inequality.
The work of Frank and Kirchberg \cite{FK}   shows that under this only
assumption 
the index element of the conditional expectation lies in the enveloping
von Neumann algebra of the bigger
algebra. This reflects  the fact that
 the Jones
basic construction is not always possible in the $C^*$--algebraic setting.
 
We introduce in our theory an extra requirement:  
the right index element of ${}_AX_B$ should lie in  the multiplier algebra
of
$A$,
and therefore in its centre.
When this assumption is satisfied, we say that $X$ is of  finite right
index.

We prove that this property is in fact equivalent to
other properties which would seem stronger a priori, such as, e.g., 
the fact that
the left action
of $A$ on $X$ has range into the compacts $\K(X_B)$ (Theorem 2.22).

Finite bases are a useful tool in  Jones index theory, 
as they lead to a simple formula for the index element.
One of the  problems in Jones index theory for $C^*$--algebras is to
be able to establish
existence of finite bases.
In this direction, Izumi proved in \cite{I} that a Pimsner--Popa
conditional expectation from a simple, unital $C^*$--algebra admits 
a finite quasi--basis in the sense of \cite{W}. 

We obtain,
as a consequence of  Theorem 2.22, the following more general result:
 if $X$ is a bi-Hilbertian $A$-$B$
bimodule of finite right numerical
index
and if $A$ is unital, 
  the right index element of $X$
 belongs $A$ if and only  
if $X$ is finitely generated as a right $B$-module.

More generally, we show that bimodules with finite right index over
$\sigma$--unital
$C^*$--algebras admit countable, unconditionally convergent, bases
(see Prop. 1.6 and Cor. 2.24). Thus
the the right index 
element of $X$
is the strict limit
of the sum of left inner products with themeselves
of elements of any  countable right basis.

In the general case we shall deal with {\it generalized bases} in the
sense
of Def. 1.8, which always exist, as shown in Prop. 1.9.

Now the assumption that the right index element be a multiplier of $A$
 guaranties 
   the existence of the
 Jones basic construction (see Theorem 2.30), which takes the form of
a positive, $A$--bilinear, strictly continuous map $F:\K(X_B)\to A$
satisfying a Pimsner--Popa inequality.
Since left $A$-action lies in $\K(X_B)$, $F$ extends uniquely to a
$A$-bilinear map
$\hat{F}:\L(X_B)\to M(A)$ between the corresponding multiplier algebras.
The right index element of $X$  coincides with
 $\hat{F}(I)$ and it
can be reached by the strict
limit of the image under $F$ of an
approximate unit of $\K(X_B)$.

We illustrate our approach to index theory with a
typical example of an 
inclusion of commutative unital
$C^*$--algebras satisfying a Pimsner--Popa inequality,  for which a 
finite quasi--basis in the sense of \cite{W} does not exist. This class of
examples
arises  from branched coverings, or orbifolds. It was first pointed
out in \cite{W}, 2.8 and later analyzed by Frank and Kirchberg in
\cite{FK}.
 We show that this inclusion is in fact determined by
  a canonical nonunital subinclusion  
of finite right index in our sense (cf. Example 2.35).

Let us go back to our aim of comparing Jones index
theory for Hilbert bimodules with conjugation theory.
  One of the main result of this paper is  that these two approaches
are equivalent  
(cf. Theorems 4.4 and 4.14).
We show that a bi-Hilbertian $C^*$--bimodule
has finite Jones index in our sense if and only if it  has a 
conjugate object in the 
2-$C^*$--category of right Hilbert $C^*$--bimodules with nondegenerate
left actions.
A choice of the operators $R$ and $\overline{R}$ satisfying the
conjugate equations
leads to a specification of a left $A$-valued inner product on $X$ making
it into a finite index bi-Hilbertian bimodule.
Conversely, a finite index bi-Hilbertian structure on $X$ yields a
solution of the conjugate equations. 
 The square of the
 Longo-Roberts dimension relative to a solution ($R$,
$\overline{R}$) coincides with the corresponding
numerical index of the bimodule.

Imprimitivity bimodules are finite index bimodules, with left and right 
index elements equal to the identities.
They
 can be characterized, among
general right Hilbert $C^*$--bimodules
 as those objects with trivial 
minimal dimension
  (Cor. 4.16).

We show two applications of
our characterization theorem. The first one is that 
if ${}_AX_B$ is of finite index, the set
of Hilbert module mappings (i.e. those with an adjoint)
on $X_B$ commuting with the left action coincides,
as an algebra, with the set  of Hilbert module mappings on ${}_AX$ commuting with
the right action (Cor. 4.6). Moreover, each one of these $C^*$--algebras
is a 
 continuous bundle of finite dimensional $C^*$--algebras over a
compact space, in the sense of \cite{KW} (Theorem 3.3).

As a second application, we show that the
tensor product of two bi-Hilbertian bimodules
of finite  (resp. numerical) index is still of finite  (resp. numerical)
index
(Theorem 5.2).

An index theory for Hilbert bimodules
turns out to be  more general than for conditional
 expectations in the case where the
algebras are not $\sigma$-unital.
In fact, it is known that if $E:B\to A$ is a conditional expectation 
satisfying a Pimsner-Popa inequality,   an approximate unit 
of $A$ must be an approximate unit of $B$ as well. Therefore if $A$
is $\sigma$-unital, $B$ must be $\sigma$-unital as well.
In particular, a unital $C^*$--algebra and a non-$\sigma$-unital one 
can not be linked by a Pimsner-Popa conditional expectation.  
However, a $II_1$ factor can be strongly Morita
equivalent to a non-$\sigma$-unital $C^*$--algebra (see \cite{BGR}).

In Section 6 we will discuss further examples of bimodules of finite index 
arising from locally finite directed graphs and topological
correspondences.

As Franks  pointed out to the 
third--named  author,  our notion of basis is  not  standard  terminology
in the theory of Banach spaces. 
Frames are replacement for 
bases in Hilbert spaces and  naturally arise in wavelet theory.  
Franks and Larson  have recently studied 
the concept of module frames for Hilbert $C^*$--modules in  \cite{FL1} and  
\cite{FL2}.  Module (quasi-)bases or frames are useful in 
Jones index theory \cite{BDH}, \cite{FK}, \cite{W}. 
Our concept of basis corresponds to
their  the notion of  {\it standard normalized tight frame\/}.
Since we will not consider orthogonal module bases in full generality,
we will just use the term {\it basis\/}, for simplicity.

 This paper is an extended version of an appendix contained in the 
draft of  \cite{KPW1}.

\heading
1. Countable bases and generalized bases
\endheading

Let $A$ be a \cst-algebra and $X = X_A$  a right  Hilbert
\cst--module over $A$. We denote by  
$\L(X_A)$ the $C^*$--algebra of   $A$--module maps on $X$ with an
adjoint.

A finite subset $\{ u_i\}_i$  of $X$ is called a {\it finite basis} if 
    $x = \sum _i u_i(u_i|x)_A $ for  $x \in X$. Our  
aim in this section is to generalize this notion  to comprehend {\it countable 
bases} or, more generally, {\it generalized bases} in a sense that will be
explained (see Definitions 1.1 and 1.8).
These are infinite bases, and they  will be a good substitute of
finite bases in the case where  the finite generation property does not
hold.
Indeed, we will show, generalizing slightly  Dixmier's proof
of existence of approximate units in a $C^*$--algebra \cite{Di},
 that a right Hilbert module $X$ always admits a generalized basis, and,
in the
case where $X$ is countably generated, it actually admits a countable
basis.
 
 We denote by 
$\theta^r_{x,y}$ the {\it rank one operator\/} on $X$ defined by 
$\theta^r_{x,y}(z) = x(y|z)_A$.  The  linear span of rank one
operators is denoted by $FR(X_A)$ and called
the ideal of {\it finite rank operators}. Its norm closure, $\K(X_A)$, is the
 $C^*$--algebra of 
{\it compact operators\/}, which is a closed ideal in $\L(X_A)$.  

For a left  Hilbert $A$--module $X$, we define 
the rank one operators  by 
$\theta^l_{x,y}(z)={}_A(z|y)x$, and the spaces of finite rank operators 
$FR({}_AX)$,
compact operators $\K({}_AX)$ and
  adjointable left 
 $A$--module maps  $\L({}_AX)$ are defined similarly.

The right  Hilbert module $X_A$ has a finite basis if and only if  
$\L(X_A) = 
\K(X_A)$. If in addition $A$ is unital, $\L(X_A)=\K(X_A)$ 
if and only if $X_A$ is finite projective as a right module.

We will often make use of  the following formula, derived for example 
in Lemma 2.1 of \cite{KPW1}. If
$X_A$ is a right 
Hilbert $A$--module, 
for any $x_1,\dots,x_n$, $y_1,\dots,y_n\in X_A$,
$$
 \| \sum_{i=1}^n \theta^r_{x_i,y_i}\| 
 = \|((x_i|x_j)_A)_{ij}^{1/2}((y_i|y_j)_A)_{ij}^{1/2} \|,
$$ 
where the norm
  at the right hand side is evaluated in the matrix $C^*$--algebra
$M_n(A)$.
\medskip

\noindent{\bf 1.1 Definition} 
Let $X$ be a right Hilbert $A$--module.
 We say that  a subset $\{u_i\}_{i \in \Lambda} \subset X$ is an ({\it 
unconditionally convergent} {\it right}) {\it basis} for $X$ 
 if for any $x \in X$ and for any  $\varepsilon >0$ 
 there
exists a finite subset $F(\varepsilon, x)$ of $\Lambda$ such that for 
every finite
 subset $F'$ with $F' \supset F(\varepsilon, x)$,   
$
 \| x - \sum_{i\in F'} u_i(u_i|x)_A \| < \varepsilon.$

In other words the net $F\in\{\text{finite subsets of }\Lambda\} 
\to\sum_{i\in F}u_i(u_i|x)_A$ should converge  $x$, for all $x\in X$.

One can easily show that if $\{u_i\}_{\Lambda}$ is an unconditionally 
convergent right basis then
$$\|u_i\| \le 1,\quad i\in{\Lambda}$$ and 
$$\|\sum_{i\in F} u_i(u_i|x)_A \| \le \|x\|$$ for any finite subset $F$ of
$\Lambda$.  In 
fact, 
for $x \in X$, 
$$
 0\leq  (x|\sum_{i\in F}  u_i(u_i|x)_A)_A\le(x|\sum_{i\in F\cup
F(\varepsilon, 
x)} u_i(u_i|x)_A)_A$$
$$   \le (x|x)_A+\varepsilon\|x\|,
$$ 
therefore
$$
 0 \leq (x|\sum_{i\in F}  u_i(u_i|x)_A)_A   \le (x|x)_A,
$$ 
and this shows that $0 \le \sum_{i\in F}\theta^r_{u_i,u_i} \le 1$, which 
implies 
 $$(x|(\sum_{i\in F}\theta^r_{u_i,u_i})^2(x))_A \le (x 
| x)_A. $$ 
Hence we have  $\|\sum_{i\in F} u_i(u_i|x)_A \| \le \|x\|$. Moreover 
$$\|u_i\| = \| \theta^r_{u_i,u_i} \|^{1/2} \le 1.$$  
For any  unconditionally convergent 
right 
 basis $\{u_i\}_{i \in \Lambda} \subset X$,  the net 
 $F\to
\sum_{i\in F}\theta^r_{u_i,u_i}
$   
is an approximate unit for $\K(X_A)$.  In fact, for any $x,y \in X$, 
$$
\|(\sum_{i\in F}\theta^r_{u_i,u_i})\theta^r_{x,y} - \theta^r_{x,y}\| 
= \|\theta^r_{\sum_{i\in F}  u_i(u_i|x)_A - x, y} \| 
\leq \|\sum_{i\in F}  u_i(u_i|x)_A - x\|\| y\| \rightarrow 0.
$$ 
This shows the claim as $\|\sum_{i\in F}\theta^r_{u_i,u_i}\|\leq1$. 
For any fixed $n \in {\Bbb N}$, the $n \times n$ 
operator matrix $((u_i|u_j)_A)_{i,j}$ is a positive contraction as
$ \|((u_i|u_j)_A)_{i,j}\| = \| \sum _{i=1}^n \theta^r_{u_i,u_i} \| $.
\par
In this paper an unconditionally convergent right (or left) basis will be 
simply called a  {\it right (or left) basis}.
\par
A sequence $\{u_i\}_{i \in \bN} \subset X$ is called a {\it weak basis}
 for  $X$ if 
for any $x \in X$ we have $x = \lim_{n \to \infty}\sum_{i=1}^n 
u_i(u_i|x)_A$. 
One can similarly show that  $\|u_i\| \le 1$ and 
$\|\sum_{i=1}^n u_i(u_i|x)_A \| \le \|x\|$ for all $n\in{\bN}$.
 A countable subset $\{u_i\}_{i \in \Lambda} \subset X$ will also be called 
 a weak basis if it is a weak basis with respect to some
bijective correspondence identifying $\Lambda$ with ${\bN}$.

One can easily  show, using Kasparov's stabilization trick, that
 any countably generated Hilbert $C^*$--module $X$ over a 
 $\sigma$--unital \cst--algebra $A$ admits   
  a weak 
basis.  We shall show that that $X$ actually admits  
an unconditionally convergent countable basis.  
\medskip

\noindent{\bf 1.2 Lemma} {\sl 
Let $X$ be a right Hilbert  $C^*$--module over $A$, and 
 $\{u_i\}_{i \in \Lambda}$  a  subset of $X$.  
 If there exists a  subset  $X_0$
 in $X$ such that $X_0A$ is total in $X$ and such that for every $x \in 
X_0$, 
$\sum_{i \in \Lambda}u_i(u_i|x)_A =x$
  then $\{u_i\}_{i \in \Lambda}$ is
 a right basis for  $X$.  
Similarly, if $\Lambda={\bN}$ and if 
$\sum_{i=1}^\infty u_i(u_i|x)_A =x$ for $x\in X_0$, then
$\{u_i\}_{i\in{\bN}}$
is a weak basis.}\medskip

\noindent{\it Proof} By right linearity of the inner product, 
$x=\sum_{i\in \Lambda} 
u_i(u_i|x)_A$ for all $x$ in the linear span of $X_0A$, which is dense in 
$X$. 
Let $F$ be a finite subset of $\Lambda$. We then have 
$$(x|\sum_{i\in F}\theta_{u_i,u_i}(x))_A \le (x|x)_A$$ for $x$ in the 
linear span of $X_0A$,
therefore
this inequality holds for all $x \in X$. This shows that
 $\|\sum_{i \in F}\theta_{u_i,u_i}\| \le 1$. 
 A 3$\varepsilon$ argument will conclude the proof.\medskip

\noindent{\it Remark} Any countable unconditionally 
convergent basis is a 
weak basis. 
We shall show that the converse is true under
   suitable
conditions.
\medskip

\noindent{\bf 1.3 Definition} 
 We say that a weak basis 
 $\{u_i\}_{i\in \Lambda}$ for  $X$ satisfies the {\it finite intersection 
property} 
 if for every $i$
 there exists a finite subset $G_i \subset \Lambda $  such that for 
 any $j \in G_i^c$  we have $(u_i|u_j)_A = 0$.  
\medskip

\noindent{\bf 1.4 Lemma} {\sl
 Let  $\{u_i\}_{i \in \Lambda} \subset X$  
 be a countable  weak basis for  $X$.  Suppose that  either
  \roster
\item $\{u_i\}_{i \in \Lambda}$ satisfies the finite intersection
property, or
\item the set $\{\theta_{u_i,u_i}, i \in \Lambda\}$ is commutative.

\endroster 
Then
       $\{u_i\}_{i\in \Lambda}$ is an unconditionally convergent  basis.}\medskip

\noindent{\it Proof} (1) Suppose that  $\{u_i\}_i$ satisfies the finite intersection 
property. Then for every $i$ there exists a finite subset
$G_i\subset\Lambda$ such
that $u_i=\sum_{j\in G_i} u_j(u_j|u_i)_A$, therefore $u_i=\sum_{j\in\Lambda}u_j(u_j|u_i)_A$.
We can now appeal to Lemma 1.2 with $X_0=\{u_i, i\in\Lambda\}$.
\par
\noindent
(2) Suppose now that the rank one operators $\theta_{u_i,u_i}, {i \in \Lambda}$ commute pairwise.
Let us fix an identification of 
$\Lambda$ with  $\Bbb N$.  For any $x \in X$ and for any 
$\varepsilon > 0$
 there exists $N$ such that for any $n\ge N$, 
$$
 \| x - \sum_{k=1}^n u_k(u_k|x)_A \| < \varepsilon.
$$		 
Setting $F_0=\{1,2,\cdots, N\}$, for any finite subset $F$  of
 $\Lambda $ with $F \supset F_0$,  we have 
$$
0 \le  I - \sum_{k \in F}\theta_{u_k,u_k} \le I - \sum_{k \in 
F_0}\theta_{u_k,u_k}. 
$$
Since $\{\theta_{u_k,u_k}\}$ commute each other, we have
$$ (I - \sum_{k \in F}\theta_{u_k,u_k})^2 
 \le (I-\sum_{k \in F_0}\theta_{u_k,u_k})^2 .$$  
Thus we have 
 $ (x|(I - \sum_{k \in F}\theta_{u_k,u_k} )^2 x)_A 
 \le (x|(I - \sum_{k \in F_0}\theta_{u_k,u_k} )^2 x)_A 
$.
This implies that  
$$
 \| x - \sum_{k \in F}u_k(u_k|x)_A \| 
\le \| x - \sum_{k \in F_0}u_k(u_k|x)_A \| < \varepsilon.
$$
Our next aim is to show that countable unconditionally convergent bases
exist under some countability generation property.
\medskip

\noindent{\bf 1.5 Lemma} {\sl Let $A$ be a $\sigma$--unital \cst--algebra.
 Then there exists a sequence $\{u_j\}_{j}$ of positive elements 
of $A$ such that for every $n\in{\Bbb N}$,
 $v_n:=\sum_{j=1}^n(u_j)^2$ is a contraction, the sequence $\{v_n\}_n$ is a
 countable approximate unit of $A$,  and  $u_mu_n = 0$ 
 for any pair of positive integers $m$ and $n$ with $|m - n| \geq 2$.}\medskip  
   
\noindent{\it Proof} Since $A$ is $\sigma$--unital, there exists a 
strictly
 positive contraction $h$ in $A$.     For each positive integer $n$
we define a positive continuous function $f_n(x) \in C_0([0,1])$ as
follows: 
For $n \ge 2$, let

$$f_n(x) =
\cases
   0 &\quad   0 \le x < \frac{1}{2^n},\\ 
   2^nx -1 &\quad  \frac{1}{2^n} \le x < \frac{1}{2^{n-1}},\\ 
   -2^{n-1}x + 2 &\quad  \frac{1}{2^{n-1}} \le x < \frac{1}{2^{n-2}},\\
   0 &\quad    \frac{1}{2^{n-2}} \le x \le 1.
\endcases
$$
 We set 
$f_1(x) = 0$ if $0 \le x \le \frac{1}{2}$ and  
$ \ f_1(x) = 2x -1$ if $\frac{1}{2} \le x \le 1$.   Then we have 
$f_nf_m = 0$ if $|n - m | \ge 2$, $\ f_n(0) = 0$ and  
$\sum_{n=1}^{\infty} f_n(x) = 1$ for $x \in (0,1] \ $.    
Define $\ u_n= f_n(h)^{1/2}$, so that  $v_n = \sum_{i=1}^n f_j(h)$. 
Then $u_nu_m = 0$  if  $|n - m | \ge 2$. Moreover $\{v_n\}_n$ is a 
countable approximate unit of $A$ since the norm
 closure of $Ah$ and $hA$ is equal to $A$ (see, e.g., 
 Prop. 12.3.1 in \cite{B}). \medskip

We are now ready to prove the existence of  unconditionally convergent
bases.\medskip

\noindent{\bf 1.6 Proposition} {\sl  Let $X$ be a countably generated 
right
 Hilbert  \cst--module over a $\sigma$--unital \cst--algebra $A$.
 Then $X$ has an  unconditionally convergent right basis indexed by 
${\bN}=\{1,2,\dots\}$.}\medskip

\noindent{\it Proof}
As a first step, we consider the case where $X = A$ with the inner 
product
$(x|y)_A = x^*y$ for $x, y \in A$.  
 We choose a sequence $\{u_n\}_n$  of $A$  as in Lemma 1.5.  
Then $ \sum_n u_n(u_n|x)_A = \sum_n (u_n)^2 x = \lim_n v_n x = x$. 
Thus $\{u_n\}_{n \in \bN}$ is a weak basis with the finite intersection 
property.  
Hence $\{u_n\}_{n \in \bN}$ is an unconditionally convergent right basis
for $X = A$, by
Lemma 1.4.  
\par 
Next we consider the case where $X=\ell^2(A): = \ell^2(\bN)\otimes_{\Bbb
C} A$. 
Let $\{e_n\}_{n \in \bN}$ be the canonical orthonormal basis for
$\ell^2(\bN)$ and 
 let $\varphi : \bN^2 \rightarrow \bN$ be a bijection. 
  Put $g_{(i,j)} = e_i \otimes u_j \in \ell^2(\bN)\otimes A$.  
  Define $w_n = g_{\varphi^{-1}(n)} \in X$.      
Let $X_0 = \{(x_1,x_2,\cdots) \in X \ ; \ x_i=0 \text{ except for 
finitely many } 
\  i \}$.   Then  $X_0$ is dense in $X$,  and for any $x \in X_0$ we 
have 
$ \| x - \sum_{k=1}^{n}w_k(w_k|x)_A \| \to 0$ as $n \to \infty$.
By Lemma 1.2, this holds for every $x \in X = \ell^2(A)$.  
Thus $\{w_n\}_{n \in \bN}$ is a weak basis with finite intersection 
property and therefore a basis by
Lemma 1.4.
\par
In the general case   $X \simeq p\ell^2(A)$ for some projection $p$ in 
$\L(\ell^2(A))$ by 
 Kasparov's stabilization theorem (see, e.g., \cite{B}).  
We set $s_i = pw_i \in X$.  For $x \in X \subseteq \ell^2(A)$, we have 
$$
  x = \lim_{F \subset \Lambda \text{:finite}} 
\sum_{i \in F}w_i(w_i|x)_A,
$$
thus 
$$
 x = px =  \lim_{F \subset \Lambda \text{:finite}} \sum_{i \in 
F}pw_i(w_i|px)_A
   = \lim_{F \subset \Lambda \text{:finite}} \sum_{i \in F}s_i(s_i|x)_A  
$$
Therefore $\{s_i\}_{i \in \bN}$ is an unconditionally convergent
right basis for $X$.  
\medskip

We do not know whether the tensor product of two bases is still a basis.  
But for some  bases this is true, as the following proposition shows.
\medskip
 
\noindent{\bf 1.7 Proposition} {\sl Let $A$ and $B$ be  $\sigma$--unital 
\cst--algebras.  
Let $X$ and $Y$ be   countably generated right Hilbert  \cst--modules over 
$A$ and $B$ respectively.
Let
 $\phi : A \rightarrow \L(Y_B)$ be a $^*$--homomorphism.  Then there exist 
bases 
$\{s_i\}_i$ for $X$ and $\{t_j\}_j$ for $Y$ such that
  $\{s_i\otimes t_j\}_{i,j}$ is a basis  for the right Hilbert 
$B$--module 
$X\otimes_A Y$. }\medskip  

\noindent{\it Proof}  We first suppose that $X=\ell^2(A)$.  Then, as in
the proof of Prop. 1.6, we have  a basis
 $\{w_i\}_{i \in \Lambda_1}$  for $\ell^2(A)$ with the finite intersection 
property,
   where $\Lambda_1 = \bN$.  
 We similarly have  a basis  $\{t_j\}_{j \in \Lambda_2}$ for $Y$.  
\par
Set $ Z_0 = \{z \in X\otimes_A Y ; z=\sum_{i=1}^m w_ia_i \otimes y_i, 
a_i \in A, y_i \in Y, m= 1,2,... \}$.  
Then $Z_0$ is  dense in $X \otimes_A Y$.  Using the finite intersection 
property 
one can show  that for any  $z=\sum_{i=1}^m w_ia_i \otimes y_i\in Z_0$ 
and 
any  $\varepsilon >0$ there exists a finite
 subset $F_0$ of $\Lambda_1\times \Lambda_2$ such that for every finite
 subset $F$ containing $F_0$, we have 
$$
 \| z - \sum_{(k,j) \in F}w_k\otimes t_j(w_k \otimes t_j | z)_B \| < 
\varepsilon. 
$$
In fact, for any $i=1,...,m$ there exists a finite subset 
$G_i \subset \Lambda_1 $  such that for 
 any $j \in G_i^c$  we have $(w_i|w_j)_A = 0$.   
Put $H_1 = \cup_{i=1}^m G_i$ .  Let $n$ be the cardinality of $H_1$. 
 For $k \in G_i$ there 
exists a finite set $L_{i,k}  \subset \Lambda_2$ 
such that for any finite subset 
$L$ of $\Lambda_2$ with  $L_{i,k} \subset L$ we have 
$$
\| (w_k|w_i)_Aa_iy_i - \sum _{j\in L}t_j(t_j|(w_k|w_i)_Aa_iy_i)_B \| 
< \varepsilon/n .
$$
Put $H_2 = \cup_{i=1}^m \cup_{k \in G_i} L_{i,k}$.  Then it suffices to
choose
 $F_0 = H_1\times H_2$
By Lemma 1.2
 $\{w_j \otimes t_k\}_{(j,k) \in \Lambda_1 \times \Lambda_2}$ is basis
 for $X\otimes_A Y$.   
In the general case  there exists a projection $p \in \L(\ell^2(A)_A)$
such 
that 
$X \simeq p\ell^2(A)$.  Define $s_i = pw_i \in X$.  Then
$\{s_i\otimes t_j\}_{i,j}$ is a basis for 
 $X\otimes_A Y \simeq p\ell^2(A) \otimes_A Y$.  
\medskip

In the case where the right Hilbert $A$--module $X$ is not countably 
generated, unconditionally convergent countable bases will be replaced in the sequel
by {\it generalized bases}, in the following sense.\medskip

\noindent{\bf 1.8 Definition} 
Consider a set $\Lambda$ and, for each finite subset $\mu\subset\Lambda$,
let $u_\mu$ be a finite subset of $X$ with $|u_\mu|=|\mu|$. Let us endow the set
of finite subsets of $\Lambda$ with the partial order defined by inclusion.
The net $\mu\to u_\mu$
will be called a {\it generalized (right) basis} of $X$ if 
\roster
\item for all $x\in X$,  $\sum_{y\in u_\mu}(x| y)_A(y|x)_A\leq 
\sum_{y\in u_\nu}(x| y)_A(y|x)_A$, if $\mu\subset \nu$,
\item $x=\lim_\mu\sum_{y\in u_\mu} y(y|x)_A$.
\endroster
\medskip

Let $\Lambda$ be a set, and $\mu\subset\Lambda\to u_\mu\subset
X$  a net with $|u_\mu|=|\mu|$.
One can easily see that
  $\mu\to u_\mu$ is a generalized basis if and only if
$\mu\to T_\mu:=\sum_{y\in u_\mu}\theta^r_{y,y}$ is an increasing 
approximate unit
of $\K(X_A)$ with norm $\leq 1$.\medskip

\noindent{\bf 1.9 Proposition} {\sl Any right (or left) Hilbert
$C^*$--module $X$ admits
a generalized right (or left) basis}\medskip

\noindent{\it Proof}
This is essentially the proof of existence
 of an approximate unit of a $C^*$--algebra
with entries in a dense ideal (cf. Prop. 1.7.2 in   \cite{Di}).
Let $\Lambda$ be $X$, and, for each finite subset $\mu=\{x_1,\dots, x_n\}$ of
$\Lambda$ let $u_\mu=\{(1/n+\sum_1^n\theta_{x_j,x_j})^{-1/2}x_i, i=1,\dots,n\}$.
Then the proof of the cited result shows that the net  
$T_\mu:=\sum_{y\in u_\mu}\theta_{y,y}=
(\sum_1^n\theta_{x_j,x_j})(1/n+\sum_1^n\theta_{x_j,x_j})^{-1}$ 
is increasing, with norm $\leq 1$ and satisfies 
$\sum_1^n(T_\mu-I)\theta_{x_i,x_i}(T_\mu-I)\leq\frac{1}{4n}$. So for
 $i=1,\dots,n$,
$$\|(T_\mu-I)x_i\|^2\leq 
\|\sum_1^n(T_\mu-I)\theta_{x_i,x_i}(T_\mu-I)\|\leq\frac{1}{4n}.$$
\medskip

\heading
2. Bimodules of finite index
\endheading

In this section we study the notion of 
$C^*$--bimodules of finite right index. We start with a weak
notion of finite  index, based only on a Pimsner-Popa-type
inequality,
and we construct the right index element of ${}_AX_B$ as a positive
central
element in the enveloping von Neumann algebra $A''$. Later on we shall
concentrate on those 
 bimodules for which the index element lies in the multiplier
algebra of
$A$, and we perform, in this case, the analogue of the Jones basic
construction with nice
properties. We also prove 
that in the case where
$A$ is unital, the index element belongs to $A$ if and only if 
$X$ is finitely generated as a right $B$--module. 
This result will be stated
in a more general form, which includes the non 
unital case.
\medskip 

\noindent{\bf 2.1 Bimodules of finite  right numerical index}
\bigskip

\noindent{\bf  2.1 Definition} 
 Let $A$ and $B$ be \cst--algebras and $X = {}_AX_B$ a bimodule
 over the complex algebras $A$ and $B$. We say that $X$ is a  {\it  right  
Hilbert $A$--$B$ bimodule\/} 
if
\roster
\item $X$, as a right $B$-module, is endowed with a 
$B$-valued inner product making it into a 
  right
Hilbert $B$--module,
\item
 for all $a\in A$,  
the map $\phi(a): x\in X\mapsto ax\in X$ is adjointable, with adjoint
 $\phi(a)^*=\phi(a^*)$. 
\endroster
Therefore
  $\phi: a\in A\to \phi(a)\in \L(X_B)$ is  a $^*$--homomorphism from $A$
to the algebra
$\L(X_B)$ of right adjointable maps on $X_B$.
The map $\phi$ will be
referred to as the left action of $A$ on $X$. 

 We  introduce the notion of {\it left\/} Hilbert $A$--$B$ 
bimodule in a similar manner.
Thus,  if $X$ is a left Hilbert $A$--$B$  bimodule, the map $\psi: 
B\to\L({}_AX)$, 
$\psi(b): x\in X\mapsto xb\in X$, for all $b\in B$, 
and referred to as  the right action of $B$ on $X$, is a
$^*$--antihomomorphism from
$B$ to the algebra $\L({}_AX)$ of left adjointable maps on ${}_AX$. 
\medskip

Notice that left and right actions on a right (or left) Hilbert bimodule
are  not assumed to
be faithful. In the following proposition we give a sufficient condition.
Recall that a closed ideal $J$ in a $C^*$--algebra $B$ is called 
{\it essential\/}
if each nonzero closed ideal of $B$ has a nonzero intersection with $J$ 
(see 3.12.7 in \cite{P}).\medskip

\noindent{\bf  2.2 Proposition} {\sl Let $X$ be a right
 pre-Hilbert $B$--module (resp. left pre-Hilbert
$A$--module). 
If the closed linear span in $B$ (resp. $A$) of inner products $(x| y)_B$
(resp. ${}_A(x| y)$) $x,y\in X$)
is an essential ideal of $B$ (resp. $A$), the equation
$Xb=0$  for some $b\in B$ 
(resp. $aX=0$
for some $a\in A$)
implies $b=0$ (resp. $a=0$).}\medskip

\noindent{\it Proof} Let $J$ denote the closed linear span in $B$
 of right inner products.  If $b\in B$ satisfies $yb=0$ for all 
$y\in X$ then
$(x| yb)_B=(x| y)_Bb=0$ for all $x,y\in X$, hence $jb=0$ for all $j\in J$. 
By Prop. 3.12.8 in \cite{P} the natural 
injection of $J$ into its multiplier algebra $M(J)$
extends to an  injection $B\to M(J)$. Thus reading the above equation in 
 $M(J)$, we get  $jb=0$ for all $j\in M(J)$ and this implies that $b=0$ 
since $M(J)$ is unital.\medskip

\noindent{\bf 2.3 Definition}
A $A$--$B$  bimodule $_{A}X_B$ will be called 
{\it bi-Hilbertian} if it is endowed with a
right as
well as a left Hilbert $A$--$B$ $C^*$--bimodule structure
in such a way that the two Banach space norms arising from the two inner
products are equivalent. In other words,
$X$ is a left Hilbert module over $A$ and a right Hilbert module over $B$
such that the left $A$--action $\phi$ and the right $B$--action
$\psi$ are $^*$--preserving maps into 
the algebras of right adjointable and left adjointable
operators, respectively.
Furthermore
there  should exist two constants $\lambda$, $\lambda'>0$
such that, for $x\in X$,
$$\lambda'\|(x|x)_B\|\leq\|{}_{A}(x|x)\|\leq\lambda\|(x|x)_B\|.$$
\medskip

The inequality at the left hand side always extends to
 finite sums, in the sense of the following proposition.
 \medskip

\noindent{\bf 2.4 Proposition} {\sl Let 
${}_AX_B$ be a bi-Hilbertian $C^*$--bimodule, 
and let $\lambda'>0$ satisfy 
$\lambda'\|(x|x)_B\|\leq\|{}_{A}(x|x)\|$, $x\in X$. Then for all
$n\in{\Bbb N}$ and for all $x_1,\dots,x_n\in X$ we have
$$\lambda'\|\sum_1^n\theta^r_{x_i,x_i}\|\leq
\|\sum_1^n{}_{A}(x_i|x_i)\|.$$}\medskip

\noindent{\it Proof} Let $T\in M_n(B)$ be the positive matrix 
whose $(i,j)$-th entry
is $(x_i|x_j)_B$. Notice that
$$\|\sum_1^n\theta^r_{x_i,x_i}\|=
\|T\|=\sup\|\sum_{i,j=1}^n{b_i}^*(x_i|x_j)_Bb_j\|$$
where the supremum is taken over all the $n$-tuples $(b_1,\dots,b_n)$ with 
elements
in $B$ such that $\|\sum_j{b_j}^*b_j\|=1$. Now the norm at the right hand
side
coincides with the norm of $(y|y)_B$, where $y=\sum_j x_jb_j$, therefore
$$\lambda'\|(y|y)_B\|\leq\|{}_A(y|y)\|=
\|\sum_{i,j}{}_A(x_ib_i{b_j}^*|x_j)\|\leq
\|\sum_i{}_A(x_i|x_i)\|,$$ and the proof is now complete.
\medskip

On the contrary,  there may exist no
$\lambda>0$ for which
$\|\sum_{i=1}^n{}_A(x_i|x_i)\|\le
 \lambda\| \sum_{i=1}^n \theta^r_{x_i,x_i}\|$ for all $n\in{\Bbb N}$ and
all $x_1,
\dots,x_n\in X$, as the following elementary example shows.\medskip

\noindent{\bf  2.5 Example}
Let $A = B = \comp$ and let $H = \ell^2(\bN)$ be an infinite dimensional 
Hilbert 
space, regarded  as a bi-Hilbertian
${\Bbb C}$--${\Bbb C}$ bimodule in the natural way.  
Let 
${e_1,e_2,....}$ be a countable orthonormal subset of $H$.  
Then for all $n\in{\bN}$, $ \| \sum_{i=1}^n \theta^r_{e_i,e_i}\| = 1$, while   
 $\| \sum_{i=1}^n{}_{\Bbb C}(e_i|e_i)\| = n$. \medskip 

In fact, the existence of such a constant $\lambda$ will lead us to
the notion of {\it finite right numerical index} of $X$. 
We
anticipate a lemma.\medskip

\noindent{\bf 2.6 Lemma} {\sl Let ${}_AX_B$ be a right Hilbert $C^*$--bimodule,
and let $x,y\in X\to{}_A(x|y)$ be an $A$--valued, biadditive, left $A$--linear,
right $A$--antilinear form on $X$ such that
${}_A(x|y)^*={}_A(y|x)$ and ${}_A(x|x)\geq 0$ for all $x,y\in X$.
If this form is continuous,  in the sense that there is $\lambda>0$ such
that $\|{}_A(x|x)\|\leq\lambda\|(x|x)_B\|$ for $x\in X$,
and if the right $B$--action is adjointable
with respect to this form (i.e.${}_A(xb|y)={}_A(x|yb^*)$, $x,y\in X$, $b\in B$),
 there exists a unique additive map
$F: FR(X_B)\to A$ such that $F(\theta^r_{x,y})={}_A(x|y)$.
$F$ satisfies the following properties:\roster
\item (positivity) $F(T^*T)\geq0$, for  $T\in FR(X_B)$,
\item (reality) $F(T^*)=F(T)^*$, for  $T\in FR(X_B)$,
\item ($A$-bilinearity) $F(\phi(a)T)=aF(T)$, 
$F(T\phi(a))=F(T)a$
for $a\in A$, $T\in FR(X_B)$,
\item (Pimsner-Popa inequality) if $X$ is bi-Hilbertian and if
$\lambda'>0$
satisfies $\lambda'\|(x|x)_B\|\leq\|{}_A(x|x)\|$ for all $x\in X$ then
$\|F(T)\|\geq\lambda'\|T\|$ for any $T\in FR(X)$ that can be written as a 
finite sum of operators of the form $\theta^r_{x,x}$.
\endroster} \medskip    

\noindent{\it Proof} Uniqueness is obvious.
Let 
$\mu\subset\Lambda\to u_\mu$ be a generalized right basis
of the right Hilbert module $X_B$, which exists by Proposition 1.9.
Consider the linear map $F_\mu: T\in\L(X_B)\to\sum_{y\in u_\mu}{}_A(Ty|y)
\in A$.
Note that $$F_\mu(\theta^r_{x,z})=\sum_{y\in u_\mu}{}_A(x(z|y)_B|y)=
{}_A(x|\sum_{y\in u_\mu}y(y|z)_B).$$ Since $\lim_\mu\sum_{y\in
u_\mu}y(y|z)_B=z$
in the norm $\|(\cdot|\cdot)_B\|^{1/2}$, and since $\|{}_{A}(x|x)\|\leq
\lambda\|(x|x)_{B}\|$, we also
have that $\lim_\mu\sum_{y\in u_\mu}y(y|z)_B=z$ in the seminorm
defined by the left inner product. Thus $\lim_\mu
F_\mu(\theta^r_{x,z})={}_A(x|z)$.
Let us define $F$ as the pointwise norm limit of the net $\mu\to F_\mu$
on $FR(X_B)$. Obviously  this limit does not depend on the generalized
right basis.
 (1) follows from the fact that any element of the form $T^*T$, with $T\in
FR(X)$, can be written
as a finite sum of elements of the form $\theta^r_{x,x}$.
Properties (2) and (3) are easy to check.
  (4) follows from Prop. 2.4.
\medskip

The map  $F$ will be referred to as
the {\it additive extension of the form} ${}_A(\cdot|\cdot)$ to the 
 finite rank operators
on $X_B$.

Notice that a bimodule satisfying the properties 
of the previous lemma is almost bi-Hilbertian.
The only missing properties are the fact that the seminorm coming from the
left-linear $A$--valued form
is in fact a norm, and completeness of $X$ with respect to this norm.
\medskip

\noindent{\bf 2.7 Proposition} {\sl Let $X$ be a right Hilbert $A$--$B$ 
$C^*$--bimodule and let $x,y\to{}_A(x|y)$ be an $A$--valued form
on $X$ satisfying the same properties as in the previous lemma
(with right seminorm not necessarily a Banach space norm).
Then following properties are equivalent.
\roster
\item There exists $\lambda>0$ such that for all $n\in{\Bbb N}$ and for all
 $x_1,\dots,x_n\in X$,
$$\|\sum_1^n{}_A(x_i|x_i)\|\leq\lambda\|\sum_1^n\theta^r_{x_i,x_i}\|,$$
\item there exists $\lambda>0$ such that for all $n\in{\Bbb N}$ and for all
 $x_1,\dots,x_n$, $y_1,\dots, y_n\in X$,
$$\|\sum_1^n{}_A(x_i|y_i)\|\leq\lambda\|\sum_1^n\theta^r_{x_i,y_i}\|,$$
\item 
$F(T)\geq0$ for any $T\in FR(X_B)\cap\K(X_B)^+$ and 
$\sup_\mu\|F(\sum_{y\in u_\mu}\theta^r_{y,y})\|$ is finite for some
 generalized right  basis
$\mu\to u_\mu$ of $X_B$.
\endroster
If one of these conditions is satisfied, the smallest constants for which
(1) and (2) hold, 
 coincide and equal, in turn, $\sup_\mu\|F(\sum_{y\in
u_\mu}\theta^r_{y,y})\|$. In particular, the latter does not depend
on the generalized right basis. 
}\medskip

\noindent{\it Proof} (1)$\Rightarrow$ (2)
Let 
$\mu\subset\Lambda\to u_\mu$ be a generalized basis
of the right Hilbert module $X_B$.
Consider the linear map $F_\mu: T\in\L(X_B)\to
\sum_{y\in u_\mu}{}_A(Ty|y)\in A$, already considered
 in the
proof of Lemma 2.6.
 We claim that 
$\|F_\mu\|\leq\lambda$ for any $\mu$.
 We show the claim. For any $T\in\L(X_B)$,
$$\|F_\mu(T)\|=\|\sum_{y\in u_\mu} {}_A(Ty|y)\|
\leq\|\sum_{y\in u_\mu} {}_A(Ty|Ty)\|^{1/2}
\|\sum_{y\in u_\mu}{}_A(y|y)\|^{1/2}$$
by the Cauchy--Schwarz inequality of the left inner product
(see, e.g., \cite{B} 
Prop. 13.1.3). Now by our assumption the last term is bounded above by
$$\lambda
\|\sum_{y\in u_\mu} \theta^r_{Ty,Ty}\|^{1/2}
\|\sum_{y\in u_\mu} \theta^r_{y,y}\|^{1/2}=$$
$$\lambda\| T\sum_{y\in u_\mu}\theta^r_{y,y}T^*\|^{1/2}
\|\sum_{y\in u_\mu} \theta^r_{y,y}\|^{1/2}\leq\lambda\|T\|.$$
We have already seen that
 $\lim_\mu F_\mu(\theta^r_{x,z})={}_A(x|z)$.
Since, for any $x_1,\dots,x_n$, $z_1,\dots,z_n\in X$,
 $\|F_\mu(\sum_1^n\theta^r_{x_i,z_i})\|\leq\lambda
\|\sum_1^n\theta^r_{x_i,z_i}\|$,  the proof is completed taking the norm limit
at the left hand side.

(2)$\Rightarrow$(3)
We know that the linear extension $F$ of the left inner product to the
finite rank operators on $X_B$ is positive on $FR(X)$ (Lemma 2.6). We
first show
that $F$ is still positive on $FR(X)\cap\K(X_B)^+$. 
By (2) $F$ is norm continuous, therefore, if $T\in FR(X)\cap\K(X_B)^+$
and if $\mu\to u_\mu\subset X$ is a generalized basis of $X$, the
net $T^{1/2}\sum_{y\in u_\mu}\theta^r_{y,y}T^{1/2}$ converges to $T$ in
norm, therefore $F(T)=\lim_\mu \sum_{y\in u_\mu}F(\theta^r_{T^{1/2}y,
T^{1/2}y})\in A^+$. Furthermore for all $\mu$, $\|F(\sum_{y\in
u_\mu}\theta^r_{y,y})\|\leq\lambda\|\sum_{y\in
u_\mu}\theta^r_{y,y}\|\leq\lambda$.

(3)$\Rightarrow$(1) Let $x_1,\dots,x_n$ be elements of $X$, and set
$T=\sum_{1}^n\theta^r_{x_i,x_i}$.
Let $\mu\to u_\mu$ be a generalized right basis of $X$. Since
$\sum_{y'\in u_\mu}\theta^r_{y',y'}$ is an approximate unit of $\K(X_B)$
and since the left inner product is continuous with respect to the right
one, for all $\mu$, the net $\mu'\to\sum_{y\in u_\mu, y'\in u_{\mu'}}
{}_A(\theta^r_{y',y'}Ty|y)$ converges to $\sum_{y\in u_\mu}{}_A(Ty|y)$
in norm. On the other hand this net coincides with
$F((\sum_{y'\in u_{\mu'}}\theta^r_{y',y'})T(\sum_{y\in
u_\mu}\theta^r_{y,y}))$.
 The form $S,T\in FR(X_B)\to F(ST^*)$
is left $A$-linear, right $A$-antilinear, symmetric and positive and
therefore 
 it
satisfies the Cauchy-Schwarz inequality
$\|F(ST^*)\|^2\leq\|F(SS^*)\|\|F(TT^*)\|$.
It follows that
$$\|F((\sum_{y'\in u_{\mu'}}\theta^r_{y',y'})T(\sum_{y\in
u_\mu}\theta^r_{y,y}))\|^2\leq$$
$$\|F(
(\sum_{y'\in u_{\mu'}}\theta^r_{y',y'})TT^*(\sum_{y'\in
u_{\mu'}}\theta^r_{y',y'})
)\|\|F(
(\sum_{y\in
u_\mu}\theta^r_{y,y})^2
)\|.$$
Now
$$(\sum_{y'\in u_{\mu'}}\theta^r_{y',y'})TT^*(\sum_{y'\in
u_{\mu'}}\theta^r_{y',y'})\leq\|T\|^2(\sum_{y'\in
u_\mu'}\theta^r_{y',y'})$$
and 
$(\sum_{y\in
u_\mu}\theta^r_{y,y})^2\leq(\sum_{y\in
u_\mu}\theta^r_{y,y})$, so, applying $F$, we deduce that the above term is
bounded above by
$\|T\|^2{\lambda_0}^2$ 
where $\lambda_0=
\sup_\mu\|F(\sum_{y\in u_\mu}\theta^r_{y,y})\|$.
Passing first to the limit over $\mu'$ and then over $\mu$ we deduce that
(1) holds with $\lambda=\lambda_0$.
 
It is now clear from the proof that if one of these three equivalent
conditions holds, the best constants satisfying (1) and (2) coincide, an
coincide in turn with
$\sup_\mu\|F(\sum_{y\in u_\mu}\theta^r_{y,y})\|$.\medskip
 
\noindent{\bf 2.8 Definition}
A  $C^*$--bimodule satisfying one of the equivalent
properties described in the previous proposition will be 
called of {\it  finite right numerical index}. The 
corresponding smallest
positive constant
  will be called the {\it right numerical  index} of $X$, and denoted
$r-I[X]$.
\medskip

 Let $X$ be an  $A$--$B$ bimodule. The {\it 
contragradient bimodule\/}  of $X$ is 
the $B$--$A$ bimodule  $ \overline{X} = \{\overline{x} ; x \in X \}$ with
complex conjugate
vector space structure and 
  bimodule structure
  given by
 $$b\cdot \overline{x} = \overline{xb^*}, \quad 
 \overline{x}\cdot a = \overline{a^*x},\quad b\in B, a\in A.$$
If $X$ is a  right (left) Hilbert $A$--$B$ $C^*$--bimodule,
$\overline{X}$ becomes a left (right) Hilbert  $B$--$A$ $C^*$--bimodule
with inner product given by:
$${}_B(\overline{x}|\overline{y}) = (x|y)_B$$ 
 $$((\overline{x}|\overline{y})_A = {}_A(x|y).)$$
Therefore if ${}_AX_B$ is bi-Hilbertian, ${}_B\overline{X}_A$ is
bi-Hilbertian as well.
\medskip

\noindent{\bf 2.9 Definition}
We will say that
    ${}_AX_B$ is  of {\it finite left numerical index}\/ if   
the contragradient bimodule ${}_B\overline{X}_A$ is of finite right
numerical 
index. 
Its left numerical index 
is defined by $\ell-I[X]:=r-I[\overline{X}]$.

A bi-Hilbertian bimodule of   finite left and right numerical
indices will be simply called of finite numerical index. Its
numerical index is
defined
by
$I[X]:=(r-I[X])(\ell-I[X])$.
\medskip

\noindent{\bf 2.10 Corollary} {\sl 
Let ${}_AX_B$ be a bi-Hilbertian $C^*$--bimodule, and, for $n\in{\bN}$,
let us
consider $Y_n:=\oplus_1^n X$ as a $M_n(A)$-$B$ bimodule in the natural way.
Endow $Y_n$ with the following forms:
${}_{M_n(A)}(\underline{x}|\underline{y})=({}_A(x_i|y_j))$, $(\underline{x}
|\underline{y})_B=\sum_1^n(x_i|y_i)_B$,
where $\underline{x}=(x_1,\dots,x_n)$, $\underline{y}=(y_1,\dots,y_n)$.
\roster
\item
If 
for some $\lambda>0$, $\|{}_A(x|x)\|\leq\lambda\|(x|x)_B\|$ then
also $\|{}_{M_n(A)}(\underline{x}|\underline{x})\|\leq\lambda\|
(\underline{x}|\underline{x})_B\|$. 
\item If $X$ is of finite right numerical index,
$Y_n$ becomes a $C^*$--bimodule of finite right numerical index and
$r-I[Y_n]=r-I[X]$ for all $n\in{\Bbb N}$.
\item
If  $X$ is of finite  numerical index, $Y_n$ is
bi-Hilbertian
and of finite left numerical index (and hence of finite numerical index, 
by (1)):
$\|(\underline{x}|\underline{x})_B\|\leq \ell-I[X]\|{}_{M_n(A)}
(\underline{x}|\underline{x})\|$ and
$\ell-I[X]=\ell-I[Y_n]$, for all $n\in{\Bbb N}$.\endroster}\medskip

\noindent {\it Remark} Notice that the  constants comparing the two norms
on $Y_n$
do not depend on $n$.
\medskip

\noindent{\it Proof} (1) It is easy to check that $Y_n$ is a right Hilbert
$C^*$--bimodule.
Now $\|{}_{M_n(A)}
(\underline{x}|\underline{x})\|=\|\sum_1^n\theta^l_{x_i,x_i}\|$
and $\|(\underline{x}|\underline{x})_B\|=\|\sum_1^n(x_i|x_i)_B\|$. 
Therefore
if 
for some $\lambda>0$, $\|{}_A(x|x)\|\leq\lambda\|(x|x)_B\|$ then
also $\|{}_{M_n(A)}(\underline{x}|\underline{x})\|\leq\lambda\|
(\underline{x}|\underline{x})_B\|$ by Prop. 2.4 applied to the
contragradient bimodule. 
The remaining properties of the 
 left $M_n(A)$--valued form are easily checked. 

(2)
There is a natural 
identification
of 
$FR(Y_n)$ with $M_n({\Bbb C})\otimes FR(X)$ taking
$\theta^r_{\underline{x},
\underline{y}}$ to the
matrix
$(\theta^r_{x_i,y_j})$. Under this identification, the map $F_n:FR(Y_n)\to
M_n(A)$ obtained
extending additively the left form of $Y_n$ identifies with ${\text
id}_{M_n}\otimes F$,
where
$F$ is the additive extension of the left form of $X$. Notice that
$FR(Y_n)\cap\K(Y_n)^+$ identifies with 
$(M_n\otimes FR(X))
\cap(M_n\otimes\K(X))^+$, so if 
$T=(T_{i,j})\in(M_n\otimes FR(X))\cap(M_n\otimes \K(X))^+$,
$\sum_{i,j}{a_i}^*F(T_{i,j})a_j=
F(\sum_{i,j}\phi(a_i)^*T_{i,j}\phi(a_j))\in A^+$
for all $a_1,\dots,a_n\in A$.
Hence $F_n$ takes positive values on 
$FR(Y_n)\cap\K(Y_n)^+$. If now $\mu\to u_\mu$ be a generalized basis
of $X$,
$\mu'=\mu\oplus\dots\oplus\mu\to{u'}_{\mu'}=
\{(y,0,\dots,0),\dots(0,\dots,y),y\in
u_\mu\}$ is a generalized
basis of $Y_n$ with corresponding approximate unit of $\K(Y_n)$ given by
the diagonal matrix of $\sum_{y\in u_\mu}\theta^r_{y,y}$.
Therefore the norm of the evaluation of $F_n$ over this diagonal matrix
coincides with $\|\sum_{y\in u_\mu}{}_A(y|y)\|$, and this concludes the
proof of (2).

(3)
Since $\|{}_{M_n(A)}
(\underline{x}|\underline{x})\|=\|\sum_1^n\theta^l_{x_i,x_i}\|$
and $\|(\underline{x}|\underline{x})_B\|=\|\sum_1^n(x_i|x_i)_B\|$, 
if $X$ is of finite left numerical index,
$\|(\underline{x}|\underline{x})_B\|\leq (\ell-I[X])\|{}_{M_n(A)}
(\underline{x}|\underline{x})\|$, therefore, taking into account (1), 
$Y_n$ is bi-Hilbertian.
The compacts on $Y_n$ with respect to the left inner product identify with
$\K({}_AX)$ via
$\theta^l_{\underline{x},\underline{y}}\to\sum_1^n\theta^l_{x_h,y_h}$.
This shows that $$(\ell-I[X])\|\sum_{j=1}^m\theta^l_{\underline{x^j},
\underline{x^j}}\|=
(\ell-I[X])\|\sum_{j,k}\theta^l_{x^j_k,x^j_k}\|\geq$$
$$\|\sum_{j,k}(x^j_k|x^j_k)_B\|=
\|\sum_j(\underline{x}^j|\underline{x}^j)_B\|$$
so $Y_n$ has finite left numerical index and  $\ell-I[Y_n]=\ell-I[X]$.

\medskip

The following result is the first step towards the Jones basic 
construction.
\medskip

\noindent{\bf 2.11 Corollary} {\sl If $X$ is a bi-Hilbertian $A$-$B$
$C^*$--bimodule of
finite right numerical  index, the additive extension $F$ of the left
inner
product to $FR(X_B)$ extends uniquely to a norm continuous
map $F:\K(X_B)\to A$. One has: $\|F\|=r-I[X]$. This extension, still
denoted by $F$,  is   positive,
$A$-bilinear (in the sense that $F(\phi(a)T)=aF(T)$ and
$F(T\phi(a))=F(T)a$ for
$a\in A$, $T\in\K(X_B)$) and has range contained in the closed ideal of
left inner products. Moreover one has $\phi(F(T))\geq\lambda' T$
for all
$T\in\K(X_B)^+$, where $\lambda'$ is the best constant for which
$\lambda'\|(x|x)_B\|\leq\|{}_A(x|x)\|$.}\medskip

\noindent{\it Proof} The only assertion that
is not obvious yet
is the inequality $\phi(F(T))\geq \lambda' T$ for $T\in \K(X_B)^+$.
Now  part (4) in Prop. 2.6 implies that $\|F(T)\|\geq\lambda'
\|T\|$
for $T\in\K(X_B)^+$. Left $A$--action is faithful on the norm closed
ideal
$\J$
generated by left inner products: $\phi(j)=0$ for some $j\in \J^+$
implies $0={}_A(x|jy)={}_A(x|y)j$ for $x,y\in X$ and therefore $j=0$.
Since $F(T)\in\J$ for all $T\in\K(X_B)$, $\phi$ is isometric on $\J$,
therefore for all  $T\in\K(X_B)^+$,
$\|\phi(F(T))\|=\|F(T)\|\geq\lambda'\|T\|$.
Arguing as in \cite{FK}, with the map $\phi\circ F$ in place of
a conditional expectation, we deduce the desired inequality.
\medskip

Pimsner--Popa conditional expectations provide typical examples of
bimodules of finite right numerical index.\medskip

\noindent{\bf 2.12 Proposition} {\sl
 Let $A\subset B$ be an inclusion of $C^*$--algebras
and
let $E: B\to A$ be a conditional expectation with fixed point set $A$.
Assume that  $\|E(b)\|\geq \lambda \|b\|$, 
for all positive elements $b\in B$ and for some $\lambda>0$.
\roster
\item
Consider ${}_BX_A=B$ as a $B$--$A$ bimodule in the natural way, and with
inner  products $(x|y)_A=E(x^*y)$, ${}_B(x|y)=xy^*$. 
Since $\|(x|x)_A\|\leq\|{}_B(x|x)\|\leq{\lambda}^{-1}\|(x|x)_A\|$,
 $X$ is bi-Hilbertian. By  \cite{FK},
 there is a constant $\lambda'>0$ such that $E-\lambda'$ is
completely positive. Let us choose the best such $\lambda'$.
Then $X$ has finite right numerical index and
$r-I[X]={\lambda'}^{-1}$.
\item Consider now ${}_AY_B=B$ as a $A$--$B$ bimodule with inner products
$(x|y)_B=x^*y$ and ${}_A(x|y)=E(xy^*)$. Then the $B$--$A$
antilinear map $X\to Y$  induced by the $^*$--involution of $B$ 
identifies $Y$ with the contragradient $\overline{X}$ of $X$.
Therefore $X$ is of finite left numerical index and
$\ell-\text{Ind}[X]=1$.
\endroster
}\medskip

\noindent{\it Proof}
(1) For
all $n\in{\bN}$, and all
$x_1,\dots,x_n\in X$,
$${\lambda'}^{-1}\|\sum_1^n\theta^r_{x_i,x_i}\|=
{\lambda'}^{-1}\|(E({x_i}^*x_j))_{i,j}\|_{M_n(A)}\geq$$
$$\|({x_i}^*x_j)_{i,j}\|_{M_n(B)}=\|\sum_1^n{}_B(x_i|x_i)\|,$$
i.e. $X$ is of finite right numerical index in our sense and
$r-I[X]\leq{\lambda'}^{-1}$. 
On the other hand, let 
$\mu\to u_\mu$ be a generalized basis of $X_A$, and 
set, for every $\mu$ and every $x\in X$, $x_\mu:=\sum_{y\in u_\mu}
yE(y^*x)$. Then
$${x_\mu}^*x_\mu\leq
\|\sum_{y\in u_\mu}yy^*\|
\sum_{y\in u_\mu}E(x^*y)E(yx^*)=$$
$$E(x^*x_\mu)\|\sum_{y\in u_\mu}yy^*\|\leq
\sup_\mu\|\sum_{y\in u_\mu}yy^*\|E(x^*x_\mu).$$ Taking the limit
over $\mu$, we are led
to the
inequality
$r-I[X]\geq{\lambda}^{-1}$.
Consider now the inclusion $M_n\otimes A\subset M_n\otimes B$ and the
conditional
expectation $E_n:=\text{id}\otimes E$, which satisfies
$E_n(b)\geq\lambda'b$, $B\in M_n(B)^+$. Cor. 2.10 shows
that $\oplus_1^n B$ is a $M_n(B)$-$A$ bimodule with the same right index
as $B$, hence, combining with the above argument, we deduce that
$r-I[X]\geq{\lambda'}^{-1}$.

The proof of part (2) is easy, therefore we omit it.
\medskip

\noindent{\it Remark} If ${}_BX_A$ and ${}_AX_B$ arise from a
Pimsner--Popa conditional
expectation $E$, as the previous proposition, the
 corresponding map  $F_Y$ constructed in Cor. 2.11 reduces to $E$
itself. More interestingly, 
$F_X:\K(X_A)\to B$  
 is related to the construction
of the   {\it dual} conditional
expectation.
However, if the index of $E$, as an element of $Z(B'')$ (cf
Def. 2.17),
does not belong to the multiplier algebra of $B$, $F_X$ is not a multiple
of a conditional expectation.
\bigskip

\noindent{\bf 2.2 Tensoring bi-Hilbertian
$C^*$--bimodules}\bigskip

\noindent In this subsection we analyze the behaviour of bi-Hilbertian
bimodules
under taking their tensor products.
We show that the algebraic tensor product
$X\odot_B Y$ of bi-Hilbertian $C^*$--bimodules can be made into a
bi-Hilbertian
bimodule in a natural way if  $X$ is of finite right
numerical index
 and $Y$ is of finite left numerical index, and that this is also
a necessary condition in general.

The problem of studying conditions under which $X\otimes_B Y$ is 
of finite index will be considered in section 5 (cf. Theorem 5.2).

Let ${}_AX_B$ and ${}_BY_C$ be bi-Hilbertian $C^*$--bimodules. Then the
algebraic tensor product
$X\odot_B{Y}$ is an $A$-$C$ bimodule in a natural way,
also endowed with a right and a left
pre-bi-Hilbertian structure:
$$
(x_1\otimes {y_1}|x_2\otimes {y_2})_C 
= ({y_1}|(x_1|x_2)_B{y_2})_C,
$$
$$
{}_A(x_1\otimes {y_1}|x_2\otimes {y_2})=
{}_A(x_1{}_B({y_1}|{y_2})| x_2). 
$$
Therefore $X\odot_B{Y}$ can be made into a right Hilbert
 $C$--module
$X\otimes^r_B{Y}$ completing with respect to the first inner
product 
and also
into a left Hilbert  $A$--module $X\otimes^\ell_B Y$ completing
with respect to the
second inner product    
(always after dividing out by vectors of  seminorm zero). 

Under which conditions these two
seminorms are
equivalent on the algebraic tensor product $X\odot_B Y$?

Choosing for $Y$ the strong Morita equivalence
$_B\ell^2(B)_{\K\otimes B}$
with inner products 
${}_B(\underline{b}|\underline{b'})=
\sum_j{b_j}{{b'}_j}^*$, 
$(\underline{b}|\underline{b'})_{\K\otimes B}=\sum
\delta_{i,j}\otimes {b_i}^*{b'}_j$, $X\otimes^r_B \ell^2(B)$ 
identifies with
$\ell^2(X)$ with inner products
${}_A(\underline{x}|\underline{x'})=
\sum_j{}_A(x_j|{x'}_j)$, 
$(\underline{x}|\underline{x'})_{\K\otimes B}=\sum
\delta_{i,j}\otimes ({x_i}|{x'}_j)_B$.
Therefore the left and right seminorms on $X\odot_B\ell^2(B)$
are equivalent if and only if $X$ is of finite right numerical index.
Similarly, the left and right seminorms on
$\overline{\ell^2(B)}
\odot_B Y$, with 
$\overline{\ell^2(B)}$ the inverse  strong Morita equivalence,
are  equivalent if and only if $Y$ is of finite left numerical index.
We show that these necessary conditions on $X$ and $Y$ are also
sufficient.
\medskip

 \noindent{\bf 2.13 Proposition} {\sl Let ${}_AX_B$ and ${}_BY_C$ be
bi-Hilbertian $C^*$--bimodules. Assume that $X$ is of finite right
numerical index and that $Y$ is of finite left numerical index. 
Let $F_X:\K(X_B)\to A$, $F_{\overline{Y}}:\K(\overline{Y}_B)\to C$ be
the corresponding maps constructed in Cor. 2.11.
Then
\roster
\item
the two seminorms arising from the left and right
inner products on $X\odot_B{Y}$ as above are equivalent.
Therefore $X\otimes^r_B {Y}=X\otimes^\ell_B{Y}(=:X\otimes_B Y)$
and it is
 a bi-Hilbertian $A$-$C$ bimodule.
\item
Consider
$\K(\overline{Y}_B, X_B)$ as a
  $A$--$C$ bimodule with left and right 
  inner products 
 ${}_A(T|S) = F_X(TS^*)$ and $(T|S)_C = F_{\overline{Y}}(T^*S)$.  Then
$\K(\overline{Y}_B, X_B)$
is complete in any of the induced norms, and becomes in this way a
bi-Hilbertian $C^*$--bimodule. 
\item
The map $x\otimes{y}\in X\otimes
{Y}\to\theta^r_{x,\overline{y}}\in\K(\overline{Y}_B,X_B)$ extends 
to a 
 bijective $A$-$C$ bimodule map
 $U : {}_AX_B \otimes_B {}_B{Y}_A \rightarrow \K(\overline{Y}_B, X_B)$ 
preserving the left and right inner products.
\endroster
}\medskip

\noindent{\it Proof}
Consider $X$ and $\overline{Y}$ as right Hilbert $B$--modules, and
define the map $U:X\odot_BY\to FR(\overline{Y}, X)$ associating
$\theta^r_{x,\overline{y}}$ to the simple tensor $x\otimes y$.
$U$ is a well defined $A$-$C$-bimodule map.
For any $x_1, x_2\in X, y_1, y_2 \in Y$, we have 
$$
(x_1\otimes {y_1}|x_2\otimes {y_2})_C 
= ({y_1}|(x_1|x_2)_B{y_2})_C
= {}_C(\overline{y_1}|\overline{y_2}(x_2|x_1)_B)
$$
and

$$ 
\align
 (\theta^r_{x_1,\overline{y_1}}|\theta^r_{x_2,\overline{y_2}})_C 
 & = F_{\overline{Y}}(\theta^r_{x_1,\overline{y_1}}\ ^*
\theta^r_{x_2,\overline{y_2}})
 = F_{\overline{Y}}(\theta^r_{\overline{y_1}(x_1|x_2)_B,\overline{y_2}})
\\
 & = {}_C(\overline{y_1}(x_1|x_2)_B|\overline{y_2}) 
= {}_C(\overline{y_1}|\overline{y_2}(x_2|x_1)_B) .
\endalign
$$ 
Similarly we have 
$$
{}_A(x_1\otimes {y_1}|x_2\otimes {y_2})
={}_A (\theta^r_{x_1,\overline{y_1}}|\theta^r_{x_2,\overline{y_2}}).
$$
Since, when $X$ and $Y$ are bi-Hilbertian, $F_X$ and $F_Y$ are faithful
maps (see Cor. 2.11), the two seminorms have the same vectors of
length
zero (therefore $U$ is an injective map). Furthermore the
two norms 
$\|F_{\overline{Y}}(T^*T)\|^{1/2}$
and $\|F_X(TT^*)\|^{1/2}$ 
on $\K(\overline{Y}_B,X_B)$ 
are both 
equivalent to the operator norm,  still by
Cor. 2.11, and therefore they are equivalent. We have thus shown
that
$X\otimes^r_BY$ and $X\otimes^\ell_B Y$ are isomorphic as Banach spaces.
It is
now straightforward
to check that right and left actions are adjointable, and therefore
$X\otimes_B Y$ is bi-Hilbertian. Since $U$ is a bijective map
which preserves both inner products, it extends to a bijective $A$-$C$
bimodule
map $U: X\otimes_B Y\to\K(\overline{Y}_B, X_B)$ still preserving the inner
products, and the proof is now complete.
 \medskip

\noindent{\bf 2.3 Nondegeneracy of the left action}
\bigskip

\noindent The following nondegeneracy property will be relevant for our
purposes.\medskip

\noindent{\bf 2.14 Definition} The 
left action $\phi$ of a $C^*$--algebra $A$ on a right Hilbert
$C^*$--module $X_B$ will be called {\it nondegenerate\/} if $AX$ is
total in $X$.
\medskip

We recall the following characterization of 
nondegeneracy, due essentially to Vallin \cite{V}, see
also Prop. 2.5 in \cite{L}.
\medskip

\noindent{\bf 2.15  Proposition} {\sl Let $A$ and $B$ be $C^*$--algebras and
$X$ a right Hilbert $B$--module. For a $^*$--homomorphism
$\phi:A\to\L(X_B)$ the following conditions are equivalent.
\roster
\item $\phi$ is nondegenerate, 
\item $\phi$ is the restriction to $A$ of a unital $^*$--homomorphism
$\hat{\phi}: M(A)\to\L(X_B)$, strictly continuous on the unit
ball,
\item for some approximate unit $(u_\alpha)_\alpha$ of $A$,
$(\phi(u_\alpha))_\alpha$ converges strictly to the identity map
on $X$.
\endroster}\medskip

Note that if $\phi$ is  nondegenerate, 
 $(3)$ must hold for all approximate
units of $A$.
We  show that the left action of  a bi-Hilbertian $C^*$--bimodule is
automatically nondegenerate.\medskip

\noindent{\bf 2.16 Proposition} {\sl Let $_AX_B$ be a bi-Hilbertian $A$--$B$
$C^*$--bimodule. Then
the left (right) action of $A$ ($B$) on the underlying right 
Hilbert $C^*$--module $X$ is nondegenerate. }\medskip

\noindent{\it Proof}  
If $\{u_\alpha\}$ is an approximate unit of the closed ideal
of $A$ generated by the left inner products, $u_\alpha x$ converges
to $x$ for all $x\in X$, in the norm arising from the left inner product. 
Therefore
$AX$ is total in $X$ with respect to the norm defined by the left inner
product. Since
the two
norms
on $X$ defined by the right and left inner product are equivalent,
we also have that $AX$ is total with respect to the norm arising from the
right inner
product.
\medskip

\noindent{\bf 2.4 The index element and the Jones basic
construction}
\bigskip

\noindent If $X$ is bi-Hilbertian and of finite right numerical index, one
can
extend
the maps $\phi: A\to\L(X_B)$, $F:\K(X_B)\to A$ uniquely
 to  normal positive maps $\phi'': A''\to\K(X_B)''$, $F'':\K(X_B)''\to
A''$ between the
corresponding enveloping von Neumann algebras. Since $\phi$ is nondegenerate,
and the inclusion $M(A)\subset A''$ is unital,
$\phi''$ is unital homomorphism. The same does not hold for $F''$: 
 $F''(I)$
is, in general,  neither  the identity, nor invertible. 
\medskip 

\noindent{\bf 2.17 Definition} If ${}_AX_B$ is of finite right numerical index,
the {\it  right index element of } ${}_AX_B$, denoted, 
$r-\text{Ind}[X]$ is the element $F''(I)$ of $A''$.

If in particular ${}_BX_A$ is the bimodule arising from a conditional
expectation $E:B\to A$ as in Prop. 2.12, the corresponding right index
element will be denoted by $\text{Ind}[E]$. 
(We will give in Cor. 4.9 an alternative definition of
$\text{Ind}[E]$.)

If ${}_AX_B$ is of finite left numerical index, the left index element
of $X$ is, of course,
 $\ell-\text{Ind}[X]:=r-\text{Ind}[\overline{X}]$.

Notice that the numerical indices and the index
elements are
related by $$\|r-\text{Ind}[X]\|=r-I[X],$$
$$\|\ell-\text{Ind}[X]\|=\ell-I[X].$$

If one of $r-\text{Ind}[X]$ and
$\ell-\text{Ind}[X]$
is a scalar, or if $A=B$, the index element  of $X$ is 
 $\text{Ind}[X]:=
(r-\text{Ind}[X])(\ell-\text{Ind}[X])$.\medskip

Our next aim is to define an index element $\text{Ind}[X]$ 
in the general case.
We notice that for $c\in Z(B)$, the map
$\psi(c):
x\in X\to xc\in X$ has the map $x\in X\to xc^*\in X$ as an adjoint with
respect to the right inner product of $X$. Furthermore
$\psi(c)$ commutes obviously with all the alements of $\L(X_B)$, therefore
$\psi(Z(B))\subset
Z(\L(X_B))$. On the other hand it is not difficult to see that
$\psi(Z(B))=Z(\L(X_B))$.
We  need to consider an extension  of this
right
action of $Z(B)$ on $X$ to
the centre of $B''$. Therefore we anticipate the following lemma. 
\medskip

\noindent{\bf 2.18 Lemma} {\sl Let $X$ be a right Hilbert $B$--module,
and let $\psi: Z(B)\to Z(\L(X_B))$ denote the right action of $Z(B)$
on $X$. Then there is a canonical extension of $\psi$ to a
unital surjective
$^*$--homomorphism $\psi_0: Z(B'')\to Z(\K(X_B)'')$ with
$\text{ker }\psi_0=(1-q)Z(B'')$,
 where $q$ the central projection
of $B''$ corresponding to the weak
closure in $B''$ of the ideal generated by right 
inner products.
}\medskip

\noindent{\it Proof}
 Let $\pi$ denote a Hilbert space representation of $B$
on  $H_\pi$.
Consider the Stinespring induced representation $\tilde{\pi}$ 
of $\K(X_B)$ on the Hilbert space
$K_\pi:=X_B\otimes_B H_\pi$, defined by 
$T\to T\otimes 1_{H_\pi}$. 
Since $B$ and $\K(X_B)$ are strongly Morita equivalent, it is well
known that the map $\pi\to \tilde{\pi}$ is a bijective correspondence 
between representations of $B$ and representations of $\K(X_B)$.
Therefore the representation
$\rho=\oplus_{\pi}{\tilde{\pi}}'':\K(X_B)''\to\L(\oplus_\pi
K_\pi)$ is faithful and normal. It follows that
$\rho(\K(X_B)'')=\rho(K(X_B))''$.  On the other hand 
$Z(B'')$ acts on each $K_\pi$, and therefore on their direct sum,
by
${\psi'}_0: c\in Z(B'')\to 1_X\otimes
\pi''(c)\in\L(K_\pi)$. Clearly, ${\psi'}_0(c)\in\rho(\K(X_B))'$. 
If $A$ is a bounded operator on $\oplus_\pi K_\pi$
commuting
elementwise with $\tilde{\pi}(\K(X_B))$ then $A(x(y|z)_B\otimes\xi)=
\theta^r_{x,y}\otimes 1_{\oplus_\pi H_\pi}{A(z\otimes\xi)}$ for all
$x,y,z\in X$, $\xi\in H_\pi$, $\pi\in\text{Rep}(B)$.
 Choosing an approximate unit of the closed ideal of
right inner products of $B$ constituted by finite sums of elements
of the form $(y|y)_B$, we see that $A$ is of the form $1_X\otimes a$,
with $a\in \L(\oplus_\pi H_\pi)$.  
Approximating now $\pi''(c)$ strongly with a norm bounded net in $\pi(B)$,
shows that $T$ and ${\psi'}_0(c)$ commute. Therefore
${\psi'}_0(Z(B''))\subset
Z(\rho(\K(X_B)''))$. Notice that, ${\psi'}_0(c)=0$ if and only if
each $\pi''(c)$ annihilates the  subspace 
$qH_\pi$, or, equivalently $\pi''(cq)=0$ for all $\pi$, i.e.
$cq=0$. On the other hand if $A\in Z(\rho(\K(X_B)''))$ and we write 
$A=1_X\otimes a$, a one moment thought shows that $a\in Z(B'')$. The
$^*$--homomorphism
$\psi_0:=\rho^{-1}{\psi'}_0$ 
is then the desired extension of $\psi$.
\medskip

\noindent{\bf 2.19 Proposition} {\sl Let $X$ be of finite right numerical
index.
Then for any generalized right basis 
$\mu\to u_\mu$ of $X$,  the net $\mu\to\sum_{y\in
u_\mu}{}_A(y|y)$ is increasing and it converges
strongly in  $A''$ to $r-\text{Ind}[X]$.
This
limit is therefore independent on the choice 
of the  basis, and  belongs to the centre of $A''$.
If in addition $X$ is bi-Hilbertian one has
$\lambda r-\text{Ind}[X]\geq p$ where $p$ is the support projection of
$r-\text{Ind}[X]$ in $A''$, and
$\lambda$ is the best constant for which
$\lambda\|{}_A(x|x)\|\geq\|(x|x)_B\|, x\in X.$
Furthermore, if $\I$ denotes the weak closure in  $A''$
 of the span 
of left inner products ${}_A(x|y), x,y\in X$, one has
\roster
\item
$\text{ker }\phi''=(I-p)A'',$
\item
$\I=pA'',$
\item
the range of $F'':\K(X)''\to A''$ is $pA''$,
\item
if $z'$ denotes the inverse of $r-\text{Ind}[X]$ in $p A''$,
and $E'':=z' F'': \K(X_B)''\to pA''$. Then 
 $\phi''\circ E'': \K(X_B)''\to \phi''(A'')$ is a conditional
expectation with range $\phi''(A'')$ satisfying
$$\lambda\phi''(r-\text{Ind}[X]E''(T))\geq T, \quad T\in{\K(X_B)''}^{+}.$$
\endroster}\medskip

\noindent{\it Proof} The strong limit $z$ defined as in the statement
is independent of the generalized basis since it coincides with $F''(I)$.
Since $F''$ is still $A$-bilinear, $F(\phi(a))=F''(I)a=aF''(I)$ for all
$a\in A$, which shows that $F''(I)$ lies in the centre of $A''$.
Let us assume $X$ bi-Hilbertian, and let $\lambda$ be as in the
statement. 
The estimate 
$$\|T\|\leq \lambda\|F''(T)\|\leq r-I[X]\|T\|, \quad T\in
{\K(X)''}^+,$$
still holds, therefore if $T=\phi''(a)$, with $a\in {A''}^+$, 
$$\|\phi''(a)\|\leq\lambda\|za\|\leq r-I[X] \|\phi''(a)\|.$$
If we consider the restriction $\phi_0$ of $\phi''$ to the centre of
$A''$, we deduce that $\text{ker }\phi_0$ coincides with the weakly closed 
ideal generated by $I-p$. Therefore for a positive central element $a$ of
$A''$,
$\|\phi_0(a)\|=\|pa\|\leq\lambda\|za\|$.
Let us identify the centre of $A''$ with some $L^\infty(\Omega,\nu)$. We
claim that 
 for 
every $\epsilon>0$, the function $z-(\lambda^{-1}-\epsilon) p$ can not
take
negative
values on a measurable subset of $p\Omega$ with positive measure. Indeed,
if
$Y\subset\Omega$ where such a set, we would have, for some
$\epsilon<\lambda^{-1}$,
$$\lambda^{-1}-\epsilon=(\lambda^{-1}-\epsilon)\|\xi_Y\|\geq\|z\xi_Y\|\geq$$
$$\lambda^{-1}\|p\xi_Y\|=\lambda^{-1},$$
where $\xi_Y$ is the characteristic function of $Y$.
Therefore $z\geq\lambda^{-1} p$. 
Now if $a\in A''$, $\phi''(a)=0$ if and only if $\phi''(aa^*)=0$ and this
holds if and only if $paa^*=0$, i.e. $pa=0$, so $\text{ker
}\phi''=(I-p)A''$, and (1) is proved. Let $\I$ be the waekly closed ideal
of $A''$ defined as in the statement. Since the range of $F''$ is
contained in $\I$ and since $F''(\phi''(I))=z$, $p$ must belong to $\I$,
and therefore $pA''\subset\I$.
Conversely,
there exists an increasing, norm bounded net 
 $\alpha\to\sum_{w\in\alpha} {}_A(w|w)$,
indexed by the set of finite subsets of $X$,
which is a bounded approximate unit of the norm closed ideal generated by
the
left
inner products. Its weak limit, say $q$, is the unit of $\I$.
By Proposition 2.4 the net
$\alpha\to
\sum_{w\in\alpha}\theta^r_{w,w}$ is norm bounded.   
We have $\sum_{w\in\alpha} \phi(a)\theta^r_{w,w}\phi(a)^* \le
\lambda_0 \phi(aa^*)$
for some $\lambda_0>0$ and for all $a\in A$.  
Applying $F$ we obtain  
$$
\sum_{w\in\alpha} a{}_A(w|w)a^* \le \lambda_0 aza^* \ \
.
$$  
Thus ${\lambda_0}^{-1} aqa^* \le  aza^*$.  It follows that 
${\lambda}_0^{-1} q \le z$, hence $q\leq p$. Therefore $\I\subset pA''$,
and the proof of (2) is complete. (3) We are left to show that $pA''$
is contained in the range of $F''$. Let $z'$ be the inverse of 
$r-\text{Ind}[X]$ in $p A''$. For $a\in A''$,
$F''(\phi''(z'a))=r-\text{Ind}[X]z'a=pa$.
(4) It is now clear that $\phi'' E''$ is a conditional expectation
with range $\phi''(p A'')=\phi''(A)$. 
(5) Since $\lambda\|F''(T)\|\geq  \|T\|$ for $T$ positive 
is $\K(X_B)''$ and since $\phi''$ is isometric on $\I$, and therefore
on the range of $F''$, we see that
$\lambda\|\phi''(r-\text{Ind}[X]E''(T))\|=\lambda\|\phi''F''(T)\|=\lambda
\|F''(T)\|\geq\|T\|$. Therefore $\lambda\phi''(r-\text{Ind}[X]E''(T))\geq
T$ (cf.  \cite{FK}).
\medskip

\noindent{\bf 2.20 Corollary} {\sl Let ${}_AX_B$ be a bi-Hilbertian
$C^*$--bimodule of finite right numerical index.
Then the following properties are equivalent.
\roster
\item
$r-\text{Ind}[X]$
is invertible,
\item
$\phi''$ is faithful,
\item
the linear span  of the left inner products is weakly dense in $A''$.
\endroster}\medskip

\noindent{\it Remark} Notice that if ${}_BX_A$ is the bimodule arising
from a Pimsner--Popa conditional expectation, as in Prop. 2.12, 
$r-\text{Ind}[X]=\text{Ind}[E]$
must be invertible since the left inner product is full.
\medskip

\noindent{\bf 2.21 Definition} 
 Notice that
$\phi''(r-\text{Ind}[X])$ is invertible
in $Z(\phi(A)'\cap\K(X_B)'')$ and
$\psi_0(\ell-\text{Ind}[X])$ is  invertible in $Z(\K(X_B)'')$.
Therefore we define the index element of $X$ as 
an element of $\K(X_B)''$, in fact central in 
$\phi(A)'\cap\K(X_B)''$,  by
$\text{Ind}[X]:=\psi_0(\ell-\text{Ind}[X])\phi''(r-\text{Ind}[X])$.

\bigskip

\noindent{\bf 2.5 On the condition $r-\text{Ind}[X]\in M(A)$ and
existence
of finite bases}\bigskip

\noindent Under which conditions the index element 
lies in $M(A)$? We first discuss this question in the case where
$A$ is commutative, and derive a general criterion afterwards.
\medskip

Let $A=C_0(\Omega)$ be a commutative
$C^*$--algebra, and
let ${}_AX_B$ be a bi-Hilbertian $C^*$--bimodule with finite right
numerical index. We have seen in Prop. 2.12 that such bimodules arise,
e.g., from conditional expectations $E: A\to B$ satisfying a Pimsner-Popa
inequality.
The  right index element, 
$r-\text{Ind}[X]$,
 lies in the
commutative von Neumann algebra $A''$. However,
being the strong limit in the universal representation of $A$ of a net of
continuous, vanishing at infinity, functions over $\Omega$,  
$r-\text{Ind}[X]$
is the class function of a bounded, positive, lower
semicontinuous function on $\Omega$. 
By Dini's theorem,
$r-\text{Ind}[X]$
is continuos (or, in other words, lies in the multiplier algebra $M(A)$)
if and only if the
approximating net is uniformly convergent
on compact subsets of $\Omega$, that is, if and only if
 this net is
strictly convergent in $M(A)$.
\medskip

The key idea,  in the case
 where $A$ is not commutative, is to  replace the spectrum of $A$ with
the quasi-state space $Q$ of $A$. By Kadison's
function representation theorem (see, e.g., \cite{P}), the real Banach
space
$A_{\text{sa}}$ identifies isometrically with the real Banach space of 
continuous, vanishing at $0$, affine functions on $Q$. Now 
an analysis 
similar to the commutative situation leads to the following
characterization of the property  $r-\text{Ind}[X]\in M(A)$.\medskip

\noindent{\bf 2.22 Theorem} {\sl Let ${}_AX_B$ be a bi-Hilbertian
$C^*$--bimodule of finite right numerical index. Then the following
properties are equivalent:
\roster
\item $r-\text{Ind}[X]\in M(A)$ (and hence it is a central element of
$M(A)$),
\item there is a 
generalized right basis $\mu\to u_\mu\subset X$ of $X$ such that
the net $\mu\to\sum_{y\in u_\mu}{}_A(y|y)$ is strictly convergent in $A$,
\item for any  
generalized right basis $\mu\to u_\mu\subset X$ of $X$
the net $\mu\to\sum_{y\in u_\mu}{}_A(y|y)$ is strictly convergent
in $A$,
\item
the range $\phi(A)$ of the left action is included in $\K(X_B)$.
\endroster
If one of these conditions is satisfied, 
$r-\text{Ind}[X]=\lim_\mu
\sum_{y\in u_\mu}{}_A(y|y)$ in the strict topology of $A$.}\medskip

\noindent{\it Proof}
(3)$\Rightarrow$(2) is obvious.
(4)$\Rightarrow$(3)
Let $\mu\to u_\mu$ be a generalized basis of $X$. Then 
$\sum_{y\in
u_\mu}{}_A(y|y)$ 
converges strictly if and only if
for all $a\in A$,
$$\sum_{y\in
u_\mu}a{}_A(y|y)=\sum_{y\in
u_\mu}{}_A(\phi(a)y|y)=F(\phi(a)\sum_{y\in u_\mu}\theta_{y,y})$$
converges in norm. Therefore if (4) holds, convergence of the above net
follows from norm continuity of $F$ and the fact that $\sum_{y\in
u_\mu}\theta^r_{y,y}$ is an approximate unit of $\K(X_B)$.
(2)$\Rightarrow$(4)
 By Cor. 2.11, for any
$T\in\K(X_B)^+$,
$\lambda\|T\|\leq \|F(T)\|$, where
$\lambda$ is the best positive constant for which
$\lambda\|(x|x)_B\|\leq\|{}_A(x|x)\|$. If now (2) holds
for some generalized basis $\mu\to u_\mu\subset X$, $\phi(a)^*\sum_{y\in
u_\mu}{}_A(y|y)\phi(a)=F(
\phi(a)^*(
\sum_{y\in
u_\mu}\theta^r_{y,y})\phi(a)
)$ is an increasing, norm converging net, and
therefore 
the net $\phi(a)^*(
\sum_{y\in
u_\mu}\theta^r_{y,y})\phi(a)$ is increasing and norm converging in
$\K(X_B)$.  It follows that $\phi(a^*a)\in\K(X_B)$
for all $a\in A$.
(2) $\Rightarrow$ (1) is obvious, since $r-\text{Ind}[X]$ is  the 
strict limit of a strictly convergent net.
(1)$\Rightarrow$ (3) By Kadison's function representation (see, e.g.,
\cite{P} the
selfadjoint part of $A''$
identifies isometrically, as a real Banach space,
 with the real Banach space $B_0(Q)$ of affine, bounded
functions
on the quasi-state space $Q$ of $A$ vanishing in $0$. 
Under this
identification,
the  selfadjoint elements of $A$ correspond to the continuous
functions. If $r-\text{Ind}(A)\in M(A)$, for all $a\in A$,
$a^*((r-\text{Ind}[X])-\sum_{y\in u_\mu}{}_A(y|y))a\in A$.
On the other hand the net
$(r-\text{Ind}[X])-\sum_{y\in u_\mu}{}_A(y|y)$ decreases weakly to
$0$ in $A''$,
therefore for any $\phi\in Q$,
$$\phi(a^*((r-\text{Ind}[X])-\sum_{y\in u_\mu}{}_A(y|y))a)$$
decreases to $0$. 
By Dini's theorem, this net
 converges uniformly to 0 on $Q$, and therefore
$\|(r-\text{Ind}[X]-\sum_{y\in u_\mu}{}_A(y|y))^{1/2}a\|^2\to 0$,
which implies
$\|(r-\text{Ind}[X]-\sum_{y\in u_\mu}{}_A(y|y))a\|^2\to 0$
as the net $\sum_{y\in u_\mu}{}_A(y|y)$ is norm bounded.

Assume now that one of these equivalent properties holds, and let
$z:=\lim_\mu\sum_{y\in u_\mu}{}_A(y|y)\in M(A)$ for some generalized right
basis
$\mu\to u_\mu$ of $X$.
Since, for all $a\in A$, $za=F(\phi(a))$, $z$ is independent of the choice
of the basis. 
\medskip

\noindent{\bf 2.23 Definition}
A bi-Hilbertian $A$--$B$ $C^*$--bimodule $X$
  will be called  {\it of finite right index}
 if
\roster 
 \item $X$ is of finite right  numerical index, 
 \item $r-\text{Ind}[X]\in M(A)$ (and hence $r-\text{Ind}[X]\in Z(M(A))$).
\endroster

\noindent{\it Remark} Notice that property (2) above can be replaced by
any
of the equivalent conditions in Theorem 2.22. 
\medskip

Similarly, $X$ is {\it of finite left index} if the contragradient
bimodule
${}_B\overline{X}_A$ is of finite right index.

${}_AX_B$ will be called {\it of finite index} if it is of finite right
as well as left indices. \medskip

We study the special case where the $C^*$--algebras
are $\sigma$--unital or unital. \medskip

\noindent{\bf 2.24 Corollary}
{\sl Let  ${}_AX_B$ be a bi-Hilbertian $C^*$--bimodule of finite right
index.
 \roster
\item
If $A$  is
$\sigma$--unital,  $X$ is countably
generated as a right  Hilbert module,
\item 
if  $A$ and $B$ are $\sigma$--unital, 
   $X_B$ admits
an unconditionally convergent countable right basis $\{u_i\}_{i\in{\bN}}$,
 therefore
$$r-\text{Ind}[X]=\sum_{i\in{\Bbb N}}{}_A(u_i|u_i)$$ in the strict
topology
of $A$
and
$F(T)=\sum_{i\in{\Bbb N}}{}_A(Tu_i|u_i), T\in\K(X_B)$ in norm.
\endroster}\medskip

\noindent{\it Proof} (1)
Let  $(u_i)_{i\in{\Bbb N}}$ be
a countable
approximate unit of  $A$. By nondegeneracy of the left
action, and the fact that the left action
has range in $\K(X_B)$,  $\phi(u_i)\in\K(X_B)$
is a countable approximate unit for $\K(X_B)$, so $\K(X_B)$
is $\sigma$--unital and this shows that $X$ is countably generated as
a right Hilbert $B$--module (see, e.g., \cite{B}). (2) If in addition
$B$ is $\sigma$-unital,
$X_B$ admits a countable right basis by Lemma 1.6, therefore the formulas
for $r-\text{Ind }[X]$ and for $F$ follow. 
 \medskip

\noindent{\bf 2.25 Corollary} {\sl Let ${}_AX_B$ be a bi-Hilbertian
bimodule 
of finite right numerical index, and let $A$ be a unital $C^*$--algebra.
The following are equivalent:
\roster 
\item
$X$ admits a finite right basis,
\item $r-\text{Ind}[X]\in A$.
\endroster}\medskip

\noindent{\it Proof} (1)$\Rightarrow$ (2) follows from the definition of 
$r-\text{Ind}[X]$.
Conversely, assume
that  (2)
holds. By nondegenereracy of the left action,
$\phi(I)$
must be the identity map on $X$. Since, by Theorem 2.22, the range of
$\phi$
is included in the compacts, $I\in\K(X_B)$, so $X_B$ admits a finite basis.
\medskip

\noindent{\bf 2.26 Corollary} {\sl Let ${}_AX_B$ be a bi-Hilbertian
$C^*$--bimodule with finite right numerical index. If $A$ is simple
then $r-\text{Ind}[X]$ is a scalar, and therefore
$X$ is of finite right index.}\medskip

\noindent{\it Proof} 
Since  $A$ is a simple $C^*$--algebra, the only positive elements 
of
$Z(A'')$ arising as strong limits of increasing nets in $A$
must be scalar (see, e.g., Lemma 3.1 in \cite{I}), so
$r-\text{Ind}[X]$ is a scalar, and therefore it belongs to $M(A)$.
\medskip

The following result has been obtained by Izumi \cite{I} in the case
where
$B$ is a simple $C^*$--algebra.\medskip

\noindent{\bf 2.27 Corollary} 
{\sl Let $A\subset B$ be an inclusion of unital $C^*$--algebras,
and let $E:B\to A$ be a conditional expectation satisfying a Pimsner-Popa
inequality. 
Then $E$ admits a finite quasi--basis in the sense of \cite{W} if and
only if $\text{Ind}[E]\in
B$.}\medskip

\noindent{\bf 2.6 The Jones basic construction}\bigskip

\noindent{\bf 2.28 Proposition} {\sl
Let ${}_AX_B$ be a bi-Hilbertian $C^*$--bimodule
of finite right index. Then 
\roster
\item 
the map $F:\K(X_B)\to A$ extends uniquely
to a strictly continuous map $\hat{F}:\L(X_B)\to M(A)$. One has
$\hat{F}(I)=r-\text{Ind}[X]$, $\|\hat{F}\|=r-I[X]$. $\hat{F}$ is still
positive,
$M(A)$--bilinear and  satisfies 
$$\lambda'T\leq \hat{\phi}\hat{F}(T),\quad T\in\L(X_B)^+,$$
where $\lambda'$ is the best constant for which
$\lambda'\|(x|x)_B\|\leq\|{}_A(x|x)\|$, $x\in X$,
\item
the
support projection $p$ of $r-\text{Ind}[X]$ in $A''$ lies in fact the
centre of
$M(A)$ and satisfies $r-\text{Ind}[X]\geq\lambda'p$,
\item
$\text{ker }\phi=(I-p)A$, $\text{ker }\hat{\phi}=(I-p)M(A)$,
\item the norm closed subspace of $A$
generated by the left inner products coincides with $pA$,
\item the range of $F:\K(X_B)\to A$ is $pA$.
\endroster}\medskip

\noindent{\it Proof} (1) Let us restrict the map $F'':\K(X_B)''\to
A''$
to a positive map $\hat{F}:\L(X_B)\to A''$. If $T\in\L(X_B)$, and $a\in
A$, both $\phi(a)T$ and $T\phi(a)$ lie in $\K(X_B)$, hence
$a\hat{F}(T)=F(\phi(a)T)\in A$, $\hat{F}(T)a=F(T\phi(a))\in A.$ Therefore
$\hat{F}(T)\in
M(A)$. $\hat{F}$ is clearly a strictly continuous extension of $F$, and it
is uniquely determined by  this property. In particular, 
$\hat{F}(I)=r-\text{Ind }[X]$.
The remaining properties follow from the corresponding properties of
$F''$.
(2) Let us restrict $\hat{F}$ to the image under $\hat{\phi}$ of the
centre
of $M(A)$. One has $\hat{F}(\hat{\phi}(a))=r-\text{Ind}[X]a$ for all
$a\in Z(M(A))$.
By the
inequality
in (1), $$\text{ker }\hat{\phi}\upharpoonright_{Z(M(A))}=\{a\in Z(M(A)):
a(r-\text{Ind}[X])=0\}.$$ In particular, regarding $Z(M(A))$ as the
algebra
of continuous function over its spectrum,
for all $a\in Z(M(A))$, the norm of 
$\hat{\phi}(a)$ coincides with the norm of the restriction
of $a$ on the support $K$ of $r-\text{Ind}[X]$.
 If $\xi$ is an element  of the spectrum of $Z(M(A))$
such that
$(r-\text{Ind}[X])(\xi)\neq0$, and $\epsilon>0$, we can find an open
set  $U$ containing 
$\xi$ such that $(r-\text{Ind}[X])(\xi')<(r-\text{Ind}[X])(\xi)+\epsilon$
for all $\xi'\in U$. Let $a$ be a continuous function on the spectrum of 
$Z(M(A))$ with support in $U$, such that $0\leq a\leq 1$ and taking value 
$1$ on a compact set containing  $\xi$. We then have
$$(r-\text{Ind}[X])(\xi)+\epsilon>
\|r-\text{Ind}[X]a\|\geq\lambda'\|\hat{\phi}(a)\|=
\lambda'\|a\upharpoonright_{K}\|=\lambda',$$
so
$$r-\text{Ind}[X](\xi)\geq\lambda'.$$ Therefore the support of
$r-\text{Ind}[X]$ is an open and closed subset of the spectrum 
of $Z(M(A))$, which implies that its characteristic function $p$
belongs to $Z(M(A))$, and $r-\text{Ind}[X]\geq\lambda'p$.
(3) $\text{ker }\phi=\text{ker }{\phi''}\cap A=(I-p)A''\cap A=(I-p)A$ 
by Prop. 2.19. Smilarly, $\text{ker }\hat{\phi}=(I-p)M(A)$.
(4) Let $\J$ be the norm closed subspace generated by the left inner
products,
with weak closure $\I$ in $A''$.
We have $\J\subset \I\cap A=pA''\cap A=pA$. Let $z'$ be the inverse of 
$p(r-\text{Ind}[X])$ in $pZ(M(A))$  regarded as an element of $Z(M(A))$.
Then for all $a\in A$, $pa=r-\text{Ind}[X]z'a=F(\phi(z'a))\in \J$
since the range of $F$ is contained in $\J$. Thus $pA=\J$. The last
property is now clear.
\medskip

\noindent{\bf 2.29 Corollary} {\sl If $X$ is a bi-Hilbertian $A$-$B$
$C^*$--bimodule of finite right index, 
 the following properties are equivalent.
\roster
\item The left inner product is full,
\item $r-\text{Ind }[X]$ is invertible,
\item $\phi$ is faithful.
\endroster}\medskip

We next construct the analogue of the Jones basic construction
in the $C^*$--algebra setting.\medskip 

\noindent{\bf 2.30 Corollary} {\sl Let ${}_AX_B$ be a bi-Hilbertian
$C^*$--bimodule
of finite right index. Consider
the positive $A$-bilinear map $E:T\in\K(X_B)\to z'F(T)\in pA$,
with $z'$ the inverse of $(r-\text{Ind}[X])p$ in $pZ(M(A))$.
Then $\phi E:\K(X_B)\to \phi(A)$ is a conditional expectation
with range $\phi(A)$
which satisfies 
$${\lambda}\phi (r-\text{Ind}[X]E(T))\geq T,\quad
T\in\K(X_B)^+$$
where $\lambda$ is the best constant for which 
$\lambda\|{}_A(x|x)\|\geq\|(x|x)_B\|$.}\medskip

\noindent{\it Proof} $\phi E$ is clearly a positive, $\phi(A)$--bilinear map
with range $\phi(pA)=\phi(A)$. For all $a\in A$,
$\phi E(\phi(a))=\phi(z'za)=\phi(pa)=\phi(a)$ by (3) of Prop. 2.28,
hence
$\phi E$ is a conditional expectation.
The remaining inequality follows from Cor. 2.11 and $\phi\circ
F=\phi(r-\text{Ind}[X])\phi\circ E$.
\medskip

\noindent{\bf 2.7 Examples}
\bigskip

\noindent We conclude this section  with few   examples of finite
index
bimodules already known in
the literature. More examples will be discussed in section 6.

The first example arises from compact quantum groups,
or, equivalently, conjugation in finite dimensional Hilbert spaces,
\cite{Wo},
where the bi-Hilbertian structure is  usually 
described in terms of  antilinear invertible
mappings between Hilbert spaces implementing the 
conjugation structure.\medskip

\noindent{\bf  2.31 Example}
Let  $H = \comp ^n$ be a finite dimensional Hilbert 
space and 
 $T$ be a positive invertible linear map on $H$.
$H$ is an $\comp$--$\comp$ bimodule in the obvious way. We endow $H$ with a 
bi-Hilbertian bimodule structure by setting
 $(x|y)_{\comp} = \sum_{i} \overline{x(i)}y(i)$ and
${}_{\comp}(x|y) = \sum_{i} (Tx)(i) \overline{y(i)}=(y|Tx)_\comp$. These 
inner products
induce
equivalent norms on $H$ since $T$ is invertible,  making $H$
 into a bi-Hilbertian bimodule. Since $H_\comp$ and
${}_\comp H$ are  finite dimensional Hilbert spaces, 
  $\K(H_\comp)=\L(H_\comp)$ and $\K({}_\comp H)=\L({}_\comp H)$,
therefore
$H$ will be of finite index if  it is of finite left 
and right numerical indices. 
Let $\text{Tr}$ be the nonnormalized trace on $M_n(\comp) = \K(H)$. 
Since $\text{Tr}(T\theta^r_{x,y}) = {}_{\comp}(x|y)$, 
for $x_1,\dots,x_p\in H$,
$$\|\sum_1^p{}_\comp(x_j|x_j)\|\leq\text{Tr}(T)\|\sum_1^p
\theta^r_{x_j,x_j}\|,$$ therefore
$r-\text{Ind}[X]=\text{Tr}(T)$.
Since
$(x|y)_\comp={}_\comp(T^{-1}y|x)$, 
we deduce that
$\ell-\text{Ind}[X]=\text{Tr}(T^{-1})$.
\medskip

\noindent{\bf 2.32 Example} More generally, let $\Omega$ be a 
locally compact
Hausdorff space, $H=(\Omega, H(\omega)_{\omega\in\Omega}, \Gamma)$ a 
continuous field of Hilbert spaces, with $\Gamma$ the space of continuous 
sections of $H$. Let $$X=\{x\in\Gamma:\omega\to\|x(\omega)\|\in 
C_0(\Omega)\}$$ be the associated  right Hilbert bimodule over 
$C_0(\Omega)$. Let us consider a field $\omega\to 
T(\omega)\in\L(H(\omega))$ of positive, trace-class operators  on each 
$H(\omega)$ 
defining an element $T$ of $\L(X_{C_0(\Omega)})$ (e.g. $T\in FR(X)$).
We can then define a left inner product on $X$, continuous with respect to 
the right one,  by
$${}_{C_0(\Omega)}(x|y)(\omega):=(y(\omega)|T(\omega)x(\omega)).$$
Writing the left inner product in the form
$${}_{C_0(\Omega)}(x|y)(\omega)=\text{Tr}T(\omega)\theta_{x(\omega), 
y(\omega)}, $$  shows  that $X$ is of finite right numerical index  
if and only if $\sup_\omega \text{Tr}T(\omega)$ is finite. In this 
case, $r-\text{Ind}[X](\omega)=\text{Tr}T(\omega)$. Therefore 
$X$ is of finite index if and only if $\omega\to\text{Tr}T(\omega)$
is a bounded, continuous function on $\Omega$ (e.g. $T\in
FR(X)\cap\K(X)^+$).
Notice that the set of linear combinations of elements $T\in \K(X)^+$ for 
which $T(\omega)$ is trace-class 
and $\omega\in\Omega\to \text{Tr} T(\omega)$ is continuous, is a 
$^*$-ideal of $K(X)$ by
4.5.2 in \cite{Di}, norm dense in $\K(X)$. Assume from now on that 
$\sup_{\omega}\text{dim}H(\omega)$
is finite. Then $T=I$ defines a bi-Hilbertian bimodule of finite
right (and left) numerical
index, with $r-\text{Ind}[X](\omega)=\text{dim}H(\omega)$. However,
it is
not of finite index, unless the dimension function  is 
continuous. However, if $T\in\K(X)^+$, $\omega\to\text{Tr}T(\omega)$ is 
always bounded and continuous  by \cite{F}, and therefore $T$ does define
a 
finite right  index structure on $X$.
\medskip

The next example concernes with conditional expectations between 
unital $C^*$--algebras, and was introduced in \cite{W} as 
a $C^*$--algebraic analogue  of Jones index theory for finite
subfactors.\medskip

\noindent{\bf 2.33 Example}
Let $A \subset B$ be an inclusion of unital $C^*$--algebras and let
$E : B \rightarrow A$ be a conditional expectation of finite index in 
the sense of  \cite{W}. We thus have elements    
$\{u_1,...,u_n\}$ in $B$ such that
$ x = \sum_{i=1}^n u_iE({u_i}^*x)$  for any $x \in B$. Such
 elements were called a quasi-basis of $E$, and the index of $E$ was
defined
as $\text{Ind}[E]=\sum_i u_i{u_i}^*$
in \cite{W}. 
Consider $X =B$ 
as a 
$B$--$A$ bimodule in the obvious way, and define on
 the left $B$-valued inner product ${}_B(x|y) = xy^*$ 
and 
right $A$-valued inner product $(x|y)_A = E(x^*y)$. 
   Since $\theta^l_{x,y} = 
y^*x$, $\K({}_BX)=B$ acting on $X$ by right multiplication. Thus
  the right $A$-action has range in 
$\K({}_BX)$. 
Proposition 2.6.2 in \cite{W} shows that  $E$ satisfies the inequality  
 $E(x^*x) \ge \|\text{Ind} [E] \|^{-1}x^*x$.
 Therefore $X$ is bi-Hilbertian
and of finite left index:
$\ell-\text{Ind}[X]=I_A$.
\par Its contragradient bimodule is
 $Y = {}_AB_B$, as a  $A$--$B$ bimodule, with  
 $A$-valued inner product ${}_A(x|y) = E(xy^*)$ and 
 $B$-valued inner product $(x|y)_B = x^*y$.    $Y$ is a 
bi-Hilbertan
$B$--$A$ bimodule of finite left index.  In fact the algebra  
$\K({}_AY)$ is isomorphic to the 
$C^*$--basic construction $C^*\langle B,e_A\rangle$
which contains the image of the left action of $B$
via $b=\sum_i\theta^l_{ {u_i}^*b, {u_i}^*}$. Therefore $X$ is of finite
right index as well and
$r-\text{Ind}[X]=\ell-\text{Ind}[Y]=\text{Ind}[E]$.
 Consider the dual conditional expectation 
 $F: C^*\langle B,e_A \rangle \rightarrow B$, 
$F(xe_Ay) = (\text{Ind} [E])^{-1}xy$.  
 The Pimsner-Popa inequality for $F$ shows that $Y$ is bi-Hilbertian,
while the fact that $F$ is contractive gives
$$
 \| \sum_{i=1}^n{}(x_i|x_i)_B\|
 \le \|\text{Ind} [E] \| \| \sum_{i=1}^n \theta^l_{x_i,x_i}\|. 
$$
\par 
Consider now the bimodule
 $Z = {}_AB_A$ as a bi-Hilbertian $A$--$A$ Hilbert bimodule  with  
right $A$-valued inner product $(x|y)_A = E(x^*y)$ and 
left $A$-valued inner product ${}_A(x|y) = E(xy^*)$.  Then $Z$ is 
isomorphic to 
$Y\otimes_B X$  so $Z$ is a Hilbert 
$A$--$A$ bimodule of finite right index.
\medskip

The following example is a generalization 
of index theory to finitely generated Hilbert bimodules, studied
in \cite{KW1}.\medskip

\noindent{\bf 2.34 Example}
Let $A$ and $B$ be unital \cst--algebras  and let $X$ be a Hilbert 
$A$--$B$ 
bimodule 
such that both left and right actions are unit preserving.  Then 
 $X$ is a bi-Hilbertian bimodule of finite index if and only if 
$X$ is of finite type  in the 
sense of \cite{KW1}.  In fact, assume that $X$ is bi-Hilbertian
and of finite index, then  
 $X$ is necessarily finitely generated  projective  as a 
right module (or as a left module),
since 
$\K(X_B)$ (or $\K({}_AX)$) contains the identity map.
The two norms defined by the two  inner 
products of $X$ are equivalent, thus $X$ is of finite type.  Conversely, 
assume that 
$X$ is of finite type in the sense of \cite{KW1}. Then it is clear that
$X_B$ is bi-Hilbertian
and that the left 
$A$--action
on $X$ has range into $\K(X_B)=\L(X_B)$. Furthermore
$X_B$ is of finite right
numerical index by Lemma 1.26 in \cite{KW1}.  Thus $X$ is of
finite right index.  
Similarly,
$X$ is of finite left index and therefore of finite index.

We conclude this section with a discussion of a 
a Pimsner-Popa conditional expectation with no finite quasi--basis.
This  example was pointed out in \cite{W}. Later it was
considered also in  in \cite{FK}.
We show that this inclusion is determined by a natural $\sigma$--unital
 subinclusion   of finite index in the sense of Def. 2.23.
\medskip

\noindent{\bf 2.35 Example}
Consider the 
$C^*$--algebra $C([-1,1])$ of continuous functions over the interval
$[-1,1]$ and the $C^*$--subalgebra
$C([-1,1])_e=\{f\in C([-1,1]): f(-x)=f(x)\}$ of even functions. 
The conditional expectation $E: C([-1,1])\to C([-1,1])_e$
associating 
 to $f\in C([-1,1])$ the function $\frac{1}{2}(f(x)+f(-x))$ does not have a finite quasi--basis in the sense of 
\cite{W} since 
$C([-1,1,])$ is not a finite  projective module over $C([-1,1])_e$.  It follows that the bi--Hilbertian
bimodule ${}_{C([-1,1])}C([-1,1,])_{C([-1,1])_e}$ with inner products 
$${}_{C([-1,1])}(f|g)=f\overline{g},\quad 
(f|g)_{C([-1,1])_e}=E(\overline{f}g)$$ is not of finite right index in the
sense of Def. 2.23 because the identity operator over 
$C([-1,1])_{C([-1,1])_e}$
is not compact.
However, 
the Pimsner--Popa inequality
$E(f)\geq \frac{1}{2}f$ holds for any
$f\in C([-1,1])^+$.   One  can treat this example by our 
methods passing to a  subinclusion in the following way.
Consider the $\sigma$--unital $C^*$--subalgebra $C_0([-1,1])=\{f\in C([-1,1]): f(0)=0\}$.
 Then the restriction of $E$ still  defines a conditional
expectation $E: C_0([-1,1])\to C_0([-1,1])_e$, where
 $C_0([-1,1])_e=C([-1,1])_e\cap C_0([-1,1])$, and therefore a bi-Hilbertian  $C^*$--bimodule 
$$X={}_{C_0([-1,1])}C_0([-1,1])_{C_0([-1,1])_e}$$ which we  show to be
  of finite right index. We first  show that the left action $C_0([-1,1])$
has  range included in the compacts. Set, for $f\in C_0([-1,1])$,
$$f_e(x)=E(f)(x)=\frac{f(x)+f(-x)}{2}\in C_0([-1,1])_e,$$
$$f_o(x)=\frac{f(x)-f(-x)}{2}\in C_0([-1,1])_o\in \{f\in C_0([-1,1]): f(-x)=-f(x)\}.$$ 
Clearly $f=f_e+f_o$ and $E(\overline{f}_eg_o)=0$, $f,g\in C_0([-1,1])$. Therefore
 the  right Hilbert bimodule $X_B$ splits into the direct sum
of the subspaces of even and odd functions:
$X=X_e\oplus X_o$,  $X_e:=\{f\in X: f(-x)=f(x)\}$, $X_o=\{f\in X: f(-x)=-f(x)\}$. Similarly,
as a vector space, $C_0([-1,1])=C_0([-1,1])_e\oplus C_0([-1,1])_o$. For $f,g\in X=C_0([-1,1])$,
$\theta^r_{f,g}(h_e+h_o)=f\overline{g}_eh_e+f\overline{g}_oh_o$. Therefore if, for $n\in{\Bbb N}$, $u_n$
 is a positive continuous function in $ C_0([-1,1])_e$ such that $u_n(x)=1$
for $|x|\geq\frac{2}{n}$ and $u_n(x)=0$ for $|x|\leq\frac{1}{n}$, the
sequence
$\theta^r_{f, u_n} +\theta^r_{xf,\frac{u_n}{x}}$ is 
norm converging to the multiplication operator by
$f$. Therefore by Theorem 2.22, (2) of Def. 2.23 holds. We are left to
show that $X_B$ is
of finite right numerical index, and this follows from the Pimsner--Popa
inequality. 
\medskip

\heading
3. Continuous bundles of finite dimensional $C^*$--algebras\\
arising from  bimodules of finite right index
\endheading

Let $X$ be a right Hilbert $A$-$B$ bimodule
with nondegenerate left action $\phi$, and let
us consider the extension $\hat{\phi}: M(A)\to\L(X_B)$ of  
 $\phi$ to the multiplier algebra (see Prop. 2.15). 
Restricting $\hat{\phi}$ to the centre $Z(M(A))$ of $M(A)$ yields a
unital 
$^*$--homomorphism $\hat{\phi}:Z(M(A))\to\L(X_B)$, still denoted $\hat{\phi}$.
\medskip

\noindent{\bf 3.1 Proposition} {\sl If $X$ is a right Hilbert $A$-$B$ bimodule
with nondegenerate left action (this being the case if, e.g., $X$ 
is bi-Hilbertian, by Prop. 2.16),
the range of $\hat{\phi}:Z(M(A))\to\L(X_B)$
is actually included in the centre of ${}_A\L(X_B)$,
the algebra of right adjointable maps on $X_B$ commuting with the left
action. 
Therefore
${}_A\L(X_B)$ becomes a 
$Z(M(A))$-algebra in the sense of \cite{Ka}.}\medskip

Adopting a standard procedure we can represent ${}_A\L(X_B)$ as a
 {\it semicontinuous
field of $C^*$--algebras\/} $\omega\to\L_\omega$ over the spectrum 
$\Omega$ of $Z(M(A))$
 in the sense of
 \cite{Ka}. 
Let, for $\omega\in\Omega$, $J_\omega$ be 
the closed two-sided ideal of ${}_A\L(X_B)$  generated by 
the image under $\hat{\phi}$ of $C_\omega(\Omega)$, 
the continuous functions on $\Omega$
vanishing at $\omega$.
The fiber   at $\omega$ is the
quotient $C^*$--algebra
 $\L_\omega:={}_A\L(X_B)/J_\omega$. 
We will show that  the field $\omega\to\L_\omega$ is in fact continuous  
 in the case
where $X$ is 
bi-Hilbertian and of finite  right index (see Theorem 3.3).

Let ${}_AX_B$ be bi-Hilbertian and of finite right index. In Proposition
2.28 we have constructed a $M(A)$-bilinear, positive, 
strictly
continuous map $\hat{F}:\L(X_B)\to M(A)$  satisfying a Pimsner-Popa 
inequality and with range the ideal $pM(A)$, with $p$ the support
projection of $r-\text{Ind}[X]$. Restricting $\hat{F}$ to the
$C^*$--subalgebra
${}_A\L(X_B)$ yields a map, still denoted $\hat{F}$, with the same
properties, and with range  the ideal
 $pZ(M(A))$ of the commutative $C^*$--algebra $Z(M(A))=C(\Omega)$.  
We write  $\Omega=\Omega_0\cup\Omega_1$, with $\Omega_0$
corresponding to the  projection $p$ and 
$\Omega_1$ 
to  $I-p$. 
The map $\hat{F}$ makes
${}_A\L(X_B)$ into a right Hilbert $C(\Omega)$-module (in fact a
a Hilbert $C(\Omega_0)$-module) by
$(S|T)=\hat{F}(S^*T)$. Since $\hat{F}$ is norm continuous and satisfies a
Pimsner-Popa
inequality, the operator norm and the Hilbert module norm are equivalent,
therefore ${}_A\L(X_B)$ is complete in the Hilbert module norm.
 Since the inner product is evaluated on a commutative $C^*$--algebra, we
can 
represent ${}_A\L(X_B)$ as a {\it continuous field of Hilbert spaces}
over $\Omega$
in the sense of \cite{Di}. For each $\omega\in \Omega$, 
the fiber Hilbert space at
$\omega$ is given by $H_\omega={}_A\L(X_B)/M_\omega$, where $M_\omega$ 
is the norm closed subspace
of ${}_A\L(X_B)$, in the Hilbert module norm,
generated by ${}_A\L(X_B)\hat{\phi}(C_\omega(\Omega))$.
For each $\omega\in\Omega_1$, 
$M_\omega ={}_A\L(X_B)$ since $\hat{\phi}$ annihilates 
$(I-p)Z(M(A))=C(\Omega_1)$, therefore
$M_\omega=0$, as expected.
Since the $C^*$-algebra norm and the Hilbert
module norm are equivalent, $J_\omega=M_\omega$  as vector spaces, and they
 are isomorphic as Banach spaces. In particular, $\L_\omega=0$ for 
$\omega\in\Omega_1$.
Let $\pi_\omega: {}_A\L(X_B)\to \L_\omega$ and 
$p_\omega:{}_A\L(X_B)\to H_\omega$ denote 
the corresponding quotient maps in the $C^*$-algebraic and Banach space 
space
sense. 
\medskip

\noindent{\bf 3.2 Lemma} {\sl If ${}_AX_B$ is a bi-Hilbertian bimodule of 
finite right index, for all $T\in{}_A\L(X_B)$ and for all
$\omega\in\Omega_0$,
$${\lambda'}^{1/2}\|\pi_\omega(T)\|\leq\|p_\omega(T)\|\leq
(r-\text{Ind}[X])(\omega)^{1/2}\|\pi_\omega(T)\|,$$
where $\lambda'$ is the best positive scalar for which
$\|{}_A(x|x)\|\geq\lambda'\|(x|x)_B\|$, $x\in X$.}\medskip

\noindent{\it Proof}
The positive $M(A)$-bilinear map $\hat{F}:\L(X_B)\to M(A)$ satisfies
$ \hat{\phi}F(T)\geq \lambda'T$ for all $T\in\L(X_B)^+$, by Prop. 2.28. 
Therefore if 
$T\in{}_A\L(X_B)^+$, evaluating $\pi_\omega$
on this estimate yields $\pi_\omega(\hat{\phi}F(T))\geq
\lambda'\pi_\omega(T)$
which shows that
$$\|p_\omega(T)\|^2=
(p_\omega(T),p_\omega(T))=\hat{F}(T^*T)(\omega)=$$
$$|\pi_\omega(\hat{F}(T^*T))|\geq\lambda'|\pi_\omega(T^*T)|=
\lambda'\|\pi_\omega(T)\|^2.$$
Consider the map $G_\omega:\L_\omega\to {\Bbb C}=
C(\Omega)/C_\omega(\Omega)$ associating
$\hat{F}(T)(\omega)$ to $\pi_\omega(T)$. This map is well defined:
$\pi_\omega(T)=0$
implies that $T\in J_\omega$,  therefore $\hat{F}(T)$ belongs to 
$F(J_\omega)$ which is contained
in  the closed linear span
of $C_\omega(\Omega) F({}_A\L(X_B))$ in the $C^*$--algebra norm. Clearly 
the
latter
space is contained in $C_\omega(\Omega)$.
Now $G_\omega$ is a positive functional on 
the  $C^*$-algebra $\L_\omega$ taking the unit of $\L_\omega$ to 
$(r-\text{Ind}[X])(\omega)$, 
and therefore  $\|G_\omega\|=(r-\text{Ind}[X])(\omega)$.
Thus for
 all $T\in {}_A\L(X_B)$,
$$\|p_\omega(T)\|^2=|\hat{F}(T^*T)(\omega)|=
\|G_\omega(\pi_\omega(T^*T))\|\leq$$
$$(r-\text{Ind}[X])(\omega)\|\pi_\omega(T^*T)\|=(r-\text{Ind}[X])(\omega)
\|\pi_\omega(T)\|^2.$$\medskip

We are now ready to show the following result.\medskip

\noindent{\bf 3.3 Theorem} {\sl Let
 $X$ be a bi-Hilbertian $A$-$B$ $C^*$--bimodule of
finite right 
index, and let $\Omega$ be the spectrum of $Z(M(A))$. Then for each 
$\omega\in\Omega$, the quotient $C^*$-algebra
 $\L_\omega$ is  finite dimensional, and  
 $$\text{dim}(\L_\omega)\leq [{\lambda'}^{-1} 
(r-\text{Ind}[X])(\omega)]^2,$$
where $\lambda'$ is the best constant for which 
$\|{}_A(x|x)\|\geq\lambda'\|(x|x)_B\|$ and $[\mu]$ denotes the integral
part
of the real number $\mu$. In particular, the fibers are trivial
on $\Omega_1$. Furthermore
the collection of epimorphisms 
$\pi_\omega:{}_A\L(X_B)\to \L_\omega$, $\omega\in
 \Omega$,
defines a continuous bundle of $C^*$--algebras in the
sense of \cite{KW}.}\medskip

\noindent{\it Proof} 
Let us consider the positive map 
$\hat{F}:\L(X_B)\to M(A)$, which satisfies
 $\hat{\phi}(\hat{F}(T))\geq \lambda'T$ for $T\in \L(X_B)^+$
by Cor. 2.11.
We restrict $\hat{\phi}\hat{F}$  
 to a map ${}_A\L(X_B)\to\hat{\phi}(Z(M(A))$
satisfying a corresponding inequality.
Evaluating $\pi_\omega$ on this  inequality yields
$\pi_\omega(\hat{\phi}\hat{F}(T))\geq\lambda'\pi_\omega(T)$,
$T\in{}_A\L(X_B)^+$.
On the other hand for each $\omega$ in the support projection of
 $r-\text{Ind}[X]$,
$\pi_\omega(\hat{\phi}\hat{F}(T))=G_\omega(\pi_\omega(T))$, where
$G_\omega$ is the positive
functional
 of $\L_\omega$ defined as in the proof of the previous lemma: 
$G_\omega(\pi_\omega(T))=\hat{F}(T)(\omega)$. Therefore 
$g_\omega:=((r-\text{Ind}[X])(\omega))^{-1}G_\omega$
is a state of $\L_\omega$ satisfying 
$$(r-\text{Ind}[X])(\omega)g_\omega(\pi_\omega(T))\geq
\lambda'\pi_\omega(T), 
T\in{}_A\L(X_B)^+.$$ It is well known that this condition implies that
$\L_\omega$ is a finite dimensional $C^*$--algebra with at most 
$[{\lambda'}^{-1} (r-\text{Ind}[X])(\omega)]$ minimal orthogonal
projections,
therefore
$\text{dim}(\L_\omega)\leq [{\lambda'}^{-1} (r-\text{Ind}[X])(\omega)]^2$. 

We are left to show that $\omega\in\Omega\to\pi_\omega$ is a continuous 
bundle in the sense of axioms (i)--(iii) of Def. 1.1 in \cite{KW}. 
If $T$ is positive
and satisfies $\pi_\omega(T)=0$ for all
$\omega\in\Omega$ then $T\in J_\omega$ for all $\omega$ in $\Omega$
and therefore $F(T)=0$ which implies $T=0$ by the Pimsner-Popa inequality.
This shows axiom (i). Axiom (ii) is obvious. We are left to show that
for all $T\in{}_A\L(X_B)$, the function
$\omega\in\Omega\to\|\pi_\omega(T)\|$ 
is continuous. We will appeal to the continuity criteria discussed in 
section 2
of \cite{KW}. This function is upper semicontinuous by Lemma 2.3 in \cite{KW} and it
is lower semicontinuous by Lemma 2.2 in the same paper.  Indeed, if 
$\Omega'\subset\Omega$ is a closed subset of $\Omega$ and $D\subset\Omega'$
is dense in $\Omega'$ then the condition $\pi_{\omega}(T)=0$ for all
$\omega\in D$ and some $T\in{}_A\L(X_B)^+$ implies $T\in J_\omega$ for all
$\omega\in D$, thus $F(T)(\omega)=0$ for all $\omega\in D$
and therefore for all $\omega\in\Omega'$ by continuity of the function $F(T)$.
Now evaluating $\pi_\omega$ 
on both sides of the inequality $\hat{\phi}\hat{F}(T)\geq \lambda'T$ shows
that $\pi_\omega(T)=0$ for all $\omega\in \Omega'$.
\medskip

\noindent{\it Remark} Notice that the estimate given in Theorem 3.3
can not be improved in general. In fact, if $H$ is the finite index
${\Bbb C}$--${\Bbb C}$ bimodule defined as in  Example 2.31. Then
$\Omega$
is a one point space, ${}_{\Bbb C}\L(H_{\Bbb C})=M_n({\Bbb C})$,
which is the only fiber. In this case ${\lambda'}^{-1}=\|T^{-1}\|$,
so the corresponding estimate reduces to 
$n\leq \|T^{-1}\|\text{Tr}(T)$ which becomes an equality for $T=I$.
\bigskip

\heading 
4. On the equivalence between finite index and conjugate
equations
\endheading

Our next aim is to show an equivalence between the notion of
$C^*$--bimodule of finite index in the sense of Sect. 2 and
 Longo-Roberts conjugate object in the $C^*$--category of right
Hilbert bimodules
\cite{LR}.\medskip

\noindent{\bf 4.1 The $C^*$--categories $\H_\A$, ${}_\A\H_\A$
and the $W^*$--categories ${\H_\A}^w$, ${{}_\A\H_\A}^w$}\bigskip

\noindent{\bf 4.1 Definition} Let $\A$ be a fixed set of $C^*$--algebras.
We will denote by   ${\Cal H}_\A$ the category 
 with objects and arrows defined as
follows. 
Objects of
 ${\Cal H}_\A$ are    right Hilbert $C^*$--bimodules 
$X$ over elements of $\A$ for which the left action is   
nondegenerate. 
 The
set of arrows $(X,Y)$ in $\H_\A$ between
 two objects
 ${}_AX_B$ and ${}_AY_B$
is the set
$(X,Y):=\L(X_B,
{Y}_B)$ 
of (right) adjointable maps from $X$ to $Y$. 
 Given two objects ${}_AX_B$ and ${}_BY_C$ of $\H_\A$, their
tensor product $X\otimes_B Y$ is still a nondegenerate right Hilbert
$C^*$--bimodule, and therefore it is an object of $\H_\A$. For any
$T\in(X,Y)$,
the map taking a simple tensor $x\otimes y\in X\otimes_B Y$ to
$T(x)\otimes y$, and denoted $T\otimes I_Y$, extends to an adjointable map
on $X\otimes_B Y$. 
For any $C^*$--algebra $A\in\A$,
let $\iota_A$ be  $A$, regarded as a    right Hilbert
bimodule over $A$ itself, in the natural way.  
Since
left action on $\iota_A$ is nondegenerate,  $\iota_A$ is an
object of 
${\Cal H}_\A$.
For any right Hilbert
$A$--$B$ $C^*$--bimodule $X$, the tensor product
 Hilbert bimodule  $X\otimes_B\iota_B$  identifies
naturally with $X$. In general,  
 $\iota_A\otimes_A X$
identifies   with the right Hilbert sub-bimodule of $X$ generated by
$AX$, which, in the case where  the left action is nondegenerate,
still
coincides with $X$ (cf. Def. 2.14).
Therefore    $\{\iota_A, A\in\A\}$
is the set of  left and
right units for the $\otimes$--product between objects.
One can summarize the structure of $\H_\A$, and say that $\H_\A$ is a  
 {\it semitensor}
2-$C^*$--category (in the sense of  \cite{DPZ}).

If we want a   {\it tensor} 2-$C^*$--category  
we need to restrict the arrow spaces, and consider only
 bimodule maps. Namely, let
${}_\A\H_\A$ 
be
 the subcategory
 of $\H_\A$ with the same objects
 and arrows
$(X,Y):={}_A\L(X_B,
{Y}_B)$, the set of right adjointable maps from $X$ to $Y$ commuting with
the left action. This is now a tensor 2-$C^*$--category.

In the sequel we will consider also the  $W^*$--categories 
${\H_\A}^w$
${{}_A\H_\A}^w$ 
with the same objects, and  set of arrows between two
objects $X$ and $Y$ obtained completing the corresponding
arrow spaces of $\H_\A$ and ${}_\A\H_\A$ in a suitable weak topology.
Choose, for each unit object ${}\iota_B\in\H_\A$, a state
$\omega_B$
of $B$, and let us endow  $X$ with the inner product
$$(x,x')_{\omega_B}=\omega_B((x|x')_B), x,x'\in X.$$
Completing $X$, after dividing out by vectors of seminorm zero, with
respect
to this inner product, yields  a Hilbert space
$H_{\omega_B}(X)$.
For each $T\in\L(X_B, Y_B)$, let $\F_{\omega}(T)\in\B(H_{\omega_B}(X),
H_{\omega_B}(Y))$ be the operator which acts by  left multiplication
by $T$. We get in this way a 
  $^*$--functor
$\F_{\omega}:\H_\A\to\H$ to the category $\H$
of Hilbert spaces.
Consider now the {\it universal} $^*$--functor
$\F=\oplus_{\omega}\F_\omega:\H_\A\to\H$,
where the direct sum is taken over all 
choice functions $\omega: B\in \A\to\omega_B$. $\F$  is 
 faithful on arrows and  strictly continuous on the unit
ball of each arrow space.

Define $(X,Y)$ to be the completion  of $\F(\K(X_B,Y_B))$ in the weak
topology of the bounded operators from $\oplus_\omega H_{\omega_B}(X)$
to $\oplus_{\omega} H_{\omega_B}(Y)$,
and let ${\H^w}_\A$ be the $W^*$--category with arrows these 
$W^*$-closed subspaces. 
Since any operator in $\L(X_B, Y_B)$ is the strict limit of a norm bounded
net
from $\K(X_B, Y_B)$, $\F(\L(X_B,
Y_B))\subset (X,
Y)$, therefore
$\H_\A$ becomes a  $C^*$--subcategory of ${\H^w}_\A$ under $\F$.
The universal functor enjoys the following
universality property.
\medskip

\noindent{\bf 4.2 Proposition} {\sl
A $^*$--functor  $\G:\H_\A\to\H$ to the category of Hilbert spaces, 
 strictly continuous on the unit ball of each arrow space of $\H_\A$, 
 extends uniquely to a $^*$--functor
$\G'':{\H^w}_\A\to\H$, normal on the arrow spaces.}\medskip

\noindent{\it Proof}
Let us first assume that
each  Hilbert space $\G(\iota_B)$ is cyclic
for $\G((\iota_B,\iota_B))$.
Let $\xi_B$ be a normalized cyclic
vector. Then, identifying 
$X$ with the subspace of intertwiners $\K(\iota_B,X_B)\subset(\iota_B,X)$,
$\G(X)\xi_B$ is a subspace of the Hilbert space $\G(X)$ associated to the
object $X$. 
We claim that $\G(X)\xi_B$
is the whole $\G(X)$. Let
$\eta\in \G(X)$ be a vector orthogonal to $\G(X)\xi_B$. For all $x\in X$,
$\G(x^*)\eta$ is   orthogonal to $\G((\iota_B,\iota_B))\xi_B$
and hence it is zero. Since $\G$ is strictly continuous on the unit ball
of $(X,X)$, we conclude that $\eta=0$.
We therefore have an identification of the Hilbert space $\G(X)$
with 
$\F_\omega(X)$ where $\omega:B\to\omega_{\xi_B}$, and also an
identification of $\G$ with $\F_\omega$.
Now every $^*$--functor $\G:\H_\A\to\H_{\comp}$
is the direct sum
cyclic $^*$--functors, therefore
$\G$ is a direct sum of some $\F_\omega$, and the rest now follows easily.
\medskip

In particular, if $R_Y:\H_\A\to\H_\A$ is  the  $^*$--functor
which tensors on the right by an object $Y\in\H_A$, the normal extension
of $\F\circ R_Y$, with $\F$ the universal $^*$--functor,  makes
${\H^w}_\A$ into a semitensor 
2-$W^*$-category. 

The subcategory ${}_\A{\H^w}_\A$ of ${\H^w}_\A$ with the same objects
and arrows 
$$(X,Y):=\{T\in\F(\K(X_B,Y_B))'': T\F(\phi(a))=\F(\phi'(a))T, a\in A, x\in
X\},
$$
(where $\phi$ and $\phi'$ denote respectively the left actions of $A$ on 
$X$ and $Y$, and $\F$ is the universal $^*$--functor) is now a tensor
2-$W^*$--category.
 \medskip

\noindent{\it Remark} The functor of $R_Y$ 
 may not be injective on arrows in any of these
categories. In other words,
if ${}_AX_B$,  
${}_A{X'}_B$ and ${}_BY_C$ are  right Hilbert $C^*$--bimodules,  the
natural
 $^*$--homomorphism 
$$T\in\L(X_B,{X'}_B)\to T\otimes 
I_Y\in\L((X\otimes_B Y)_C, (X'\otimes_B Y)_C)$$
may not be injective.
 In fact, if $X=X'=\iota_B$ and $b\in B\subset \L(\iota_B)=M(B)$,
under the identification of $\iota_B\otimes_B Y$ with $Y$,
$b\otimes I_Y$ corresponds to  the left action of $B$ on $Y$ 
evaluated in $b$, which may vanish.\medskip

\noindent{\bf 4.2 Conjugation in ${}_\A\H_\A$ and
${{}_\A\H_\A}^w$}\bigskip

\noindent In the sequel $\T$ will denote either ${}_\A\H_\A$ or
${}_\A{\H^w}_\A$.
Following \cite{LR}, we can introduce the notion of
conjugation in the tensor $C^*$ (or $W^*$) category $\T$.
\medskip

\noindent{\bf 4.3 Definition} Let 
 $X = {}_AX_B$ be  an object of $\T$.
An object $Y = {}_BY_A$ of $\T$
is called a {\it conjugate} of $X$  if there exist intertwiners   
$R \in(\iota_B,Y \otimes_A X)\in \T$ and 
$\overline{R} \in (\iota_A, X \otimes_B Y)\in\T$
such that 

$$
\align
 \overline{R}^*\otimes I_X\circ I_X \otimes R & = I_X \\
 R^* \otimes I_Y\circ I_Y \otimes \overline{R} & = I_Y.
\endalign
$$
We adopt the convention that the $\otimes$--product is evaluated before 
$\circ$--product.
We emphasize that, if $\T={}_\A\H_\A$, $R$ and $\overline{R}$ are
$C^*$--bimodule
maps, i.e.
they commute with left as well as right actions of the appropriate
$C^*$--algebras. Therefore in this case  $R^*R$ and 
$\overline{R}^*\overline{R}$
are 
elements of ${}_{B}\L(\iota_B)=Z(M(B))$ and
${}_A\L(\iota_A)=Z(M(A))$ respectively.
If, instead, $\T={}_\A{\H^w}_\A$, we can only conclude that $R^*R$
and
$\overline{R}^*\overline{R}$
are central elements of $B''$ and $A''$ respectively. 

The above equations will be referred to as {\it the conjugate equations\/}.
 Clearly, if $Y$ is a conjugate of $X$ then $X$ is a conjugate of $Y$.

The {\it dimension of} $X$ {\it relative to the pair} $(R, 
\overline{R})$ is defined by $\text{dim}_{R,\overline{R}} 
X=\|R\|\|\overline{R}\|$.
The {\it minimal dimension} of $X$, denoted   $\text{dim } X$, is 
the infimum of all relative dimensions  $\text{dim}_{R,\overline{R}}$. 
\medskip

\noindent{\it Uniqueness of the conjugate object.} 
Let  $Y$ be a conjugate object of $X$ in   $\T$,
 and let $R$ and $\overline{R}$
solve the corresponding conjugate equations.
  Let $U\in(Y, {Y'})$ be an invertible intertwiner in $\T$.
Set $R':=U \otimes I_X\circ R$ 
and 
$\overline{R'}= I_X \otimes {U^*}^{-1}\circ\overline{R}$.  Then    
 $(Y', R', \overline{R'})$ defines another 
conjugate of
$X$ and    every conjugate of $X$ arises in this way (see \cite{LR}). 
In the case where $U$ is a unitary, 
$\overline{R'}^*\overline{R'}=\overline{R}^*\overline{R}$ 
and ${R'}^*R'=R^*R$, so the  dimension relative to this new 
pair of intertwiners does 
not change.
In this situation we say that the conjugates $(Y,R,\overline{R})$ and 
 $(Y',R',\overline{R'})$ are {\it unitarily equivalent}. 
\medskip

\noindent{\bf 4.3 From finite index to conjugation}\bigskip

\noindent{\bf 4.4 Theorem} {\sl 
Let ${}_AX_B$ be a bi-Hilbertian $C^*$--bimodule. Then left actions on the
underlying
right Hilbert $C^*$--bimodules $X$ and $\overline{X}$ are nondegenerate,
and therefore these are objects of $_{\A}\H_\A$ and 
${}_\A\H^w_{\A}$.
\roster
\item If $X$ is of finite numerical index,  
 $\overline{X}$ is a
conjugate of $X$
in ${}_\A{\H^w}_\A$. More specifically,
if $\{u_\mu\}_\mu$ and $\{v_\nu\}_\nu$ are, respectively, a generalized
right and left basis of $X$, the nets $\overline{R}_\mu:=\sum_{y\in u_\mu}
y\otimes\overline{y}$ and $R_\nu:=\sum_{z\in v_\nu} \overline{z}\otimes z$
converge strongly under the universal $^*$--functor to intertwiners
$\overline{R}\in (\iota_A, X\otimes_A^r\overline{X})$ and 
$R\in(\iota_B, \overline{X}\otimes^r_B X)$ of 
${}_\A\H^w_{\A}$
which do not depend on the choice of the bases,
and solve the conjugate equations. The following relations also hold
for $x, x'\in X$,
$$
\overline{R}^*(\theta^r_{x,x'} \otimes I_{\overline{X}})\overline{R} 
= {}_A(x|x'),$$
$$
{R}^*(\theta^r_{\overline{x},\overline{x'}} \otimes I_{X})R 
= (x|x')_B,$$
$$R^*R =\ell-\text{Ind}[X],$$  
$${\overline{R}}^*\overline{R} =  r-\text{Ind}[X].$$

\item If $X$ is of finite index, 
$R$ and $\overline{R}$ 
 belong to ${}_\A\H_\A$. So the right Hilbert bimodule $\overline{X}$ is a
conjugate of $X$ in 
${}_\A\H_\A$. Their
right adjoint operators are given by:
$${\overline{R}}^*x\otimes\overline{x'}={}_A(x|x'),$$
$$R^*\overline{x}\otimes x'=(x|x')_B.$$
\endroster}\medskip

\noindent{\it Proof} 
By Prop. 2.16
the left (right) action on the right (left) Hilbert
$C^*$--bimodule $X$  is nondegenerate, therefore the right Hilbert
$C^*$--bimodules  $X$ and $\overline{X}$ are objects of ${}_\A\H_\A$ and
${}_\A\H^w_{\A}$.
We claim that, 
under the natural identifications of
$X\otimes^r_B\overline{X}$
with $\K(\iota_A, X\otimes^r_B\overline{X}_A)$ and 
of $\overline{X}\otimes^r_A{X}$ with $\K(\iota_B, \overline{X}
\otimes^r_A{X}_B)$, the nets
$\overline{R}_\mu:=\sum_{y\in u_\mu}y\otimes\overline{y}$ and 
$R_\nu:=\sum_{z\in v_\nu}\overline{z}\otimes z$ converge strongly in 
the univarsal $^*$--functor to operators $\overline{R}$ and $R$ which do not depend on the choice of
the bases.   It suffices  to show  that the first net is strongly Cauchy,
as, replacing $X$ with 
$\overline{X}$,  $\mu\to u_\mu$ 
changes to $\nu\to v_\nu$. Now by
Prop. 2.19,  $\sum_{y\in u_\mu} {}_A(y | y)$ 
is a positive, increasing, norm bounded net,  and it is strongly
convergent  in $A''$ to $r-\text{Ind}[X]$.
Since for $\mu<\mu'$,
$\sum_{y\in u_{\mu'}}\theta^r_{y,y}-\sum_{y\in u_\mu}\theta^r_{y,y}$ 
 is a positive
contraction, we have 
$$
 (\sum_{y\in u_{\mu'}}y \otimes \overline{y}
-\sum_{y\in u_{\mu}}y \otimes \overline{y}
 | 
\sum_{y\in u_{\mu'}} y \otimes
 \overline{y}-\sum_{y\in u_{\mu}}y \otimes \overline{y}
)_A=$$ 
$$F_X((
\sum_{y\in u_{\mu'}}\theta^r_{y,y}-
\sum_{y\in u_{\mu}}\theta^r_{y,y})^2)\leq
F_X(\sum_{y\in u_{\mu'}}\theta^r_{y,y}-
\sum_{y\in u_{\mu}}\theta^r_{y,y})=$$
$$\sum_{y\in u_{\mu'}}{}_A(y|y)-
\sum_{y\in u_{\mu}}{}_A(y|y).$$
Therefore the net 
${\overline{R}}_\mu
\in\K(\iota_A,
X\otimes\overline{X}_A)$ 
is  strongly convergent on a dense subspace of the underlying Hilbert
space.
We show that  this net is norm bounded.
We have, for $a\in A$,
$$
\align
 &  \|(\overline{R}_\mu(a)|\overline{R}_\mu(a))_A\| 
= \|(U(\sum_{y\in u_\mu}(y \otimes \overline{y})a)|
U(\sum_{y\in u_\mu}(y \otimes \overline{y})a))_A\|  \\
 =  &  \|F((\sum_{y\in u_\mu}\theta^r_{y,y}\phi(a))^*
 (\sum_{y\in u_\mu}\theta^r_{y,y}\phi(a)))\|
 \leq \| F(\phi(a)^*\phi(a))\|\\
= & \|(r-\text{Ind} [X])a^*a\| ,
\endalign
$$
where $U$ is the biunitary map defined in Prop. 2.13 (3).
Hence 
$$
\|\overline{R}_\mu\| \leq
(\sup_{a \not= 0} \frac{\|(r-\text{Ind}[X])a^*a\|}{\|a^*a\|})^{1/2} 
 = \| r-\text{Ind}[X] \|^{1/2}.
$$ 
It follows that  $\overline{R}_\mu$ is strongly convergent to an operator 
$\overline{R}\in(\iota_A, X\otimes\overline{X})\subset{\H^w}_\A$
with $\|\overline{R}\|\leq\| r-\text{Ind}[X] \|^{1/2}$.
Similarly we define a  map $R\in(\iota_B, \overline{X}\otimes X)\subset
{\H^w}_\A$ 
as the strong limit of 
 $\sum_{z\in v_\nu}(\overline{z} \otimes z)$ such that
 $\|R\| \leq \| \ell-\text{Ind}[X]\|^{1/2}$.  
In order to show that  $\overline{R}$ is  independent on the basis, we
compute its
Hilbert space adjoint. Let $\omega:B\in\A\to\omega_B$ be a choice of
states of the $C^*$--algebras
of $\A$, and let $\F_\omega:\H_\A\to\H$ be the associated
cyclic $^*$-functor to the category of Hilbert spaces. For $x,x'\in
X$, $a\in A$,
 $$(a, \F_\omega({\overline{R}}_\mu)^*x\otimes \overline{x'})_{\omega_A}=
\sum_{y\in u_\mu} (y\otimes\overline{y}a,x\otimes
\overline{x'})_{\omega_A}=$$
$$\sum_{y\in u_\mu}
\omega_A((\overline{y}a|(y|x)_B\overline{x'})_A)=\sum_{y\in u_\mu}
\omega_A({}_A(a^*y|x'(x|y)_B))=$$
$$\omega_A({}_A(a^*\sum_{y\in u_\mu}y(y|x)_B, x')),$$
Therefore $\F_\omega(({\overline{R}}_\mu)^*x\otimes\overline{x'})$ 
converges weakly to ${}_A(x|x')$, regarded as an element of the Hilbert space
$\F_\omega(\iota_A)$.
It follows that $\overline{R}^*$, 
and hence $\overline{R}$, is independent of the generalized right basis.
On the other hand the net ${\overline{R}}_\mu$, regarded as a net
in the  Hilbert space $\F_\omega(X\otimes\overline{X})$, has norm bounded
above
by $(r-I[X])^{1/2}$, therefore
$$\|\overline{R}\|=\|\overline{R}^*\|\geq
(r-I[X])^{-1/2}\|\overline{R}^*(\sum_{y\in u_\mu} y\otimes\overline{y})\|=$$
$$(r-I[X])^{-1/2}
\|\sum_{y\in u_\mu}{}_A(y|y)\|,$$
which shows that
$\|R\|=(r-I[X])^{1/2}$.

Let now $U\in M(A)$ be a unitary. For any generalized right basis
$u_\mu$, $\mu\to\{Uy, y\in u_\mu\}$ is still a generalized right basis, so 
$U\overline{R}U^*=\overline{R}$
by  independence of the operator $\overline{R}$ on the basis. Hence 
$\overline{R}\in{}_A{\H^w}_\A$.

We show that $R$ and $\overline{R}$ solve the conjugate equations.
For $x \in X, \ b \in B$, we have, in the Hilbert space associated to 
$X$ under the universal $^*$-functor: 
$$
\overline{R}^* \otimes I_X\circ I_X \otimes R(x b)
= \overline{R}^* \otimes I_X
(x\otimes \lim_\nu\sum_{z\in v_\nu}(\overline{z} \otimes z)b) 
= \lim_\nu\sum_{z\in v_\nu} {}_A(x|z)zb
= xb.
$$
Since $X\otimes_B\iota_B$ identifies with $X$ via the map 
$x\otimes b\mapsto xb$,
we obtain the  conjugate equation $ \overline{R}^*\otimes I_X\circ I_X 
\otimes R= I_X$ in ${}_\A{\H^w}_\A$.
Similarly we have 
$R^* \otimes I_Y\circ I_Y \otimes \overline{R}  = I_Y$.

For any $a \in \F(\iota_A)$ we have 
$$
\overline{R}^*\overline{R}(a) = \lim_\mu \overline{R}^*
(\sum_{y\in u_\mu} {y} \otimes \overline{y} a)
        = \lim_\mu\sum_{y\in u_\mu} {}_A(y|a^*y)
        = \lim_\mu\sum_{y\in u_\mu} {}_A(y|y)a ,
$$
so
${\overline{R}}^*\overline R=r-\text{Ind}[X]$ and
$R^*R=\ell-\text{Ind}[X]$ as well.

We  show that 
 $ \overline{R}^*(\theta^r_{x,z} \otimes I_{\overline{X}})\overline{R} =
{}_A(x|z)$
(the similar equation relative to $R$ will follow replacing $X$ with
$\overline{X}$).
For $a \in A$, 
$$
\align
 {} & \  \overline{R}^*(\theta^r_{x,z} \otimes I_Y)\overline{R}(a) 
 = \overline{R}^*(\theta^r_{x,z} \otimes I_Y)
 \lim_\mu(\sum_{y\in u_\mu} (y \otimes \overline{y})a) \\ 
 & = \lim_\mu\overline{R}^*(\sum_{y\in u_\mu}
 \theta^r_{x,z}(y) \otimes \overline{a^*y}) 
   = \lim_\mu\sum_{y\in u_\mu}{}_A(\theta^r_{x,z}(y) | a^*y) \\ 
   &  = \lim_\mu\sum_{y\in u_\mu} {}_A(x(z|y)_B | y)a 
     = {}_A(x | \lim_\mu\sum_{y\in u_\mu} y(y|z)_B )a 
     = {}_A(x|z)a.  
\endalign
$$

(2)
In the case where $X$ is of finite index,
the net $\overline{R}_\mu(a)$ 
 converges in norm for all $a\in A$,  therefore $R$  is actually mapping 
$A$ to 
$X\otimes\overline{X}$.
Furthermore $\overline{R}$ is right adjointable, in fact  its
adjoint 
$\overline{R}^* : X \otimes_B \overline{X} \rightarrow A$ is defined 
 by 
$\overline{R}^*(x \otimes \overline{x'}) = {}_A(x|x')$:
$$
\align
 (\overline{R}(a) | x \otimes \overline{x}')_A 
 & = \lim_\mu\sum_{y\in u_\mu}a^*(y \otimes \overline{y} | x \otimes
 \overline{x}')_A\\ 
 & =\lim_\mu \sum_{y\in u_\mu} a^*(\overline{y} | (y|x)_B
 \overline{x}')_A  
  = \lim_\mu\sum_{y\in u_\mu} a^* {}_A(y | x'(y|x)_B^*)\\ 
  & = a^*{}_A (\lim_\mu\sum_{y\in u_\mu} y(y|x)_B | x')
   =  a^*{}_A(x|x').
\endalign
$$
\medskip   

A first consequence of the previous theorem is the fact that
the left Hilbert bimodule structure on a finite index bimodule
is unique up to equivalence. \medskip 

\noindent{\bf 4.5 Corollary} {\sl Let ${}_AX_B$ be a bi-Hilbertian 
$C^*$--bimodule of finite index.
Any other left inner product on the underlying right Hilbert bimodule
${}_AX_B$
making it into a finite index, bi-Hilbertian
bimodule is of the form
$${}_A(x|y)'={}_A(Qx|y), x,y\in X,$$
where $Q$ is a positive invertible element of $\L_B({}_AX)$.}\medskip

\noindent{\it Proof}
Consider  another left inner product ${}_A(\cdot|\cdot)'$ making $X$
into a bi-Huilbertian, finite index $C^*$--bimodule. 
Let $X'$ denote the left Hilbert bimodule structure over $X$ with 
inner product ${}_A(\cdot|\cdot)'$.
By part (2) of Theorem 4.4, we can find
another solution ($\overline{X'}$, $R'$, $\overline{R}'$) to the conjugate
equations such that 
${}_A(x|x')'=
\overline{R'}^*(\theta^r_{x,x'} \otimes I_{\overline{X'}})\overline{R'}$.
 By uniqueness
of the conjugate object (cf 
 a remark following Definition 4.3) there is an invertible
$U\in{}_B\L_A(\overline{X},\overline{X'})$
such that
$\overline{R}'=I_X\otimes
U\circ\overline{R}$. 
Therefore 
${}_A(x|x')'=
\overline{R}^*(\theta^r_{x,x'} \otimes U^*U)\overline{R}$.
We just need to plug in the fact that
$\overline{R}=\lim_\mu\sum_{y\in
u_\mu}
y\otimes\overline{y}$ in the pointwise norm convergence topology
 and choose $Q:={J_X}^{-1}U^*UJ_X$, with $J_X: X_B\to {}_B\overline{X}$
the
natural conjugation map.
 \bigskip

\noindent{\bf 4.4 On the equality ${}_A\L(X_B)=\L_B({}_AX)$
for finite index bimodules}\bigskip

\noindent Let ${}_AX_B$ be a bi-Hilbertian $C^*$--bimodule. We can
consider the 
$C^*$--algebra ${}_A\L(X_B)$ of right adjointable maps commuting with the
left
action, but also  the $C^*$--algebra  ${}\L_B({}_AX)$ of left
adjointable maps
commuting
with the right action.
If $A$ and $B$ are unital, and $X$ is finitely generated, as a right and
left module, any bimodule map on $X$ is right adjointable and left
adjointable, therefore ${}_A\L(X_B)=\L_B({}_AX)={}_A\text{End}_B(X)$.
More generally, 
under which conditions  ${}_A\L(X_B)=\L_B({}_AX)$ as
algebras? 
The following result provides an answer.\medskip

\noindent{\bf 4.6 Corollary} {\sl If ${}_AX_B$ is a bi-Hilbertian bimodule
of finite index, any element of ${}_A\L(X_B)$ is adjointable
with respect to the left inner product, and therefore it belongs
to $\L_B({}_AX)$. Similarly,  any
element of $\L_B({}_AX)$
is adjointable with respect to the right inner product. Therefore
${}_A\L(X_B)=\L_B({}_AX)$ as algebras.}\medskip

\noindent{\it Proof} Let $R$ and $\overline{R}$ be the solution  
to the conjugate equations arising from the left and right inner products
as in the proof of the previous theorem. By Frobenius reciprocity 
there is a linear isomorphism from ${}_A\L(X_B)$ to 
${}_B\L(\iota_B,
\overline{X}\otimes X_B)$ 
given by $T\to I_{\overline{X}}\otimes T\circ R$
and an antilinear isomorphism from 
${}_B\L(\iota_B,
\overline{X}\otimes X_B)$ 
to ${}_B\L(\overline{X}_A)$ given by $S\to
S^*\otimes I_{\overline{X}}\circ I_{\overline{X}}\otimes\overline{R}$
(see \cite{LR}). A straightforward computation shows that the composition
of these maps is the map $T\in{}_A\L(X_B)\to
JTJ^{-1}\in{}_B\L(\overline{X}_A)$ where $J:{}_AX\to\overline{X}_A$
is the   conjugation map. Therefore  the map 
$T\in{}_A\L(X_B)\to T\in\L_B({}_AX)$ is a linear multiplicative isomorphism.
\medskip

\noindent{\it Remark}
In general 
the $^*$--involution of  ${}_A\L(X_B)$ differs from that of 
$\L_B({}_AX)$.
 We illustrate this phenomenon
in the particular case where $X$ comes from
a conditional expectation. 
Let
$E:B\to A$ be a finite index conditional expectation in the sense of  \cite{W},
between unital $C^*$--algebras. 
Set ${}_BX_A=B$ as a $B$-$A$ bimodule in the natural way, and
with inner products
 $(x|y)_A=E(x^*y)$, $_B(x|y)=xqy^*$, where $q\in A'\cap B$ 
is a positive invertible
element.
Since $E$ satisfies a Pimsner-Popa inequality \cite{W}, there exists a 
positive scalar
 $\lambda$ 
such that $\lambda E-\text{id}$ is completely positive, by \cite{FK}.
Therefore 
$$\lambda\|\sum_1^n\theta^r_{x_i,x_i}\|=
\lambda\|(E({x_i}^*x_j))_{i,j}\|\geq$$
$$\|\sum_1^nx_i{x_i}^*\|\geq
\|q\|^{-1}\|\sum_1^n{}_B(x_i|x_i)\|,$$ therefore $X$ is of finite 
right numerical index and also of finite index since $X_A=B$ is finitely generated
over $A$. On the other hand, since $\theta^l_{x,x}(z)=zqx^*x$,
$$\|\sum_1^n\theta^l_{x_i,x_i}\|=\|q^{1/2}\sum_1^n{x_i}^*x_iq^{1/2}\|
\geq$$
$$\|q^{-1/2}\|^{-2}\|\sum_1^n{x_i}^*x_i\|
\geq\|q^{-1}\|^{-1}\|\sum_1^n(x_i|x_i)_A\|,$$
which shows that $X$ is of finite left numerical index, and therefore of finite index
since ${}_BX$ is singly generated over $B$.
By the previous corollary ${}_B\L(X_A)=\L_A({}_BX)$ and the latter
coincides with $A'\cap B$ acting on $X$ by right multiplication. 
Therefore for any $T\in A'\cap B$, the map $x\in B\to xT$ is adjointable
with respect to the right inner product. 
The $^*$--involution of $\L_A({}_BX)$ (denoted by $T\to{}^*T$) is defined by
the equation 
$${}_B(x|{}^*T(y))={}_B(T(x)|y)=xTqy^*=xq(yqT^*q^{-1})^*,$$
for $T\in A'\cap B$. Therefore
${}^*T=qT^*q^{-1}$ where $T\to T^*$ is the $^*$--involution of 
 $B$. On the other hand the $\tilde{}$--involution 
of ${}_B\L(X_A)$ 
is defined by
$$(x|\tilde{T}y)_A=(T(x)|y)_A=(xT|y)_A=E(T^*x^*y).$$
Since $\tilde{T}$ acts as right multiplication by an element 
of $A'\cap B$, it is
determined by the equation
$$E(b\tilde{T})=E(T^*b), b\in B,$$
which shows that $\tilde{T}=\sum y_iE(T^*{y_i}^*)$, where $\{y_i\}$ is a 
finite quasi--basis of $E$.
Now if $E$ was chosen to satisfy the equation 
$$E(bT)=E(Tb), b\in B, T\in A'\cap B,$$
the $\tilde{}$--involution
(coming from the right inner product)
 and the original involution on $A'\cap B$ coincide. This is possible
if, e.g., $Z(A)$ is finite dimensional. In fact, in this case
$A'\cap B$ is finite dimensional as well, and therefore
$E'(b)=E(\int_{G}ubu^*du)$, with $G$ the unitary group of $A'\cap B$,
is still a conditional expectation from $B$ onto $A$ satisfying the
required equation.
However, the involution on $A'\cap B$ coming from the left inner
product differs from  the original involution
  if $q$ is not central in $A'\cap B$.
\medskip

\noindent{\bf 4.5 Computing the left index element of $X$}\bigskip

\noindent{\bf 4.7 Lemma} {\sl Let $X = {}_AX_B$ and $Y={}_BY_A$ be 
nondegenerate
right  
Hilbert 
$C^*$--bimodules, conjugate of each other as objects of
${}_\A{\H^w}_{\A}$, and let
  ($R$, $\overline{R}$) be a solution of the corresponding conjugate equations.
Let us regard $\K(X_B)$ as a $C^*$--subalgebra of the intertwiner space
$(X,X)\simeq\K(X_B)''$ of ${\H^w}_\A$.
Then the map $$T\in\K(X_B)\to   (\overline{R}^*\circ T \otimes I_Y
\circ\overline{R})\otimes I_X\otimes R^*R-T\in\K(X_B)''$$
is completely positive.
}\medskip

\noindent{{\it Proof} Let us take
 the adjoint of the first 
conjugate equation: $$
 I_X \otimes R^*\circ \overline{R} \otimes I_X = I_X,
$$
thus for all $n\in{\bN}$ and any positive  $T=(T_{ij}) \in M_n(\K(X_B))$,
$$
\align
 T & = (\overline{R}^* \otimes I_X\circ I_X \otimes R T_{ij}  
       I_X \otimes R^ *\circ\overline{R}\otimes I_X)_{i,j} \\
   & = (\overline{R}^* \otimes I_X (T_{ij} \otimes RR^*) \overline{R}
 \otimes I_X)_{i,j} 
   \le  ((\overline{R}^* \circ T_{ij}  \otimes I_Y\circ
\overline{R})\otimes I_X\otimes(R^*R))_{i,j}
\endalign
$$ 
since $RR^*\leq  I_{Y\otimes X}\otimes(R^*R)$ by Lemma 2.7 in \cite{LR}. 
\medskip

\noindent{\it Remark} 
Choosing  $T=I_X$, we obtain, in particular,  $\text{dim}_{R,
\overline{R} } X\geq 1$.
\medskip

Combining the previous lemma with the main theorem of \cite{FK}, 
 yields
 the following result.\medskip

\noindent{\bf 4.8 Theorem} {\sl Let ${}_AX_B$ be a 
bi--Hilbertian bimodule  of finite right numerical
index, and let $F:\K(X_B)\to A$ be the positive $A$--$A$ bimodule map
constructed in Cor. 2.11. Then $X$ is also of finite left numerical index.
  Denoting by $\phi$ and $\psi$
the left and right
actions of $A$ and $B$ on $X$ respectively, and by $q$ the support projection
of the left index element in $B''$,
$\ell-\text{Ind}[X]$ 
 is the smallest central element $c$  of $qB''$
for which the map
$\psi_0(c)
\phi\circ F-\text{id}:\K(X_B)\to\K(X_B)''$
is 
completely positive. Here $\psi_0$ denotes the extension to $Z(B'')$ of 
the right action of $Z(B)$ on $X$ defined in Lemma 2.18.
}\medskip

\noindent{\it Proof}  
We claim that $X$ is of finite left numerical index if and only if
there exists a positive real $c$ for which $c\phi\circ
F-\text{id}:\K(X_B)\to \L(X_B)$ is completely positive.
We show the claim. If $X_B$ has finite left numerical index,
we can construct a solution $R$, $\overline{R}$, $\overline{X}$ 
 to the
conjugate equations 
as in the proof of Theorem 4.4.
We have proved there that 
$R^*R=\ell-\text{Ind}[X]$ and that for
$T\in\K(X_B)$, $\overline{R}^*(T\otimes I_{\overline{X}})
\overline{R}=F(T)$.
So, recalling the definition of tensor products between 
operators in ${}_\A{\H^w}_\A$, with  $\A=\{A,B\}$, we see that
$$I_X\otimes R^*R=\psi_0(\ell-\text{Ind}[X])$$
and
$$
 (\overline{R}^*(T\otimes I_{\overline{X}}) \overline{R})\otimes I_X= 
 \phi\circ F(T),\quad T\in\K(X_B).$$
Inserting these data in the conclusion of Lemma 4.7, we deduce that
$\psi_0(\ell-\text{Ind}[X])\phi\circ F-\text{id}$ is completely positive,
as a map from $\K(X_B)$ to $\K(X_B)''$.
Therefore, with $c=\|\ell-\text{Ind}[X]\|$, $c\phi\circ
F-\text{id}:\K(X_B)\to \L(X_B)$ is completely positive. 
Conversely, if for some positive real $c$, $c\phi\circ F-\text{id}$
is completely positive on $\K(X_B)$, for $n\in {\Bbb N}$ and
for $x_1,\dots, x_n\in X$,
$$\|\sum_1^n(x_i|x_i)_B\|=\|(\theta^r_{x_i,x_j})_{ij}\|\leq
c\|(\phi({}_A(x_i|x_j)))_{ij}\|=c\|\sum_1^n\theta^\ell_{x_i,x_i}\|,$$
so $X$ is of finite left numerical index. 
On the other hand in Prop. 2.19 we have constructed a surjective
conditional expectation
$\phi''\circ E'':\K(X_B)''\to\phi''(A'')$ 
normalizing $\phi''F''$, which does satisfy 
$\mu\|\phi''E''(T)\|\geq\|T\|$ for some positive real $\mu$ and all
$T\in{\K(X_B)''}^{+}$.
By the main result of \cite{FK}, $c\phi'' E''-\text{id}$ is 
completely positive for some positive real $c$, and therefore
$c\|\phi''(r-\text{Ind}[X])^{-1}\|\phi''F''-\text{id}:\K(X_B)''\to\K(X_B)''$
is completely
positive.
Restricting this map to $\K(X_B)$ and combining  with the
claim, shows that
$X$ is of finite left numerical index.

Let now $\nu\to v_\nu$ be a generalized left basis of $X$. Choosing
$n=|\nu|$, $T=(\theta^r_{z_i,z_j})\in M_n(\K(X))^+$, we see that, if
$c$ is any central element of $qB''$ for which $T\in\K(X)\to\psi_0(c)\phi
F(T)-T\in\K(X)''$
is completely positive then 
$$(\psi_0(c)\phi({}_A(z_i|z_j))_{i,j}=
(\psi_0(c)\phi
F(\theta^r_{z_i,z_j}))_{i,j}\geq(\theta^r_{z_i,z_j})_{i,j},$$
which implies
$$\sum_{i,j}(z_i| {}_A(z_i|z_j) z_j)_Bc\geq\sum_{i,j}(z_i|z_i(z_j|z_j)_B)_B,$$
or, in other words,
${R_\nu}^*{R_\nu}c\geq (\sum_{z\in
v_\nu}(z|z)_B)^2$. Thus
 $(\ell-\text{Ind}[X])c
\geq(\ell-\text{Ind}[X])^2$, so 
$c\geq \ell-\text{Ind}[X]$.
\medskip

\noindent{\bf 4.9 Corollary} {\sl  Let $A\subset B$ be an inclusion of
$C^*$--algebras,
and
$E: B\to A$ be a conditional expectation with range $A$, for which there
is $\lambda>0$
such that
$\|E(bb^*)\|\geq\lambda\|bb^*\|$ for all $b\in B$. Let
$\text{Ind}[E]$ 
be the   index of $E$ defined as in Def. 2.17.
Then 
$\text{Ind}[E]$ 
is the smallest  central element $c$ of $B''$ for which
$cE-\text{id}$ is completely positive.}
\medskip

\noindent{\it Proof}
Let ${}_A\overline{X}_B$ be the contragradient of the $B$-$A$ bimodule $X$
associated to $E$ as part (2) of  Prop. 2.12. Clearly $\overline{X}$ is
of finite
numerical index. Since
$\ell-\text{Ind}[\overline{X}]=r-\text{Ind}[X]=\text{Ind}[E]$, and since
$\K(\overline{X}_B)=B$ and $F_{\overline{X}}=E$, $\text{Ind}[E]$ is, by
Cor. 4.9, the
smallest central element $c$ of $B''$ for which $cE-\text{id}$ is
completely positive (recall that $r-\text{Ind}[X]$ is invertible by Cor. 2.20).
\medskip

\noindent{\it Remark}
If 
$\phi E:\K(X_B)\to \phi(A)$ is the faithful conditional expectation 
defined in Cor. 2.30
then 
$$
\text{Ind}[X]\phi\circ E-\text{id}$$
is completely positive on $\K(X_B)$. In fortunate cases where
$\phi''(r-\text{Ind}[X])$ is central in $\K(X_B)''$ (e.g. either 
$A$ is simple, cf. Cor. 2.26, or $X$ arises from a conditional
expectation,
Prop. 2.12, or $A=B$ is commutative and 
right action coincides with left action)
then
$\text{Ind}[X]=\text{Ind}[\phi\circ E]$. This observation thus shows that
the index element of a conditional expectation coincides with the index element
of the dual conditional expectation.
\medskip

\noindent{\bf 4.6 From conjugation to finite index}
\bigskip

\noindent Let $X$ be a bi-Hilbertian bimodule with a conjugate in
${}_\A\H_\A$.
Our next aim is 
 to construct a left 
inner product on $X$ making it into a bi-Hilbertian 
bimodule of finite index.
\medskip

\noindent{\bf 4.10 Lemma} {\sl Let  ${}_BY_A$ be a conjugate
object of ${}_AX_B$ in the tensor 2--$C^*$--category
 ${}_\A\H_\A$, and let   $R$ and
$\overline{R}$ be a pair of intertwiners solving  the conjugate equations,
in the sense of Def. 4.3.
 There exist unique positive semidefinite left inner products 
 on $X$ and $Y$  such that 
$$
   {}_A(x|a^*x') = \overline{R}^*(\theta^r_{x,x'}\otimes I_Y) 
\overline{R}(a) \in A     \qquad \qquad   \text{for} \qquad   
a \in A, x, x'\in X,
$$
$$
   {}_B(y|b^*y') = {R}^*(\theta^r_{y,y'}\otimes I_X) 
{R}(b) \in B     \qquad \qquad   \text{for} \qquad   
b \in B, y, y'\in Y.
$$
}\medskip

\noindent{\it Proof} For $x, x' \in X$ we define an 
element ${}_A(x|x') \in M(A)
 = \L(\iota_A)$ by
$$
  {}_A(x|x')(a) = \overline{R}^*(\theta^r_{x,x'}\otimes I_Y) 
\overline{R}(a) \in A    \qquad   \text{for} \quad   a \in A.
$$
 Then $x,x'\mapsto {}_A(x|x')$ 
defines a continuous sesquilinear form on $X$ with 
values in 
$M(A)$ .  We claim that ${}_A(x|x') \in A$.  
 Let $\{u_i\}_i$ be a selfadjoint approximate unit of  $A$ with $\|u_i\|\le 1$.
 We show that  $\{\overline{R} (\theta^r_{x,x'} \otimes I_Y)
R(u_i)\}_i$ 
is a
 norm Cauchy net.  First we assume that $x'$ is of the form $y=a^*x''$ for
some 
 $a \in A$ and $x'' \in X$.  Since
$$
 \overline{R}^* (\theta^r_{x,a^*x''} \otimes I_Y) \overline{R}(u_i)
 = \overline{R}^* (\theta^r_{x,x''} \otimes I_Y) \overline{R}(au_i),
$$
and $au_i \to a$ in norm, $\{\overline{R}^* (\theta^r_{x,a^*x''} 
\otimes I_Y) \overline{R}(u_i) \}_i$ is a Cauchy net in norm.  
For a general element $x' \in X$, we choose $\tilde{x}\in AX$ sufficiently 
close to $x'$ (this being possible as left $A$-action is nondegenerate),
so
$$ 
\align
 &  \|\overline{R}^* (\theta^r_{x,x'} \otimes I_Y) \overline{R}(u_i) 
 - \overline{R}^* (\theta^r_{x,x'} \otimes I_Y) \overline{R}(u_j)\|
 \\
  \le \  & \|\overline{R}^* (\theta^r_{x,x'} \otimes I_Y) 
\overline{R}(u_i) 
 - \overline{R}^* (\theta^r_{x,\tilde{x}} \otimes I)
 \overline{R}(u_i)\| \\  
  \qquad &  + \   \|\overline{R}^* (\theta^r_{x,\tilde{x}} \otimes I_Y)
 \overline{R}(u_i) 
 - \overline{R}^* (\theta^r_{x,\tilde{x}} \otimes I_Y)
 \overline{R}(u_j)\|  \\
   \qquad & + \|\overline{R}^* (\theta^r_{x,\tilde{x}} \otimes I_Y)
  \overline{R}(u_j)  - \overline{R}^* (\theta^r_{x,x'} \otimes I_Y)
 \overline{R}(u_j)\|  \\
  \le & \|\overline{R}\|^2 \|x\|\|x'-\tilde{x}\|
 + \|\overline{R}\|^2 \|x\| \|x'-\tilde{x}\| \\
 + & \|\overline{R}^*(\theta^r_{x,\tilde{x}}\otimes I_Y)
\overline{R}(u_i)-\overline{R}^*(\theta^r_{x,\tilde{x}}\otimes I_Y)
\overline{R}(u_j)\|.
\endalign
$$
Thus $\{\overline{R}^* (\theta_{x,x'} \otimes I)  
\overline{R}(u_i) \}_i$  
is still a Cauchy net in $A$.   
\par
For  $a \in A$, we have  
$$
 \overline{R}^*(\theta^r_{x,x'} \otimes I_Y)\overline{R}(u_i)a  
 = \overline{R}^*(\theta^r_{x,x'}\otimes I_Y)\overline{R}(u_ia) 
  \to \overline{R}^*(\theta^r_{x,x'}\otimes I_Y)\overline{R}(a)
$$
in norm.  This shows that the limit of the Cauchy net  
$$\{\overline{R}^* (\theta^r_{x,x'} \otimes I_Y)  
\overline{R}(u_i) \}_i $$ in $A$  coincides with  
${}_A(x|x') \in A$ and 
does not depend on 
the choice of approximate unit $\{u_i\}_i$.
\par
It is easy to see that $(x,x') \mapsto {}_A(x|x')$ is left $A$-linear and 
right 
 conjugate $A$-linear.  
Since for $b,c \in A$ and $x,x' \in X$,   
$$
(({}_A(x|x'))^* b)^*c 
= b^*\overline{R}^*(\theta^r_{x,x'}\otimes I_Y)\overline{R}(c) $$
$$
= b^*\lim_{i\to\infty} \overline{R}^*(\theta^r_{x,x'}\otimes I_Y)
\overline{R}(u_i)c
= ({}_A(x'|x)b)^*c, 
$$
we have  ${}_A(x|x')^* = {}_A(x'|x)$.  
Since 
$$
({}_A(x|x)(a)|a)_A = 
(\overline{R}^*(\theta_{x,x}\otimes I)\overline{R}(a))|a)_A 
= ((\theta_{x,x}\otimes I)\overline{R}(a))|\overline{R}(a))_A \ge 0
$$
in the canonical  Hilbert $A$--module $\iota_A = A_A$ with $(a|b)_A =
a^*b$, we 
have 
 ${}_A(x|x) \ge 0$. 
Now exchanging the roles of $X$ and $Y$, and of $R$ and $\overline{R}$, and applying
this argument to $Y$, we deduce the existence of a left inner product on $Y$ as well.
\medskip
  
  We  next show that $X$ and $Y$ aquire a structure of left Hilbert modules.
To do so, we construct isomorphisms with the contragradient left Hilbert bimodules
$\overline{Y}$ and $\overline{X}$ respectively.
For a $A$--$B$ bimodule $X$, we shall denote by $J_X: X\to\overline{X}$
the map associating $\overline{x}$ to $x$, for any $x\in X$. 
Clearly, $J_X(ax)=J_X(x)a^*$ and $J_X(xb)=b^*J_X(x)$ 
for $a\in A$, $b\in B$, $x\in X$.
\medskip

\noindent{\bf 4.11 Lemma}
{\sl Let $X$ and $Y$ be conjugate objects of ${}_\A\H_\A$, and let us
endow them 
with left inner products defined, as in the    previous lemma, by a pair of 
intertwiners $R$ and $\overline{R}$
solving the conjugate equations. Then there exist natural
 bimodule isomorphisms $U : \overline{Y} \rightarrow X$ 
and $V:\overline{X}\to Y$ from the contragradient bimodules,
which preserve the corresponding left and right  inner products and satisfy
$$VJ_X=(UJ_Y)^{-1}.$$
In particular,
 $X$ and $Y$ become bi-Hilbertian $C^*$--bimodules.}\medskip    

\noindent{\it Proof}
For $y \in Y$ we define $\l_y : X \to Y \otimes_A X$ by 
$\l_y(x)= y \otimes x$.   
Then we have $\l_y^*(y' \otimes x) = (y|y')_Ax$.  
We notice that the set $\{a^*xb | a \in A,\,x \in X,\,b \in B\}$ is
 total in $X$ since, by assumption, left action is nondegenerate.  
\par
By the first conjugate equation
$$
a^*xb =  (\overline{R}^*\otimes I_X)(I_X \otimes R)(a^*xb) 
       =  (\overline{R}^* \otimes I_X)(a^*x \otimes R(b)).
$$
For $\tilde{a}\in A$, $x' \in X$
$$
\align
 & ((\overline{R}^* \otimes I_X)(a^*x \otimes R(b))|\tilde{a}x')_B 
  = (a^*x \otimes R(b)|\overline{R}(\tilde{a})\otimes x')_B \\
  = & (x \otimes R(b) | a\overline{R}(\tilde{a}) \otimes x')_B 
  = (x \otimes R(b) | \overline{R}(a) \otimes \tilde{a}x')_B \\
  =  & (R(b)|{\l_x}^*(\overline{R}(a)) \otimes \tilde{a} x')_B 
  = (R(b)|\l_{{\l_x}^*(\overline{R}(a))} (\tilde{a}x'))_B \\
   =  & ({\l_{{\l_x}^*(\overline{R}(a))}}^*(R(b))| \tilde{a} x')_B.
\endalign
$$
Thus for $a\in A$, $b \in B$ and $x \in X$,  
$$
 a^*xb = {\l_{{\l_x}^*(\overline{R}(a))}}^*(R(b)) \in X. 
$$ 
This shows that $\{{l_y}^*(R(b))| y \in Y, b \in B \}$ 
is total in $X$.  
Similarly, for $a\in A$, $b\in B$ and $y \in Y$,
$$
 b^*ya = {\l_{{\l_y}^*(R(b))}}^*(\overline{R}(a)) \in Y.
$$
We next show that for $y$, $y' \in Y$, $b$, $b' \in B$
$$
 {}_A({\l_{y'}}^*R(b')| {\l_y}^* R(b)) = 
{}_A(\overline{y'}b'|\overline{y}b),
$$
where the left hand side is computed with respect to
 the new inner product on $X$ 
introduced in Lemma 4.10 and the right hand side with respect to
the left inner product on $\overline{Y}$ defined in the paragraph following
Def. 2.8.
We start from the right hand side.  For $a$, $a'\in A$, 
$$  
\align 
 & {}_A(a'\overline{y'}b' | a\overline{y}b)
  = ({b'}^* y'{a'}^* | b^*ya^*)_A   
  =  ({\l_{{\l_{y'}}^* R(b')}}^* \overline{R}({a'}^*) |
 {\l_{{\l_y}^* R(b)}}^*\overline{R}(a^*))_A \\
   = & (\l_{{\l_y}^*R(b)}\l_{{\l_{y'}}^*R(b')} ^* \overline{R}({a'}^*)|
 \overline{R}(a^*) )_A  
   =  (( \theta^r_{{\l_y}^*R(b),{\l_{y'}}^*R(b')} \otimes I_Y)
 \overline{R}({a'}^*)|\overline{R}(a^*))_A \\ 
  = & ({\overline{R}}^*( \theta^r_{{\l_y}^*R(b),
{\l_{y'}}^*R(b')} \otimes I_Y) \overline{R}({a'}^*)|a^*)_A 
  =   ({}_A( \l_{y}^* R(b)| \l_{y'}^* R(b')){a'}^* | a^*)_A \\ 
 =  & ({}_A( \l_{y}^* R(b)| \l_{y'}^* R(b')){a'}^*)^* a^* 
 =    a' {}_A( \l_{y'}^* R(b') | \l_y^* R(b)) a^*.
\endalign
$$
Therefore  $U : \overline{y}b \in \overline{Y} \mapsto {\l_y}^* R(b) \in 
X$ is well
 defined and extends to a left $A$-linear  map from
 $\overline{Y}$ to $X$ preserving  left  $A$-valued inner product.  
Since  the right $A$-valued 
 inner product of $Y$ is definite, the left $A$-valued  inner product 
on $X$ constructed in Lemma 4.10 is also definite.  
Clearly
 this map is also right $B$-linear. Since $\overline{Y}$ is a left Hilbert bimodule,
so is $X$ with respect to its left inner product. 
 Since this left inner product is continuous with respect to the right
one, $X$ is bi-Hilbertian by general Banach space theory.

Similarly, $V : \overline{x}a \in \overline{X} \mapsto 
{\l_x}^*\overline{R}(a) \in Y$
  extends to a $B$--$A$ linear  map preserving the left inner 
product
 from $\overline{X}$ to $Y$ and making $Y$ into a left Hilbert bimodule.
Now   $UJ_Y$ takes $b^*y$ to ${\l_y}^*R(b)$ and $VJ_X$ takes $a^*x$ to 
${\l_x}^*\overline{R}(a)$. Therefore $UJ_YVJ_X$ takes $a^*xb$ to
$$UJ_Y({\l_{xb}}^*\overline{R}(a))=UJ_Y(b^*{\l_{x}}^*\overline{R}(a))=
{\l_{{\l_{x}}^*\overline{R}(a)}}^*R(b)$$
which we have already shown to coincide with $a^*xb$. One similarly shows that
$VJ_XUJ_Y=I_Y$. Since $U$ preserves the left 
inner products, $J_XUJ_Y=V^{-1}$, and therefore
$V$,
preserves the right inner products. For the same reason, $U$ preserves the right inner
products as well.
\medskip

\noindent{\bf 4.12 Lemma}
{\sl Let ${}_AX_B$ and ${}_BY_A$ be right Hilbert $C^*$--bimodules with
nondegenerate left actions, conjugate of each other in ${}_\A\H_\A$.
Consider $X={}_AX$ as a  left Hilbert $C^*$--bimodule
 with left inner product  defined
  by a solution ($R$, $\overline{R}$) of the conjugate 
equations as in Lemma 4.10 (cf. Lemma 4.11).
Then the range of 
the right  $B$--action on
 $X$ is contained   in $\K({}_AX)$.
Similarly, regarding $X=X_B$ as a right Hilbert $C^*$--bimodule with its 
original right
 inner product, the range of the left $A$--action is included in 
    $\K(X_B)$.   
}\medskip

\noindent{\it Proof}   
  Each element of the form $R(b)$ can be approximated, in the norm induced
by the right inner product of $Y\otimes_A X$,
 by finite sums $\{\sum_i y_i \otimes x_i\}$. In turn, as seen in the course 
of the proof of Lemma 4.11,
 each element of $Y$ can 
be
 approximated by finite sums $\{\sum_j{\l_{x_j}}^*\overline{R}(a_j)\}$.
Thus for any $b \in B$ and any positive integer $n$, there exist 
$x^{(n)}_k$, $w^{(n)}_k  \in X$ and $a ^{(n)}_k \in A$ for 
$k = 1,\ldots, N_n$ such that 
$$
r_n := \sum_{k=1}^{N_n} 
{\l_{{x^{(n)}_k}}}^*\overline{R}(a^{(n)}_k) \otimes w^{(n)}_k 
\rightarrow R(b)   \qquad \text{as} \ n \rightarrow \infty.  
$$
By the first conjugate equation, for any $x \in X$ 
$$
xb = (\overline{R}^* \otimes I_X) (I_X \otimes R)(xb)
   =  (\overline{R}^* \otimes I_X) (x \otimes R(b)).
$$
Define $Q(x):= xb \in X$ and 
$Q_n(x) := (\overline{R}^* \otimes I_X) (x \otimes r_n) \in X$.  
We  claim that $Q_n$ is  a finite rank operator on the left Hilbert module 
${}_AX$. We show the claim.
$$
\align
Q_n(x) & = (\overline{R}^* \otimes I_X)(x \otimes r_n)
         = \sum_{k=1}^{N_n} \overline{R}^*
           (x \otimes {\l_{{x^{(n)}_k}}}^*\overline{R}(a^{(n)}_k)) 
            \otimes w^{(n)}_k \\
       &  = \sum_{k=1}^{N_n} \overline{R}^*
            ((\theta^r_{x,x^{(n)}_k} \otimes 
I_Y)\overline{R}(a^{(n)}_k))
            \otimes w^{(n)}_k
          = \sum_{k=1}^{N_n} {}_A(x|x^{(n)}_k)a^{(n)}_kw^{(n)}_k \\
       & = \sum_{k=1}^{N_n} \theta^{\l}_{a^{(n)}_kw^{(n)}_k, 
x^{(n)}_k}(x).
\endalign 
$$
Using the norm on $X$ induced by the right inner product: 
$$
\|Q(x) - Q_n(x)\| \le 
\|\overline{R}^* \otimes I_X \| \|x\| \|R(b) - r_n \|.
$$
On the other hand, 
 the  norm on $X$ coming from the right inner product 
is equivalent to the norm coming from the left inner product,
 therefore
$Q$ is the norm limit of  $\{Q_n\}$ in $\L({}_AX)$, and this shows
that  right action of $ B$ on $X$ lies in $\K({}_AX)$.

Replacing now $Y$ with $X$, we deduce that the right action of $A$ on $Y$
is compact with respect to the left inner product of $Y$. But by Lemma
4.11,
the left Hilbert $C^*$--bimodule $Y$ identifies, through the map $V$, with the 
contragradient $\overline{X}$
of the original right Hilbert $C^*$--bimodule $X$, therefore the left action of 
$A$ on $X$ is compact with respect to the original right inner product of
$X$ itself.\medskip

\medskip

\noindent{\bf 4.13 Lemma}
{\sl Let $X$ be an object of ${}_\A\H_\A$ with a conjugate object $Y$ in
${}_\A\H_\A$, 
and let us make
$X$ and $Y$ into  bi-Hilbertian $C^*$--bimodules with left inner
products defined, as in Lemma 4.10, by a solution ($R$, $\overline{R}$) of
the
conjugate equations. Let  us identify $Y$, as a bi-Hilbertian $C^*$--bimodule,
with $\overline{X}$ via the biunitary map $V:\overline{X}\to Y$ defined in 
Lemma 4.11. 
Then
for any $x,x' \in X$,  
$$ \overline{R}^* (x \otimes \overline{x}') = {}_A(x|x') 
\quad \text{and} \quad 
 R^*(\overline{x}\otimes x')  = (x | x')_B.
$$}\medskip

\noindent{\it Proof} 
We shall show that  
$$
 \overline{R}^*(x \otimes V\overline{x'}) = {}_A(x|x') 
  \qquad \text {and} \qquad 
 R^*(V\overline{x} \otimes x') = (x | x')_B.
$$
The first equation follows from  
$$
  \overline{R}^*(x \otimes V(\overline{x'}a)) 
  = \overline{R}^*(x \otimes {\l_{x'}}^* \overline{R}(a)) 
 = \overline{R}^*(\theta^r_{x,x'} \otimes I_Y) \overline{R}(a) 
 = {}_A(x|x')a.
$$
Similarly, we have 
$$
 R^*(y \otimes U\overline{y}')) = {}_B(y|y'),
$$
where the operator $U$ is still defined in Lemma 4.11.
Now writing $y=V\overline{x}$ and $y'=V\overline{x'}$, and using 
the relation $UJ_YVJ_X=I_X$ obtained in Lemma 4.11, gives
$$R^*V\overline{x}\otimes x'={}_B(V\overline{x}|V\overline{x'})=(x|x')_B.$$
\medskip

We are now ready to prove a   converse of part (2) of  Theorem
4.4.
\medskip

\noindent{\bf 4.14 Theorem} 
 {\sl Let  $X$ be  a right  
Hilbert $A$--$B$ $C^*$--bimodule 
with a nondegenerate left 
action. If
 $X$ has
a conjugate object  in the 2--$C^*$--category ${}_\A{\Cal H}_\A$ of
nondegenerate
right Hilbert bimodules,   it can be given a  left $A$-valued
inner product making it into a 
finite index bi-Hilbertian $C^*$--bimodule. More precisely, any
solution $(Y, R, \overline{R})$
to the conjugate equations in ${}_\A{\Cal H}_\A$
induces a left inner product defining a finite index bi-Hilbertian structure
 on $X$
 by
$${}_A(x|ax')=
\overline{R}^*\theta^r_{x,x'}\otimes 
I_Y\overline{R}(a^*),\quad x,x'\in X, a\in A.\eqno (4.1)$$
It turns out that $Y$ is biunitarily equivalent, as a bi-Hilbertian bimodule,
to $\overline{X}$ and, under this identification, the intertwiners 
$R$ and $\overline{R}$ are defined by
$${\overline{R}}^*x\otimes\overline{x'}={}_A(x|x'),
\quad R^*\overline{x}\otimes x'=(x|x')_B.\eqno(4.2)$$}
\medskip

\noindent{\it Proof} 
 Suppose  that $X$ has a conjugate 
$Y$  defined by intertwiners $R$ and $\overline{R}$. 
So far we have proved that a solution 
$R$, $\overline{R}$ of the conjugate equations
induces a bi-Hilbertian structure on $X$ (Lemma 4.11)
 in such a 
way that the left and right actions have range into the corresponding
compact operators (Lemma 4.12). 
Also, we have been able to identify $Y$ biunitarily with $\overline{X}$
(via
the map $V$ defined in Lemma 4.11)  with $R$ and 
$\overline{R}$ acting as in Lemma 4.13. 
Since $\theta^r_{x,x}$ is positive, we have
$$
 \|R\|^{-2} \theta^r_{x,x} \le \overline{R}^*(\theta^r_{x,x}
 \otimes I)\overline{R} = {}_A(x|x),
$$ 
 by Lemma 4.7. 
 Therefore for any $x_1, \ldots , x_n \in X$,
$$
 \|\sum_{i=1}^n\theta^r_{x_i,x_i} \| \le \|R\|^2 
\|\sum_{i=1}^n{}_A(x_i|x_i)\|.
$$
On the other hand, since $\overline{R}$ is bounded, for any
$x_1, \ldots , x_n, y_1, \ldots , y_n \in X$ we have 
$$
 \|\sum_{i=1}^n{}_A(x_i|y_i) \| = \| \overline{R}^*(\sum_{i=1}^n
 \theta^r_{x_i,y_i} \otimes I) \overline{R}\| \le \|\overline{R} \|^2
 \|\sum_{i=1}^n \theta^r_{x_i,y_i} \|.
 $$ 
Therefore, taking into account Theorem 2.22, all the assumptions of
Definition 2.23 are satisfied, and this
shows that $X$
is of finite right index.
Similarly, $Y=\overline{X}$ is of finite right index, i.e. $X$ 
 is of finite left index as well, 
and therefore of finite index. 
\medskip

The arguments of the proof show that the minimal dimension 
of a 
bimodule is the infimum of the square roots of the numerical indices.
\medskip
 
\noindent{\bf 4.15  Corollary} 
{\sl $$\text{dim} \ X = \inf (r-[X])^{1/2} 
(\ell-I[X])^{1/2} 
= \inf I[X]^{1/2}, $$ 
where the infimum is taken over all possible  left inner products on
the right Hilbert bimodule 
$X$ making it into a finite index bi-Hilbertian $C^*$--bimodule.}\medskip

\noindent{\bf 4.7 A characterization of strong Morita equivalences}
\bigskip

\noindent We next characterize strong Morita equivalences among general
right
Hilbert $C^*$--bimodules  as those objects with  minimal
dimension (or numerical index) equal to $1$.
\medskip

\noindent{\bf 4.16 Corollary} 
{\sl For   a right Hilbert $C^*$--bimodule ${}_AX_B$ 
the following properties are equivalent.
\roster
 \item    $X_B$ is full and it can be given a  full left inner product 
making it into
      finite index bi-Hilbertian bimodule with respect to which 
$\text{I}[X] = 1$,
\item $X$ is an object of the category  of nondegenerate full right 
Hilbert
$C^*$--bimodules with a conjugate  such that $\text{dim} X=1$,
\item $X$  can be given a left inner product making it into
      finite index Hilbert bimodule with $r-\text{Ind}[X]=I_A$ 
and $\ell-\text{Ind}[X]=I_B$,
 \item $X$  can be given a left inner product making it into
a  strong Morita equivalence bimodule from $A$ to $B$.
\endroster
}\medskip

\noindent{\it Proof} (1)  $\Rightarrow$ (2): This implication follows 
from the previous 
theorem.
\par
\noindent
(2) $\Rightarrow$ (4):
Let $R$ and $\overline {R}$ satisfy the conjugate
equations with $\|R\|\|\ov R\|< \sqrt{2}$. Since $X$ and its conjugate are
full, the left and right indices of $X$ must be invertible by Cor.
2.29, and therefore so are $R^*R$ and ${\overline{R}}^*\overline{R}$.
The operators  $S:=R(R^*R)^{-1/2}$
and $\overline{S}:=\overline{R}(\overline{R}^*\overline{R})^{-1/2}$ are 
isometries
whose ranges generate,  as in \cite{LR}, projections
 satisfying the Jones relations with parameter $\beta$, 
where
 ${\beta}^{-1}=(\|R\| \|\ov R\|)^2<2$, thus $\beta=1$ by Jones 
fundamental
 result \cite{J}. This shows that the numerical index $I[X]$ of $X$ with
respect to the original right inner product and the left inner 
product induced by this pair,
is $1$.
Let $\phi E:\K(X_B)\to \phi(A)$ denote the  conditional
expectation defined in Cor. 2.30. Since,
by Cor. 4.9,
$I[X]\phi E(T) \ge T$ for any positive $T$ in $\K(X_B)$, and
  since 
$\phi E(\phi E(T)-T) = 0$,  
   $\phi E$, being faithful, must be  the identity map. 
Defining a new left inner product on $X$
by:
$$
{}_A(x|y)' := ({\overline{R}}^*\overline{R})^{-1}{}_A(x|y) = \theta^r
_{x,y},
$$ 
 makes $X$ into a strong Morita  equivalence bimodule.
\par
\noindent
(4) $\Rightarrow$ (3):
It is easy to show  that a strong Morita equivalence bimodule
is full as a left as well as a right Hilbert module and has
 index 1.   In fact, in this case $A=\K(X_B)$
and one has a bi-Hilbertian structure given by 
${}_A(x|x') = \theta^r _{x,x'}$, for $x,x'\in X$.
 Let 
$\{u_\mu\}_\mu$ be a generalized right basis for $X$.   
Since 
${}_A(x|x') = \theta^r_{x,x'}$, for $x,x'\in X$, we have 
$$
r-\text{Ind}[X] = \lim_\mu\sum_{y\in u_\mu} {}_A(y|y) = \lim_\mu\sum_{y\in
u_\mu} \theta^r_{y,y} = 
I_A.
$$
One similarly shows that    $\ell-\text{Ind}[X] =I_B$. 
\par
\noindent
(3) $\Rightarrow$ (1): This implication is obvious.
\bigskip

\heading 5. Tensoring finite index bimodules
\endheading

Let $A$, $B$ and $C$ be $\sigma$ unital \cst--algebras, $X$  a 
right Hilbert $A$--$B$ bimodule and $Y$  a right Hilbert $B$--$C$ bimodule.
If $X$ and $Y$ are of finite index, is $X\otimes^r_BY$ still of finite
index?
  It does not 
seem to be  easy to prove this property  directly. 
However, one can obtain a proof from our characterization
  Theorems 4.4 and  4.14.
We anticipate  the following well known lemma. \medskip

\noindent{\bf 5.1 Lemma} 
{\sl Let $Y={}_BY_C$, $Y'={}_BY'_{C}$ be right Hilbert $C^*$--bimodules and 
 $F \in {}_B\L(Y_C, Y'_C)$ be a  
$B$--$C$ bimodule homomorphism with adjoint. Then  
\roster
\item for any right Hilbert $A$--$B$ bimodule $X$,
 $\|I_X\otimes F\| \leq\|F\|$, and the 
equality holds if $X_B$ is full and the left action of 
$B$ on $Y$ is nondegenerate,
 \item for any right Hilbert $C$--$A$ bimodule $X$, 
$\|F\otimes I_X\|\leq\|F\|$ and the 
equality holds if the left $A$--action is faithful.
\endroster
}\medskip
  
\noindent{\it Proof} (1) The operator $I_X\otimes F$ is clearly well defined on the linear
span of simple vectors $x\otimes y$, $x\in X$, $y\in Y$, which is  dense in
$X \otimes_B Y_C$. One can easily check that the norm of $I_X\otimes F$ on that subspace 
is bounded above by $\|F\|$, therefore $I_X\otimes F$ extends to a bounded operator with the same norm
on the completion
$X \otimes_B Y_C$ with adjoint $I_X\otimes F^*$. It is also obvious that $I_X\otimes F$ is a bimodule map.
For the rest it suffices to assume $Y=Y'$.
The map
$F \in {}_B\L(Y_C)\to I_X\otimes F\in{}_A\L(X \otimes_B Y_C)$ is a
$^*$--homomorphism. 
 $I_X\otimes F=0$ implies
$$(x\otimes Fy| x\otimes Fy)_C=(Fy, (x|x)_BFy)_C=0,\quad x\in X, y\in Y,$$
therefore $F((x|x)_By)=0$ for all $x\in X$, $y\in Y$, which implies $F=0$ since the right inner product
of $X$ is full and the left $B$--action on $Y$ is nondegenerate. Therefore $\|I_X\otimes F\|=\|F\|$.
(2) The fact that $F\otimes I_X$ is a well defined map on
$Y\otimes X$ with norm bounded above by $\|F\|$ can be proved with arguments similar to those used above.
Now $F\otimes I_X=0$ implies $(x|(Fy|Fy)_Cx)_A=0$, $x\in X$, $y\in Y$, so the left action of $C$ on $X$
evaluated
on $(Fy|Fy)_C$ vanish for all $y\in Y$. If this action is faithful, $Fy=0$, $y\in Y$ and therefore $F=0$.
This implies that $\|F\otimes I_X\|=\|F\|$.
\medskip

We have shown in Prop. 2.13 that if ${}_AX_B$ is of finite right numerical
index and ${}_BY_C$ is of finite left numerical index, the
seminorms of $X\odot_B Y$ arising from the left and right inner products
are equivalent, therefore we can form a unique bi-Hilbertian bimodule,
$X\otimes_B Y$ completing in any of these seminorms. We now show that
this bimodule is of finite index if so are $X$ and $Y$.\medskip

\noindent{\bf 5.2 Theorem} {\sl Let $A$, $B$ and $C$  
be \cst \--algebras, and $X = {}_AX_B$ 
and $Y = {}_BY_C$ be  bi-Hilbertian $C^*$--bimodules.
If ${}_A{X}_B$ and ${}_B{Y}_C$ have finite
 index (respectively,  finite numerical index), then also $X\otimes_B
Y$ has finite index (respectively, finite numerical index)
with respect to the bi-Hilbertian structure defined in Subsect. 2.2.
}\medskip

\noindent{\it Proof}
Since $X$ and $Y$ are bi-Hilbertian and of finite numerical index 
$X\otimes_B Y$
is bi-Hilbertian by Prop. 2.13, and therefore  left and right actions
are
nondegenerate by Prop. 2.16. Since $X$ and $Y$ have finite
numerical index, the contragradient of the corresponding underlying left
Hilbert modules are their respective conjugates, by Theorem 4.4. Namely,
there are intertwiners in ${}_\A{\H^w}_\A$, with $\A=\{A,B\}$, 
$\overline{R}_1\in (\iota_A, X\otimes_B \overline{X})$,
${R}_1\in(\iota_B,\overline{X}\otimes_A\overline{X})$ ,
$\overline{R}_2\in (\iota_B, Y\otimes_C \overline{Y})$,
${R}_2\in(\iota_C,\overline{Y}\otimes_B\overline{Y})$ 
solving the corresponding conjugate equations.
We show that $\overline{Y}\otimes_B \overline{X}$ is a conjugate of
$X\otimes _B Y$ in ${}_\A{\H^w}_\A$.
We define a map  $i(R_1)$ from  
$\overline{Y} \otimes Y \simeq \overline{Y} \otimes \iota_C \otimes Y$ to 
$\overline{Y} \otimes \overline{X} \otimes X \otimes Y$, by 
$I_{\overline{Y}} \otimes R_1 \otimes I_{Y}$, and a map
$j(\overline{R}_2)$ from  
$X \otimes \iota_B \otimes \overline{X} \simeq X \otimes \overline{X}$ to 
$X \otimes Y \otimes \overline{Y} \otimes \overline{X}$ by 
$I_{X} \otimes \overline{R}_2 \otimes I_{\overline{X}}$.  
\par
We also define a $C$--$C$ bimodule homomorphim 
$R \in {}_C\L(\iota_C,(\overline{Y}\otimes \overline{X}\otimes X\otimes
Y))$ by $R=i(R_1)\circ
R_2$, 
and $A$--$A$ bimodule homomorphism  
$\overline{R} \in {}_A\L(\iota_A, X\otimes Y \otimes \overline{Y}\otimes
\overline{X}_A)$ by 
$\overline{R}=j(\overline{R}_2)\circ \overline{R}_1$.  
Then we have 
$$
\align
 \overline{R}^* \otimes I_{X\otimes Y}\circ I_{X\otimes Y} \otimes R & = I_{X\otimes Y} \\
 R^* \otimes I_{\overline{Y}\otimes \overline{X}}\circ
I_{\overline{Y}\otimes \overline{X}} \otimes \overline{R} & =
I_{\overline{Y}\otimes \overline{X}} \endalign
$$

We check  the first relation:  
$$
\align
  & \overline{R}^* \otimes I_{X\otimes Y}\circ I_{X\otimes Y} \otimes R \\
  = &\overline{R}_1^* j(\overline{R}_2^*)\otimes I_{X \otimes
 Y}\circ I_{X \otimes Y} \otimes i(R_1)R_2 \\
  = & \overline{R}_1^* \otimes I_{X} \otimes I_{Y}\circ I_{X}\otimes
 R_1 \otimes I_{Y}\circ I_{X} \otimes \overline{R}_2^* \otimes
 I_{Y} \circ I_{X} \otimes I_{Y} \otimes R_2  \\
  = &((\overline{R}_1^* \otimes I_{X}\circ I_{X}\otimes R_1) \otimes
 I_{Y}) (I_{X} \otimes (\overline{R}_2^*\otimes I_{Y}\circ I_{Y}
 \otimes R_2)) \\
  =&  I_{X\otimes Y}.
\endalign
$$

To show  the second equation we proceed in a similar way
and we use the following computation:
for
 $x \otimes y \otimes \overline{y'} \otimes b\overline{x'} \in X \otimes Y
 \otimes \overline{Y} \otimes \overline{X}$ we have 
$$  
 (I_{X} \otimes \overline{R}_2^* \otimes I_{\overline{X}} \otimes I_{X}
 \otimes I_{Y})(I_{X} \otimes I_{Y} \otimes I_{\overline{Y}} \otimes R_1
 \otimes I_{Y})(x \otimes y \otimes \overline{y'} \otimes b\overline{x'})
$$
 $$=  x \otimes {}_B(y|y')R_1(b) \otimes \overline{x'} 
 =  x \otimes R_1({}_B(y|y')b) \otimes \overline{x'}$$
$$ =  x \otimes R_1({}_B(y|y'))\otimes b\overline{x'} 
 =  (I_{X} \otimes R_1 \otimes I_{Y})(I_{X}\otimes 
\overline{R}_2^* 
 \otimes I_{Y})(x\otimes y \otimes \overline{y'} \otimes
b\overline{x'}).$$
One can easily check that the following relations:
$${}_A(z|z'):=\overline{R}^*(\theta^r_{z,z'}\otimes
1_{\overline{Y}\otimes\overline{X}})\overline{R},$$
$$(z|z')_C:={R}^*(\theta^r_{\overline{z},\overline{z'}}\otimes
1_{X\otimes Y}){R}.$$
Here $z\in X\otimes Y\to \overline{z}\in \overline{Y}\otimes\overline{X}$
is the map taking the  simple tensor $x\otimes y$ to $\overline{y}\otimes
\overline{x}$. (This map is a well defined, $A$-$C$
antilinear and  bi-antiunitary 
with respect to the corresponding 
bi-Hilbertian structures.) 
For $z_1,\dots,z_n\in X\otimes_B Y$,
$$\|\sum
{}_A(z_i|z_i)\|=\|{\overline{R}}^*(\sum_1^n\theta^r_{z_i,z_i})\otimes
1_{\overline{Y}\otimes\overline{X}}\overline{R}\|\leq
\|\overline{R}\|^2\|\sum_1^n\theta^r_{z_i,z_i}\|,$$
therefore $X\otimes Y$ has finite right numerical index.
With a similar argument,  $X\otimes Y$ has finite left
numerical index. 
If $X$ and $Y$ have finite index, $R$ and $\overline{R}$ are intertwiners 
of the $C^*$--category ${}_\A\H_\A$, by part (2) of theorem 4.4, so 
$X\otimes_BY$ has finite index
by Theorem 4.14. 
\medskip

\heading 6. Examples
\endheading

 In this  section we discuss  
  examples of Hilbert $C^*$ bimodules 
of finite index with countable  bases.
\medskip

\noindent{\bf 6.1 Finite index bimodules generating Cuntz--Krieger
algebras} 
\bigskip

\noindent In the next example we  construct a Hilbert $C^*$--bimodule of
finite index which 
generates  a countably generated Cuntz--Krieger algebra, see 
\cite{KPRR} and \cite{KPW2}.  

 Let $\Sigma$ be a countable set, and let $G=(G(i,j))_{i,j\in \Sigma}$ be 
an infinite matrix with entries in $\{0, 1\}$.
We shall assume that no row and no column of $G$ is identically zero. 
We associate to the 
matrix $G$ the directed graph ${\Cal G}=(\Sigma,E,s,r)$,  where $\Sigma$ 
is 
the set of vertices and  
$E=\{(i,j)\in \Sigma \times \Sigma | G(i,j)=1 \}$
is the set of edges.  For an edge $\gamma =(i,j)\in E$,  the 
source $s(\gamma)$ is   $i$ and the range $r(\gamma)$ is  
$j$.  
We assume that ${\Cal G}$ is locally finite, that is, for any $j\in\Sigma$,
$\{i\in \Sigma | G(i,j)=1\}$ 
and, for any $i\in\Sigma$, $\{j\in \Sigma |G(i,j)=1\}$ are finite.
 
Let  $A=c_0(\Sigma)$ be the C${}^*$--algebra of the functions on 
$\Sigma$ vanishing at infinity and let $A_0=c_{00}(\Sigma)$ be the  
dense *-subalgebra of  functions with finite support.  
We denote by $P_j$ the projection in $A$ given by $P_j(i)=\delta_{ij}$. 
Since the set of edges  $E$ is a subset of $\Sigma \times \Sigma$, 
we may regard $E$ as a set-theoretic correspondence.  
The vector space $X_0=c_{00}(E)$ of 
the function on $E$ with finite support is an $A$--$A$ bimodule by
$$
 (a\cdot f\cdot b)(i,j) = a(i)f(i,j)b(j) 
 $$
for $a$, $b\in A$, $f \in X_0$ and $(i,j)\in E$.  
We define an $A$-valued inner product on $X_0$ by 
$$
 (f|g)_A(j) = \sum_{\{i|(i,j) \in E\}}\overline{f(i,j)}g(i,j) 
$$
for $f$, $g\in X_0$.   $X_0$ becomes in this way  a right pre-Hilbert 
$A$--module.  We denote by $X$ the completion of $X_0$.  The 
left $A$-action on $X_0$ can be extended to an action $\phi: A \to
\L(X_A)$ 
on $X$
by continuity.  Since $G$ is a row finite matrix,  
$\phi(a) \subset \K(X_A)$.  
Since no column of $G$ is zero, the range map $r$ is onto.  Thus  $X_A$ 
is 
full.
Let ${\Cal O}_X=C^*\{S_x|x \in X\}$ be the Pimsner algebra  \cite{Pim}
 generated by the bimodule $X$.  
For $\alpha \in E$, $S_{\delta_{\alpha}}$ will be denoted by 
$S_{\alpha}$.

Let $F$ be the edge matrix defined by $F(\alpha,\beta)=1$ if 
$r(\alpha)=s(\beta)$ and $F(\alpha,\beta)=0$ otherwise.  
Then the generators $\{S_{\alpha}|\alpha \in E\}$ satisfy 
$$
 S_{\alpha}^*S_{\alpha} = 
\sum_{\beta}F(\alpha,\beta)S_{\beta}S_{\beta}^*.
$$
For $i\in \Sigma$, we may define $S_i=\sum_{s(\alpha)=i}S_{\alpha}
\in {\Cal O}_X$, because no row of $G$ is zero and the source map
$s$ is onto. If $\beta=(i,j)\in E$, then $S_{\beta}=S_iP_j$.  
The $C^*$--algebra ${\Cal O}_X$ is also generated 
by $\{S_i|i\in \Sigma\}$ satisfying the relations 
$$
S_i^*S_i = \sum_j G(i,j)S_jS_j^* 
$$
The $C^*$--algebra ${\Cal O}_X$ coincides with the countably generated Cuntz-Krieger 
algebra  ${\Cal O}_G$.  
\medskip

   We shall introduce  an  $A$-valued left inner product  on $X$.  We 
need an additional datum.  Assume that we are given a nonnegative 
matrix 
$T = (T_{ij})_{ij}$ such that  $T_{ij} > 0$ if and only if 
$(i,j)$ is an edge, i.e. $G(i,j) = 1$.  We  call such a matrix $T$
 a {\it weight matrix} for the graph ${\Cal G}$.   
K. Yonetani suggested that the weight matrix $T$  gives  
an $A$-valued left inner product ${}_A(\ |\ )$ on $X_0$ by 
$$
{}_A(f|g)(i) = \sum_{j} T_{ij}f(i,j)\overline{g(i,j)}
$$
for $f$, $g\in X_0$.  Then we have two associated norms 
$$
{}_A\|f\| = \sqrt{\sup_i \sum_j T_{ij}|f(i,j)|^2}
$$
and 
$$
\|f\|_A = \sqrt{\sup_j \sum_i |f(i,j)|^2}
$$

\noindent{\bf 6.1 Definition}
     A weight  matrix $T$ for the graph ${\Cal G}$ is called {\it of 
finite index} 
if 
$$
c_1 := \sup_i \sum _j T_{ij} < \infty 
$$
and
$$
c_2 := \sup_j \sum _i \frac{1}{T_{ij}} < \infty.
$$
\medskip

\noindent{\bf 6.2 Example}
Let $\Sigma = \Bbb {N}$ and  $\ G(i,j) = 1\ $ if $|i-j| = 1$ and  $\ 
G(i,j) = 0\ $ if 
$|i-j| \not= 1$.  Consider a weight matrix $T$ defined as follows: 
$T_{12} = 1$.  $\ T_{ij} = 1/2\ $ if $|i-j| = 1$ and $(i,j) 
\not= (1,2)$.  
$\ T_{ij} = 0\ $ if $|i-j| \not= 1$.  Then $c_1 = 1$ and $c_2 = 4$. 
Thus  $T$ is of finite index.
\medskip

\noindent{\bf 6.3 Example}
Let $\Sigma = \Bbb {Z}$ and  $\ G(i,j) = 1\ $ if $|i-j| = 1$ and  $\ 
G(i,j) = 0\ $ if 
$|i-j| \not= 1$.  Consider a weight matrix $T$ defined by 
$\ T_{ij} = 1/2\ $ if $|i-j| = 1$ and 
$\ T_{ij} = 0\ $ if $|i-j| \not= 1$.  Then $c_1 = 1$ and $c_2 = 4$. 
Thus  $T$ is of finite index.
\medskip

\noindent{\bf 6.4 Example}
 The homogeneous tree $\text{Tree}(n)$ of degree $n$  is the tree where
all 
vertices have degree $n$.   For example $\text{Tree}(2)$ is the graph
above 
 with  $\Sigma = \Bbb {Z}$ and  $\ G(i,j) = 1\ $ if $|i-j| = 1$ and  $\ 
G(i,j) = 0\ $ if 
$|i-j| \not= 1$.  $\text{Tree}(4)$ is the Cayley graph of the free group 
$F_2$ with respect to the generators.  We define a weight matrix $T$ 
for $\text{Tree}(n)$ by associating the value $1/n$ with each edge.  Then
$c_1= 1$ and 
$c_2 = n^2$.  Thus $T$ is of finite index.
\medskip
  
\noindent{\bf 6.5 Example}
A tree has a weight matrix of finite index if and only if it has 
bounded degree.  In general a locally finite graph has a weight 
matrix of finite index if and only if both in-  and out-degrees 
are bounded.  In fact suppose that the in-degree is unbounded.  We may 
assume that $c_1 < \infty$.  For any edge $(i,j) \in E$, 
$0 <  T_{ij} \leq \sup _i \sum _j T_{ij} = c_1$.  Then we have 
$$
c_2 = \sup _j \sum _i \frac{1}{T_{ij}} \ge \sum _i \frac{1}{T_{ij}} 
\ge \sum_i \frac{1}{c_1}.
$$ 
Since the in-degree is unbounded, the last term goes to $\infty$.  Therefore 
$c_2 = \infty$.  The rest may be similarly shown.
\medskip

\noindent{\bf 6.6 Lemma} {\sl
 Let $T = (T_{ij})_{ij}$ be a weight matrix for a graph ${\Cal G}$.  
Then the following are equivalent: 
\roster
 \item $T$ is of finite index.
 \item The two norms $ {}_A\| \ \|$ and $\| \ \|_A $ on $X_0$ are 
 equivalent. 
\endroster}  
\medskip

\noindent{\it Proof}
(1)$\Rightarrow$ (2): Suppose that $T$ is of finite index.  Then 
for any $i$, 
$$
\sum_j T_{ij}|f(i,j)|^2 \leq (\sum_j T_{ij})(\sup_j |f(i,j)|^2) 
\leq c_1(\sup_j \sum_i|f(i,j)|^2)
 = c_1 \|f\|_A^2, 
$$
hence 
$$
{}_A\|f\| = \sqrt{\sup_i \sum_j T_{ij}|f(i,j)|^2} \leq 
 c_1^{1/2} \|f\|_A. 
$$
Setting $g(i,j) = \sqrt{T_{ij}}f(i,j)$, we have 
$$
{}_A\|f\| = \sqrt{\sup_i \sum_j |g(i,j)|^2}, 
\text{\qquad and  \qquad}
 \|f\|_A = \sqrt{\sup_j \sum_i \frac{1}{T_{ij}}|g(i,j)|^2}.
$$
Thus for any $j$,
$$
\sum_i \frac{1}{T_{ij}}|g(i,j)|^2
\leq (\sum_i \frac{1}{T_{ij}})(\sup _i|g(i,j)|^2)
\leq c_2\sup _i\sum _j|g(i,j)|^2 
= c_2 {}_A\|f\|^2.
$$
Hence
$$
 \frac{1}{c_2^{1/2}}\|f\|_A 
 = \frac{1}{c_2^{1/2}} \sqrt{\sup_j \sum_i \frac{1}{T_{ij}}|g(i,j)|^2} 
 \leq {}_A\|f\|.
$$
(2)$\Rightarrow$ (1): Suppose that $T$ is not of finite type.  
Then $c_1 = \infty$ or $c_2 = \infty$.  
If $c_1 = \infty$, then for any $M > 0$ there exist a positive 
integer $k$  such that $\sum_j T_{kj} \geq M$ .  
Let $f(i,j) = 1$ if $i = k$, $(k,j)$ is an edge,   
and $f(i,j) = 0$ otherwise.  
Then ${}_A\|f\|^2 \geq \sum_j T_{kj} \geq M$ and 
$\|f\|_A  = 1$ .  Thus the two norms are not equivalent. 
If $c_2 = \infty$,  then for any $M > 0$ there exist a positive 
integer $k$  such that $\sum_i \frac{1}{T_{ik}} \geq M$ .
Let $g(i,j) = \frac{1}{T_{ij}^{1/2}}$ if $j = k$, $(i,k)$ is an edge, 
and $g(i,j) = 0$ otherwise.  
Then ${}_A\|g\| = 1$ and 
$\|g\|_A^2  \geq M$ .  Thus the two norms are not equivalent. 
\medskip

If $T$ is of finite index, we can identify the two completions of $X_0$ 
with respect to  the two norms above defined.  We 
shall denote by $X$ its  completion.  The left and right  actions
of $A$ extend to injective  $^*$-homomorphisms 
$\phi : A \rightarrow \L(X_A)$ and 
$\psi : A \rightarrow \L ({}_AX)$.  
\par
We need the following inequalities which are easily verified:
Let $Y$ be a normed space.  Then for any $y_1,...,y_n \in Y$, 
positive numbers $\lambda _1, ... , \lambda _n, c$ with 
$\sum _i \lambda _i \leq c$, we have 
$$
\| \sum _i \lambda _i y_i \| \leq c\ \sup _i \|y_i\|.
$$
For any positive operators $A, B, C, D$ with $A \leq C$, $B \leq D$, 
we have 
$$
\| A^{1/2}B^{1/2} \| \leq  \| C^{1/2}D^{1/2} \|.
$$

\noindent{\bf 6.7 Theorem} {\sl
In the above situation, if a weighted matrix $T$ is of finite index, 
then $X$ is of finite index and 
$$
\|r-\text{Ind}[X] \| = c_1 := \sup_i \sum _j T_{ij}  \qquad
{\text and}
\qquad 
\|\ell-\text{Ind}[X] \| = c_2 := \sup_j \sum _i \frac{1}{T_{ij}}.
$$
More precisely, we have 
$$
r-\text{Ind}[X]  = ( \sum _j T_{ij})_i \in \ell ^{\infty}(\Sigma)  \qquad
{\text and}
\qquad 
\ell-\text{Ind}[X]  = ( \sum _i \frac{1}{T_{ij}})_j \in \ell^{\infty}(\Sigma) 
$$
}\medskip

\noindent{\it Proof}  First we  shall show  that $X$ is of finite right index.
Since the graph is locally finite,  we have 
$\phi(A) \subset \K(X_A)$ and 
$\psi(A) \subset \K ({}_AX)$.  In fact, we have 
$$
\phi(P_i) = \sum _{\{\beta | s(\beta) = i\}} 
\theta ^{r}_{\delta_{\beta}, \delta_{\beta}} \in \K(X_A)
$$
and 
$$
\psi(P_j) = \sum _{\{\gamma | r(\gamma) = j\}} 
\frac{1}{T_{\gamma}} \theta ^{\ell}_{\delta_{\gamma}, \delta_{\gamma}} 
\in \K({}_AX)  \ \ .
$$
   
We are left to show
the inequality described in part (2) of Prop. 2.7.
 (2) of Definition 2.3.  For any $f_1, ..., f_n 
\in X_0$ and  $g_1, ..., g_n \in X_0$, we shall show that 
$$
 \| \sum_{p=1}^n{}_A(f_p|g_p)\|
 \le  c_1 \| \sum_{p=1}^n \theta^r_{f_p,g_p}\|. 
$$  
Since
$$
\align
& \| \sum_{p=1}^n \theta^r_{f_p,g_p}\| 
= \|((f_p|f_q)_A)_{pq}^{1/2}((g_p|g_q)_A)_{pq}^{1/2}\| \\
& = \sup _j \|(\sum _i \overline{f_p(i,j)}f_q(i,j))_{pq}^{1/2}
(\sum _i \overline{g_p(i,j)}g_q(i,j))_{pq}^{1/2}\|, 
\endalign
$$
we have that 
$$
\align
 \| \sum_{p=1}^n{}_A(f_p|g_p)\| 
 & = \sup _i | \sum _j T_{ij}(\sum _{p=1}^n f_p(i,j)
 \overline{g_p(i,j)}) |  \\
   & \leq   \sup _i ((\sum _j T_{ij}) \sup _j | \sum _{p=1}^n f_p(i,j)
 \overline{g_p(i,j)}) |)  \\
   & \leq   (\sup _i (\sum _j T_{ij}))(\sup _i \sup _j 
 | \sum _{p=1}^n f_p(i,j) \overline{g_p(i,j)}) |)  \\
 & = c_1 \sup _i \sup _j 
 | \sum _{p=1}^n f_p(i,j) \overline{g_p(i,j)}) |  \\
 & = c_1 \sup _j \sup _i 
 \| ( \overline{f_p(i,j)}f_q(i,j))_{pq}^{1/2}
(\overline{g_p(i,j)}g_q(i,j))_{pq}^{1/2} \|  \\
& \leq  c_1 \sup _j 
\|(\sum _i \overline{f_p(i,j)}f_q(i,j))_{pq}^{1/2}
(\sum _i \overline{g_p(i,j)}g_q(i,j))_{pq}^{1/2}\| \\
& = c_1 \| \sum_{p=1}^n \theta^r_{f_p,g_p}\|.
\endalign
$$ 
 For any $f_1, ..., f_n \in X_0$, we shall show that 
$$
 \| \sum_{p=1}^n \theta^r_{f_p,f_p}\| 
 \le  c_2 \| \sum_{p=1}^n{}_A(f_p|f_p)\|.
$$
Put $g_p(i,j) :=\sqrt{T_{ij}}f_p(i,j)$.  
Then we have that 
$$
\align
&  \| \sum_{p=1}^n \theta^r_{f_p,f_p}\| 
= \|((f_p|f_q)_A)_{pq} \|  \\
& = \sup _j \|(\sum _i \overline{f_p(i,j)}f_q(i,j))_{pq} \| 
\endalign
$$
$$
\align
& = \sup _j \|(\sum _i \frac{1}{T_{ij}} 
\overline{g_p(i,j)}g_q(i,j))_{pq} \| \\
& \le \sup _j ((\sum _i \frac{1}{T_{ij}}) 
   \sup _i \|( \overline{g_p(i,j)}g_q(i,j))_{pq} \|) \\
& \le (\sup _j \sum _i \frac{1}{T_{ij}}) 
   \sup _j \sup _i \|( \overline{g_p(i,j)}g_q(i,j))_{pq} \| 
\endalign
$$
$$
\align
& = c_2 \sup _j \sup _i \|( \overline{g_p(i,j)}g_q(i,j))_{pq} \| \\
& = c_2 \sup _i \sup _j |\sum _{p=1}^n 
    g_p(i,j) \overline{g_p(i,j)}) |  \\
& \le c_2 \sup _i |\sum _{j=1}^n \sum _{p=1}^n 
  g_p(i,j) \overline{g_p(i,j)}) |  \\
& = c_2  \sup _i |\sum _{j=1}^n \sum _{p=1}^n 
  T_{ij}f_p(i,j) \overline{f_p(i,j)}) |  \\ 
& =  c_2 \| \sum_{p=1}^n{}_A(f_p|f_p)\|.  
\endalign 
$$ 
We shall next show that $X$ is of finite left index. 
We denote by $Y_0$ be the $A$--$A$ bimodule $c_{00}(E)$ with the 
following two-sided inner products:  For $\hat{f}$, $\hat{g} \in Y_0$,  
$$
{}_A\langle\hat{f}|\hat{g}\rangle(i) = 
  \sum_j \hat{f}(i,j) \overline{\hat{g}(i,j)}
$$
and 
$$
\langle \hat{f}|\hat{g}\rangle_A(j) = 
  \sum_i \frac{1}{T_{ij}} \overline{\hat{f}(i,j)} \hat{g}(i,j).
$$
We denote by $Y$ its completion.  For $f \in Y_0$, define 
$Uf \in Y_0$ by $(Uf)(i,j) = \sqrt{T_{ij}}f(i,j)$.  
Then we have ${}_A(f|g) = {}_A\langle Uf|Ug\rangle$  and 
$(f|g)_A = \langle Uf|Ug \rangle _A$.  The map $U$ extends to 
a surjective isometry $X \rightarrow Y$ with respect to the two-sided  
inner products.  Combining the fact with the preceding argument, 
we see that $X$ is of finite left index.  
Since $\{\delta_{(i,j)}\}_{(i,j)\in E}$ is a right basis for $X$, 
$$
r-\text{Ind}[X]  =  \sum _{(i,j)} {}_A(\delta_{(i,j)} | \delta_{(i,j)}) 
= ( \sum _j T_{ij})_i \in \ell ^{\infty}(\Sigma). 
$$
Since $\{\frac{1}{\sqrt{T_{ij}}}\delta_{(i,j)}\}_{(i,j)\in E}$
 is a left basis for $X$
the formula for $\ell-\text{Ind}[X]$ is obtained similarly. The rest is
clear.
\medskip

\noindent{\bf  6.2 Crossed products of 
Hilbert $C^*$--bimodules by locally compact groups}
\bigskip

\noindent In \cite{K}, the first-named author studied continuous crossed
products of 
Hilbert $C^*$--bimodules by locally compact groups. 
Let $B$ be a unital C${}^*$--algebra, and $A$ be a C${}^*$-subalgebra of
$B$ with the same unit. Let $E: B \rightarrow A$ be a conditional 
expectation of finite index in the sense of \cite{W}.  So there exists a
finite basis 
$\{u_1,u_2,\dots,u_n\}$ of $B$ such
that $x = \sum_{i=1}^nE(xu^*_i)u_i$ for any $x \in B$.  
Let $X = {}_AB_B$ be a 
$A$--$B$ bimodule with right $B$-valued inner product $(x|y)_B = x^*y$ 
and 
left $A$-valued inner product ${}_A(x|y) = E(xy^*)$.  Then  $X$ is a 
Hilbert 
$A-B$ bimodule of finite index.

 Let $G$ be a second countable locally compact group and $\alpha$  a
continuous homomorphism from $G$ to the automorphism group of $B$ such
that $E(\alpha_g(b)) = \alpha_g(E(b))$  for every $b \in B$ and
every $g \in G$.  It can be shown that $A \rtimes_{\alpha}G$ can be 
embedded as a  C${}^*$-algebra in $B \rtimes_{\alpha}G$ in a natural way,
and it can be shown that there exists a conditional expectation $\hat{E}$ from 
$B \rtimes_{\alpha}G$ to $A \rtimes_{\alpha}G$ which extends $E$.  
Put $Y=B \rtimes_{\alpha}G$.  We define a $A\rtimes_{\alpha}G$-
$B \rtimes_{\alpha}G$ bimodule structure on $Y$ and a  left 
inner product over $A\rtimes_{\alpha}G$ and a right  inner product over $B\rtimes_{\alpha}G$
in the obvious
way using $\tilde{E}$.  Then it can be shown that $Y$ is a countably
generated Hilbert C${}^*$--bimodule of finite index (see \cite{K}).
The left  and right indices of $Y$ are essentially the same as those of $X$.  
\medskip

\noindent{\bf 6.8 Correspondences} 
 \bigskip

Let $\Omega$ be a compact Hausdorff space.
Most Hilbert \cst--bimodules over the commutative \cst-algebra  
$A = C(\Omega)$ naturally arise from  set-theoretical 
correspondences (i.e. closed subsets $\C$ of $\Omega \times \Omega$)
similar to the case of commutative von Neumann algebras as in \cite{Co}. 
We say that  a pair $(\C,\mu)$ is  a (multiplicity free) topological 
correspondence on $\Omega$ if $\C$ is a (closed) subset of 
$\Omega \times \Omega$ and 
$\mu = (\mu^y)_{y \in \Omega}$ is a family of finite regular Borel 
measure on $\Omega$ satisfying the following conditions:
 \roster
\item (faithfulness) the support $supp \mu^y$ of the measure $\mu^y$ is 
the $y$-section $\C^y:=\{ x \in \Omega | (x,y)\in \C\}$, 
\item
(continuity) for any $f\in C(\C)$, the map $y\in \Omega \to 
\int_{\C^y} f(x,y) \,d\mu^y(x) \in \comp$ is continuous.
\endroster  
The vector space $X_0=C(\C)$ is an $A$--$A$ bimodule by
$$
  (a\cdot f \cdot b)(x,y) = a(x) f(x,y) b(y) 
$$
for $a$, $b \in A$, $f \in X_0$ and $(x,y) \in \C$.  
We define an $A$-valued inner product on $X_0$ by
$$
    (f|g)_A(y) = \int_{\C^y} \overline{f(x,y)}g(x,y)\,d\mu^y(x)
$$
for $f$, $g\in X_0$.  Faithfulness and continuity of 
$\mu$ imply that $X_0$ is a right pre-Hilbert $A$--module.  
We denote by $X$ the completion $X_0$.  The left $A$-action on $X_0$ 
can be extended to a *-homomorphism $\phi : A \to {\Cal L}_A(X_A)$.  
Thus we obtain a right Hilbert $A$--$A$ bimodule with right inner products
from the correspondence $(\C,\mu)$. See \cite{D} and \cite{KW1} 
for a more precise treatment.
 
We usually assume that for any $x\in \Omega $ 
 there exists $y \in \Omega$  with$(x,y) \in \C$.  This condition 
implies that left action $\phi$ is faithful.  
We also assume that for any $y \in \Omega$
 there exists $x \in \Omega$ with $(x,y) \in \C$. The condition 
shows that right inner product on $X$ is full. In fact let 
$\omega(y) = \mu ^y(\C^y)$.  Then $\omega \in A$ is invertible. 
For any $a \in A$, put $f(x,y) = a(y)$.  Then $(I|f)_A = a\omega$. 
Hence the right inner product is full. 
\medskip

\noindent{\bf 6.9 Example}
Let us  assume that projection maps 
$$r:(x,y) \in \C \mapsto x \in \Omega \qquad \text{and} \qquad  
s: (x,y) \in \C \mapsto y \in \Omega
$$
are  local homeomorphisms.  For any $y \in \Omega$, 
let  $\mu ^y$ be the counting measure on $\C^y$. 
Then $(\C,\mu)$ is  a topological 
correspondence on $\Omega$.    
We shall show $X$ has a finite basis.  In fact, 
since $\C$ is compact,
and by 
 our assumption, there exist a finite set 
$\{(x_1,y_1), \ldots  (x_n,y_n)\} \subset \C$ and 
open neighborhoods $U_k$ of $(x_k,y_k)$ for $k = 1,\ldots, n$ such 
that the restrictions of the projection maps $r$ and $s$ to $U_k$ are 
local homeomorphisms and $\C = \cup_{k=1}^n U_k $ is an open 
covering. Let $\{f_1,\ldots, f_n \} \subset C(\C)$ be a partition 
of unity for this  open covering.  Put $g_k = f_k^{1/2} \ge 0$.  
Then for any $(x_1,y), (x_2,y) \in \C$, we have 
$$
\sum _{k=1}^n g_k(x_1,y)\overline{g_k(x_2,y)} = \delta_{x_1,x_2}.
$$
Using these equalities, for any $h \in C(\C)$,  we have that 
$\sum_{k=1}^n g_k(g_k|h)_A = h$.  Thus $\{g_1, \ldots , g_n\}$ 
is a finite basis for $X$. 
 
As an example, put  $\Omega = [0,1]$.  Let  $h_1$ be a 
map on the interval $\Omega = [0,1/2]$ given by 
$h_1(x) = 2x$ for $x \in [0,1/2]$  and  $h_2$ be a 
map on the interval $\Omega = [1/2,1]$ given by 
$h_2(x) = 2x-1$ for $x \in [1/2,1]$.  
Let $\C$ be the union of the  graphs of $h_1$ and $h_2$.  
Then $A = C([0,1])$ and  the right inner product on $X_A =  C(\C)$ is 
given by 
$$
(f|g)_A(y) =  \overline{f(y/2,y)}g(y/2,y) + 
              \overline{f(y/2+1/2,y)}g(y/2+1/2,y)
$$
for $f$, $g\in X_A$ and $y \in [0,1]$.
Thus $X_A  \cong  A \oplus A$ as a
right Hilbert $C^*$--module.  The left action
$\phi : A \to {\Cal L}_A(X_A)  
\cong M_2(A) \cong C([0,1], M_2({\Bbb C}))$ is given by  
the diagonal matrices 
$$
(\phi(a))(x) = \text{diag}(a(x/2), a(x/2 +1/2)).
$$
The associated Pimsner algebra $\O_X$ is isomorphic to the  Cuntz algebra 
$\O_2$
and the fixed point algebra by the gauge action is isomorphic to  a UHF 
algebra 
$M_{2^{\infty}}$. 
The left inner product on $X_A$ is similarly given by 
$$
{}_A(f|g)(x) =
\cases
  {f(x,2x)}\overline{g(x,2x)} & 
    (\text{if $0 \leq x \leq 1/2$})   \\
  f(x, 2x-1)\overline{g(x,2x-1)} &
    (\text{if $1/2 \leq x \leq 1$}). 
\endcases
$$              
The right and the left norms on $X_A$ are given by 
$$
\|f\|_A = \sup _{0 \leq y \leq 1} \sqrt{|f(y/2,y)|^2 + 
              |f(y/2+1/2,y)|^2 }
$$
and 
$$
{}_A\|f\| = \max \{\sup _{0 \leq x  \leq 1/2} |f(x,2x)| ,
 \ \sup _{1/2 \leq x  \leq 1} |f(x,2x-1)| \},
$$
These two norms are equivalent:
$$
{}_A\|f\| \leq \|f\|_A \leq \sqrt{2}  {}_A\|f\|  \ .
$$            
Thus the $A$--$A$ bimodule $X_A$ is of finite type in the 
sense of \cite{KW1}.  This implies that   
 $X$ is of finite index as discussed in Example 2.34.
 
Let us replace $\Omega = [0,1]$ by the circle $\Omega = {\Bbb T}$ 
and consider, similarly, the   map $h$ on $ {\Bbb T}$ such that  $h(z) = 
z^2$ for 
$z \in {\Bbb T}$.  Let $\C$ be the graph of $h$ and consider the 
Hilbert $A$--$A$ bimodule $X$.  Then 
the associated Pimsner algebra $\O_X$ is isomorphic to the purely infinite 
simple 
$C^*$--algebra with $K$--theory: $K_0(\O_X) \cong {\Bbb Z} \cong K_1(\O_X)$. 
The fixed point algebra under the gauge action is isomorphic to the 
Bunce-Deddens algebra of type $2^{\infty}$.   One can also show that $X$ 
is of finite index. 
\medskip

\noindent{\bf 6.10 Example} In general $X$ fails to have a finite basis. 
A typical example is supplied  by a tent map $h$ on the unit 
interval $\Omega = [0,1]$. This is essentially the same example as the one
discussed in 2.35.  
The map $h$ is given by 
$$
h(x) = 2x \quad (\text{if} \quad 
  0\le x \le 1/2) \quad \text{and} \quad   
h(x) = -2x+2 \quad (\text{if} \quad 
1/2 \le x \le 1).
$$  
 Let $\C$ be the graph of $h$.        
For any $y \in [0,1)$, 
let  $\mu ^y$ be a counting measure on $\C^y$, that is, 
 $\mu ^y = \delta_{y/2} + \delta_{1-y/2} $.  For $y=1$, let 
 $\mu ^1 = 2\delta_{1/2}$.  
Then $(\C,\mu)$ is  a correspondence on $[0,1]$.  
Let $X$ be the right Hilbert $A = C([0,1])$--module 
obtained from the correspondence.  The right Hilbert $A$--module 
 $X$ does not have a finite basis.  This fact is easily seen through a realization of 
it 
as an orbifold construction as follows:  Let $B = C([0,2])$ and let
$\gamma$ be a homeomorphism on $[0,2]$ with period two such that 
$\gamma(x) = 2-x$ for $x \in [0,2]$.  Then $\gamma$ induces an 
automorphism
 $\alpha$ on $B$ such that $(\alpha(f))(x) = f(2-x)$.  The fixed point 
algebra 
 $B^{\alpha}$ is isomorphic to $A = C([0,1])$.  We note that the action 
 $\alpha$ is free except for $x = 1$.  Consider the condiitonal expectation 
 $E: B \rightarrow A$ defined by $E(f) = (f + \alpha (f))/2$.  Let 
 $Y = B_A$ be a right Hilbert $A$--module  
 given by $(f|g)_A = E(f^*g)$.  Then  $X$ and $Y$ are 
isomorphic 
 as right Hilbert $A$--modules: there is a   unitary  induced 
 by the map $\phi : \C \rightarrow [0,2]$, 
$$
\phi((x,2x)) = 2x \quad  \text{if} \quad 
 0 \le x \le 1/2
$$
$$
 \phi((x,-2x +2)) = 2x \quad \text{if}  
\quad 1/2 \le x \le 1.
$$  
By a result 
 in Prop. 2.8.2 in \cite{W}, $Y$ does not have a
finite basis.  
 
 We shall construct a countable basis for $Y$ explicitly.
 Define $r_n \in C([0,1])$ by 
$$
r_n(x) = 1 \quad \text{if} \quad  0 \le x \le 1-1/n \quad 
\text{and} \quad  
 r_n(x) = -n(x-1) \quad \text{if} \quad  1-1/n \le x \le 1. 
$$ 
Put $v_1 = r_1$ and 
 $v_n = (r_n - r_{n-1})^{1/2} \in C([0,1])$ for $n \ge 2$.  
 Then 
$$
\sum_{i=1}^{\infty} |v_i(x)|^2 = \lim_{n\to\infty}r_n(x) = 1
\qquad \text{for} \quad 
 x  \not= 1
$$
 and $v_i(1) = 0$.  Define $u_0 = 1$ and $u_i \in C([0,2])$ by 
 $$
u_i(x)  = v_i(x) \quad (\text{if} \quad 0 \le x \le 1) 
\quad \text{and}  \quad u_i(x) = -v_i(2-x) \quad  (\text{if} 
\quad  1 \le x \le 2)
$$
for \ $ i = 1,2,\ldots $, and  $u_i(1) = 0$ for $i \not= 0$. 
 Then for $x \not = 1$, we have $\sum_{i=0}^{\infty} |u_i(x)|^2 = 2$ and 
 $\sum_{i=0}^{\infty} u_i(x)u_i(2-x) = 0$.  For $x = 1$, we have 
 $\sum_{i=0}^{\infty}|u_i(1)|^2 = 1$.  We claim  that 
$(u_i)_{i=0,1,\ldots}$ 
 is a basis for $Y$.   For any $f \in Y = C([0,2])$, we shall show that
$$ 
 \lim_{F:\text{finite}} \sum_{i \in F} u_i(u_i|f)_A = f .
$$ 
Let $F$ be a 
finite subset of $\{0,1,2,\ldots\}$. Put $S_F = \sum_{i\in F} u_i(u_i|f)
_A$.  
Since $\|(g|g)_A\| = \|E(g^*g)\| \le \|g\|_{\infty}^2$ for $g \in B$, it 
suffices to show that $\lim_{F:\text{finite}} \| f - S_F\|_{\infty}$ = 0, 
where $\|g\|_{\infty}$ is the sup norm. We may assume that $0 \in F$. 
For $0 \le x \le 1$, we have 
$$ 
\align
  & f(x) - S_F(x) = f(x) - \{ \frac{1}{2} \sum_{i \in F}|u_i(x)|^2f(x)
     + \frac{1}{2} \sum_{i \in F}u_i(x)u_i(2-x)f(2-x) \} \\
 = & f(x) - \{ \frac{1}{2}(1 + \sum_{i \in F\setminus \{0\}} 
v_i^2(x))f(x)
          + \frac{1}{2}(1 - \sum_{i \in F\setminus \{0\}} 
v_i^2(x))f(2-x) \}  \\
 = & \frac{1}{2}(1 - \sum_{i \in F\setminus \{0\}} v_i^2(x))(f(x) - 
f(2-x)). 
\endalign
$$      
For any $\varepsilon$ there exist $\delta > 0$ such that for any  
$x$ satisfying   $1- \delta \le x \le 1$ , 
 we have $|f(x) - f(2-x)| \le \varepsilon$ . 
Since $|1 - \sum_{i \in F\setminus \{0\}} v_i^2(x)| \le 1$, we have 
$$
|f(x) - S_F(x)| \le \varepsilon \qquad \text{for} \qquad x \in 
[1-\delta, 1].
$$
On the other hand,   $\sum_{i=1}^{\infty} v_i(x)^2 = 1$ 
  uniformly
on $[0, 1-\delta]$.  We also have  $|f(x) - f(2-x)| \le 2 
\|f\|_{\infty}$.
Therefore there exists a finite set $F_0$ such that for any finite 
subset 
$F \supset F_0$ 
and for for any $x \in [0,1]$, 
$$
|f(x) - S_F(x) | \le \varepsilon.
$$
Similar arguments work for $x \in [1,2]$.  
Therefore  $(u_i)_{i=0,1,\ldots}$
 is a countable basis for $Y$  
However, the range $\phi (A)$ of the left action is not included in  
$\K(X_B)$, because the idenitity $\phi(I_A) = I_{X_B}$ is not in 
$\K(X_B)$.  Therefore  the bimodule $X$ is not of finite right index.
\medskip

\heading
References
\endheading

\item[B] B. Blackadar, 
 K-theory for operator algebras. Mathematical Sciences Research
Institute Publications, 5 Cambridge University Press, Cambridge, 1998. 
 
\item[BDH] M. Ballet, Y. Denizeau, J.-F. Havet,
{\it Indice d'une esperance conditionelle,}  Compos. Math.
66(1988), 199-236.

\item[BGR] L. G. Brown, P. Green, M.A. Rieffel, {\it Stable isomorphism 
and strong Morita equivalence of $C^*$--algebras,} Pacific J. Math. 71
(1977), 349--363.

\item[D] V.Deaconu,
{\it Generalized solenoids and \cst--algebras,}
Pacific J.Math. 190 (1999), 247-260

\item[Di] J. Dixmier, Les $C^*$--alg\`ebres et leurs repr\'esentations.
Gauthier-Villars \'Editeur, Paris 1969.

\item[DPZ] S. Doplicher, C. Pinzari, R. Zuccante, {\it The
$C^*$--algebra of a Hilbert bimodule,} Boll. Unione Mat. Ital. Sez. B, 1
(1998), 263--281.

\item[DR1] S. Doplicher, J.E. Roberts, 
{\it Endomorphisms of $C^*$--algebras,
crossed products and duality for compact groups,} Ann. Math. 130 (1989) 
75--119.

\item[F] J. M. G. Fell, {\it The structure of algebras of operator 
fields,} Acta Math. 106 (1961) 233--280.

\item[DR2] S. Doplicher, J.E. Roberts, 
{\it A new duality theory for compact groups}
Inventiones Math. 98 (1989) 157--218.
 
\item[FK]  M. Frank and E. Kirchberg, 
{\it On conditinal expectations of finite index,} 
 J. Operator Theory 40(1998), 87-111.

\item[FL1] M. Frank and D. Larson, 
{\it A module frame concept for Hilbert $C^*$--modules, }
Contemporary Mathematics 247(1999), 207-233.

\item[FL2] M. Frank and D. Larson, 
{\it Frames in  Hilbert $C^*$--modules and $C^*$--algebras }
preprint.

\item[H] R. Haag, Local quantum physics. Fields, particles, algebras.
Second edition. Springer-Verlag, Berlin, 1996.

\item[I] M. Izumi, {\it Inclusions of simple $C^*$--algebras,}
J. Reine Angew. 547 (2002), 97--138.

\item[J] V.F.R.Jones,  {\it Index for subfactors,}
Inv. Math.  72 (1983)  1-25.

\item[K] T. Kajiwara,
{\it Continuous crossed products of Hilbert \cst--bimodules,} 
Int. J. Math. 1 (2000),  969--981.  

\item[KW] E. Kirchberg, S. Wassermann, {\it Operations on continuous
bundles of $C^*$--algebras,} Math. Ann. 303 (1995), 677--697.

\item[KPRR] A. Kumjian, D. Pask, I. Raeburn, J. Renault, 
{\it Graphs, groupoids and Cuntz-Krieger algebras,} 
J. Funct. Anal. 184(1998), 161-174.
 
\item[KPW1] T. Kajiwara, C. Pinzari and Y. Watatani
{\it Ideal structure and simplicity of the C${}^*$--algebras generated 
by Hilbert bimodules,} J. Funct. Anal. 159(1998), 295-322.

\item[KPW2] T. Kajiwara, C. Pinzari and Y. Watatani
{\it Hilbert \cst--bimodules and countably generated Cuntz-Krieger 
algebras,} 
 J. Operator Theory 45(2001), 3-18.

\item[KW1] T. Kajiwara and Y. Watatani
{\it Jones index theory by Hilbert C${}^*$--bimodules 
and K-theory,} Trans. Amer. Math. Soc. 352(2000), 3429-3472.  

\item[KW2] T. Kajiwara and W. Watatani
{\it Crossed product of Hilbert \cst--bimodules by countable discrete 
groups,} Proc. Amer. Math. Soc. 126 (1998), 841-851.

\item[KW3]  T.Kajiwara and Y. Watatani
{\it Hilbert \cst--bimodules and continuous Cuntz-Krieger algebras,} 
 J. Math. Soc. Japan. 54 (2002), 35-59.

\item[Ka] G.G. Kasparov, {\it Equivariant KK-theory and the Novikov
Conjecture}, Invent. Math. 91 (1988), 147--201.

\item[L]  E.C. Lance,
{\it Hilbert \cst--modules : A toolkit for operator algebraists, } 
London Math. Soc. Lecture Note Ser. 210, Cambridge U.P., 1995

\item[LR] R. Longo, J.E. Roberts, 
{\it A theory of dimension,} K-theory 11(1997),103-159. 

\item[P] G.K. Pedersen, $C^*$--algebras and their automorphism groups,
Academic Press
1979.

\item[Pim] Pimsner, M., {\it A class of
$C^*$--algebras generalizing both
Cuntz--Krieger algebras and crossed products by ${\Bbb Z}$,} in
Voiculescu, D. (ed.)
{\it Free probability theory\/}, AMS, 1997, 189-212.

\item[PiPo] M. Pimsner, S. Popa,
{\it Entropy and index for
subfactors,}  Ann. Sci. Ec. Norm. Sup. 19 (1985),  57--106.

\item[R1] M. Rieffel,
{\it Induced representations of $C^*$--algebras,}
Adv. Math. 13 (1974), 176-257.

\item[Y]  S. Yamagami
{\it Tensor categories for operator algebraists,} preprint.

\item[V] J.-M. Vallin, {\it $C^*$--alg\`ebres de Hopf et
$C^*$--alg\`ebres de Kac,} Proc. London Math. Soc. 50 (1985), 131--174.

\item[W]  Y. Watatani,
{\it Index  for \cst--subalgebras,}
Memoir Amer. Math. Soc. 424 (1990).

\item[Wo] S. L. Woronowicz, {\it Tannaka--Krein duality for compact matrix
pseudogroups. Twisted $SU(N)$ groups,} Invent. Math. 93 (1988), 35--76.

\end{document}